\documentclass[11pt]{amsart}
\usepackage{amsmath,amsfonts,amssymb,graphicx,amsthm}
\usepackage{epstopdf}
\usepackage{verbatim,bbm,latexsym,indentfirst,cases,mathtools}
\usepackage{enumerate}
\usepackage{epsfig,color}
\usepackage[noadjust]{cite}
\usepackage{nomencl}
\makenomenclature
\newcommand{\goto}{\rightarrow}
\newcommand{\ind}{\mathbf{1}}
\newcommand{\e}{\mathbb{E}}

\def\refl{reflection}
\def\expP{exponential \ \Par \ \refl} 
\def\Par{Parisian}

\def\for{\forall} \def\mc{\mathcal}

\def\mT{{\mathcal T}}  
  
 \def\mCB F{${\mathcal CB F}$}

\def\Rui{\Psi} \def\sRui{\ovl{\Rui}} 
 \def\bf{\bfseries}  \def\rm{\textrm} \def\it{\itshape}
\def\mG{\mathcal G}  
\newcommand{\be}{\begin{equation}}
\newcommand{\ee}{\end{equation}}
\newcommand{\ba}{\begin{array}}
\newcommand{\ea}{\end{array}}
\newcommand{\baa}{\left[\begin{array}}
\newcommand{\eaa}{\end{array}\right]}

\def\BEN{\begin{enumerate}}  \def\BI{\begin{itemize}}
\def\EEN{\end{enumerate}}   \def\EI{\end{itemize}}
\def\BM{Brownian motion }
\def\beq{\begin{eqnarray}}\def\eeq{\end{eqnarray}}
\def\bea{\begin{eqnarray*}}
\def\eea{\end{eqnarray*}}
\def\le{\left} \def\ri{\right}
\def\T{\widetilde}  \def\H{\widehat}
 \def\WH{Wiener-Hopf }

\def\I{\infty}  \def\a{\alpha}

\def\g{\gamma}   \def\de{\delta}
\def\z{\zeta}  \def\th{\theta}
\def\e{\epsilon} \def\k{\kappa} \def\l{\lambda}
\def\lm{\Pi}

\def\x{\xi }
  \def\r{\lambda} \def\s{\sigma}
\def\t{\tau}   \def\f{\phi}   
\def\c{c} \def\w{\omega}  \def\q{q} \def\qu{\quad} \def\D{\Delta}
\def\F{\Phi}  \def\L{L} \def\U{U}

\def\la{\label} \def\fr{\frac} \def\im{\item}
\def\Tl{\tilde{\lambda}} \def\Tq{\tilde{\q}}
 \def\corr{corresponding }

\def\z {\zeta}

\def\mbb{\mathbb} \def\mbf{\mathbf} 
\def\mc{\mathcal}   

\newtheorem{Thm}{Theorem}[section]
\newtheorem{Con}{Conjecture}
\def\beCo{\begin{Con}
  }
  \def\eeCo{\end{Con}}
  \def\C{{\mathbb C}} \def\E{{\mathbb E}}  \def\R{{\mathbb R}}
\def\beT{\begin{Thm}
  }
  \def\eeT{\end{Thm}} \newtheorem{Qu}{Problem}
\def\beQ{\begin{Qu}}
  \def\eeQ{\end{Qu}}
\newtheorem{Lem}{Lemma}
\newtheorem{Exa}{Example}
\newtheorem{Cor}{Corollary}[section]
\newtheorem{Def}{Definition}
\newtheorem{Rem}{Remark}[section]
\newtheorem{Exe}{Exercice}
\newtheorem{Ass}{Assumption}
\def\beXe{\begin{Exe}} \def\eeXe{\end{Exe}}
\def\eeD{\end{Def}} \def\beD{\begin{Def}}
\def\beXa{\begin{Exa}} \def\eeXa{\end{Exa}}
\def\beR{\begin{Rem}} \def\eeR{\end{Rem}}
\def\beL{\begin{Lem}} \def\eeL{\end{Lem}}
\newtheorem{Pro}{Proposition}
\def\beP{\begin{Pro}} \def\eeP{\end{Pro}}
\def\beC{\begin{Cor}}
\def\eeC{\end{Cor}}
\def\bc{\begin{cases}
  }
\def\ec{\end{cases}}
\def\BEN{\begin{enumerate}} \def\EEN{\end{enumerate}}

\long\def\symbolfootnote[#1]#2{
\begingroup
\def\thefootnote{\fnsymbol{footnote}}\footnote[#1]{#2}
\endgroup}
\def\fn{\symbolfootnote}
\def\ssec{\subsection} \def\sec{\section}
 \def\Y{Y}
\def\H{\widehat} 
\def\no{\nonumber} \def\bs{\bigskip} 
\def\bep{\begin{pmatrix}}
  \def\eep{\end{pmatrix}}

\newcommand{\Z}[2]{\ensuremath{Z^{(#1)}(#2)}}

  \def\Eq{\Leftrightarrow} \def\OU{Ornstein-Uhlenbeck}
  \def\CL{Cram\'er-Lundberg }
  \def\PK{Pollaczek-Khinchine }
   \def\valf{value function }
   \def\Frt{Furthermore, }
  \def\und{\underline} \newcommand{\md}{\mathrm d}
   \def\Mp{More precisely, }
    \def\P{{\mathbb P}}

 \def\Lra{\Longrightarrow}
\def\SA{Sparre-Andersen }
\def\nne{nonnegative } \def\nny{nonnegativity }

  \def\ci{capital injections} 
\def\Lm{L\'evy measure}

\def\1{\mathbbm{1}}
\def\Tl{\T{\lambda} }
 \def\AY{Azema-Yor }
\def\ki{\overleftarrow{\k}}
 \def\expoc{exponential claims}

\def\mbb{\mathbb}

\def\sn{spectrally negative }
\def\sp{spectrally positive }
\def\LE{Laplace exponent } \def\LEs{Laplace exponents }\def\LT{Laplace transform } \def\itf{it follows that }
\def\LTs{Laplace transforms }

\def\kil{\mathbf e}

  \def\PH{phase-type }

\newcommand{\red}{\textcolor[rgb]{1.00,0.00,0.00}}

\def\lev{L\'evy } \def\LTs{Laplace transforms }   
  
   \def\eq{equation }
 \def\Xt{(X_t)_{t \geq 0}}
 \def\thr{therefore }
  
\def\deF{de Finetti }  \def\MAP{spectrally negative Markov additive processes }
\def\GS{Gerber-Shiu } \def\vars{variables} \def\ssr{self similar}
\newcommand{\diff}{{\rm d}} \def\fp{first passage }
\def\ovl{\overline}
  \def\Ea{\E^{[a}}  \def\Eb{\E^{b]}}
 \def\Pb{P^{b]}}   \def\Ez{\E^{[0}} 
    \def\Ezb{\E^{[0,b]}} \def\Pzb{P^{[0,b]}}
    \def\Eazb{\E^{|0,b]}}

     \def\ola{\overleftarrow}
   \def\MAP{Markov additive processes } \def\vt{\vartheta}
    \def\mD{{\mathcal D}}
    \def\zrb{z^{b[}} \def\vfe{ } \def\asy{asymptotic }
   \def\wrb{w^{b[}}  \newcommand{\prf}{\noindent{\bf Proof:\ }}  \def\tbz{T^{[0}_b} \def\tba{T^{[a}_{b}}
    \def\ith{it holds that } \def\mD{{\mathcal D}}
    \def\mZ{{\mathcal Z}}  \def\mW{{\mathcal W}} 
\def\t{\tau} \def\ta{T_{a,-}}  \def\tb{T_{b,+}}
\def\tz{T_{0}}

 \def\iffa{}  
  \def\Fqp{\Phi_\q'} \def\fe{for example }
 \def\TW{W^{(\Fq)}_0} \def\HW{\H{W}^{(\Fq)}}
 
\newcommand{\Fq}{\Phi_q} \newcommand{ \Fqr }{\Phi(\q+\r)}
 \def\strs{strong Markov processes}
 \def\cP{compound Poison model } \def\cmy{complete monotonicity } \def\Bf{barrier function } \def\Dis{\mathfrak D}  \def\equ{equation } \def\fund{fundamental } \def\satg{satisfying }
  \def\rw{random walk }   \def\rv{random variable } \def\ind{independent } \def\upm{up to a multiplicative constant} 
 \def\str{strong Markov property} \def\ts{two-sided } \def\wk{well-known}
   \def\Mar{Markov } \def\eqr{\eqref}  
  \def\funs{functions }   \def\Thr{Therefore, }
 \def\prob{problem }  \def\pro{probability }\def\sur{survival } \def\sats{satisfies } \def\bco{boundary condition }  \def\resp{respectively} \def\how{however}
 \def\boun{boundary } \def\con{condition }  \def\Fr{Furthermore, } \def\frt{furthermore }  \def\nonh{nonhomogeneous }\def\exc{excursion} \def\tzb{T_{0}^{b]}}  \newcommand{\cd}{(\cdot)} 
 \def\exp{exponential }  \def\cP{compound Poisson } \def\LZ{Lokka-Zervos alternative }
\def\pros{probabilities } \def\fund{fundamental } \def\gen{generalized }
 \def\ol{\overline}  \def\gene{generalization } \def\SLG{Shreve, Lehoczky and Gaver }
\def\Y{Y} \def\tzh{T_{\{0\}}} \def\RHS{right hand side } \def\td{\tau}
 \def\OU{Ornstein-Uhlenbeck }  
\def\procs{processes} \def\proc{process} \newcommand{\sme}[2]{\frac{W_\q(#1)}{W_\q(#2)}}
\newcommand{\nsme}[2]{Z_\q(#1)-\frac{W_\q(#1)}{W_\q(#2)}Z_\q(#2)}

\newcommand{\pd}[2]{\frac{\partial #1}{\partial #2}}

\def\nun{ten } \def\Nun{Ten } \def\run{twelve } 
\newcommand{\figu}[3]{
\begin{figure}[!h]
\begin{center}
{\includegraphics[width=13 cm, height=8 cm]{#1}}
\end{center}

\vspace{-0.2cm}
\caption{\hspace{0.25cm}#2\label{f:#1}}
\end{figure}
}
\newcommand{\ninseps}[3]{
\begin{figure}[!h]
\begin{center}
{\includegraphics[width=13 cm, height=8 cm]{#1}}
\end{center}

\vspace{-0.2cm}
\caption{\hspace{0.25cm}#2\label{f:#1}}
\end{figure}
}
\begin{document}
\title[The $W,Z$ scale functions kit for  spectrally negative L\'evy
processes]
{The $W,Z$ scale functions kit for first passage problems of  spectrally negative L\'evy
processes, and applications to control problems}

\author{Florin Avram}  \thanks{Laboratoire de Math\'ematiques Appliqu\'ees, Universit\'e de Pau,  France, florin.avram@univ-Pau.fr}
 \author{Danijel Grahovac} \thanks{Department of Mathematics,
University of Osijek, Croatia, dgrahova@mathos.hr}
\author{Ceren Vardar-Acar} \thanks{ Department of Statistics, Middle East Technical University, Ankara, Turkey, cvardar@metu.edu.tr}

\date{\today}

\begin{abstract}
In the last years there appeared  a great variety of identities for first passage  problems of spectrally negative  L\'evy  processes, which can all be expressed in terms of  two ``$q$-harmonic functions" (or scale functions) $W$
and $Z$. The reason behind that is that there are two ways of exiting an interval, and thus two fundamental
  ``two-sided exit'' problems from an interval  (TSE). Since many other problems can be reduced to TSE, researchers developed in the last years
 a kit of   formulas expressed  in terms of the ``$W,Z$ alphabet''.
It is important to note -- as is currently being shown -- that these identities apply equally to other spectrally negative Markov 
processes, where however the $W,Z$ functions are typically much harder to compute.
 We collect below our favorite recipes from the  \lev ``$W,Z$ kit", drawing from various applications in   mathematical finance,  risk, queueing,  and inventory/storage theory.
 A small sample of applications concerning  extensions of the classic de Finetti
dividend  problem is offered.
An interesting  use of the  kit is for recognizing
relationships between  problems involving  behaviors apparently unrelated at first sight (like reflection, absorption, etc). Another is expressing results in a standardized form, improving thus the possibility to   check when a formula is already known.

\end{abstract}
\maketitle
\tableofcontents


\sec{Introduction}
From our  biased point of view, the $W,Z$ scale functions kit is a new set of clothes for the classic first passage  theory used in  risk, queueing, mathematical finance and related fields, which was developed over the last $40$ years. A recent explosion of new contributions to this topic, notably to processes with Parisian ruin and reflection -- see Section \ref{s:Par}, and the extension to spectrally negative \Mar processes -- see Section \ref{s:dd}, suggested  the utility
of offering a new review. We attempted to pack in our ``cookbook'' a  possibly overwhelming quantity of  results; the best way for the reader to get an idea of what's to be found here might be to have first a quick look at the List of notations Section \ref{s:n}.

In this section we introduce the \CL risk \proc, we define first passage times and some main quantities of interest for  the control
and optimization of risk processes.

{\bf Origins}. The origins of our  field lie in the ruin \prob for  the {\bf Cram\'{e}r-Lundberg} or {\bf compound Poisson}  risk model \cite{lundberg1903approximerad,AA}
\begin{align} \la{CL}
X_t = x  -\Big(\sum_{i=1}^{N^{(\l)}_t} C_i- c t\Big).
\end{align}
 Here $c$ is the premium rate, $C_i, i=1,2,...$ are i.i.d.~{nonnegative} jumps with distribution $F(dz)$,  arriving after independent
exponentially distributed times with mean $1/\l$, and $N^{(\l)}$ denotes the associated Poisson process counting the arrivals.
Note that the process in parenthesis,  called ``cumulative loss'',  is used also to model the workload process of the M/G/1 queue.

{\bf First passage theory}  concerns the first passage times above and below, and the hitting  time of a level $b$. For any process $\Xt$, these are defined by
\begin{equation}\la{fpt}
\begin{aligned}
\tb &= \tb^{X}=\inf\{t\geq 0: X_t> b\},\\
T_{b,-}&=T_{b,-}^X=\inf\{t\geq 0: X_t< b\},\\
T_{\{b\}}&=T_{\{b\}}^X =\inf\{t\geq 0: X_t=b\},
\end{aligned}
\end{equation}
with $\inf\emptyset=+\infty$. The upper script $X$ will be typically  omitted, as well as the signs $+,-$, when they are clear from the context.

 First passage  times are important in the control of reserves/risk
 processes. The rough idea is that when below low levels $a$,  reserves processes should be replenished at some cost, and when above high levels $b$, they should be partly invested to yield income -- see for example \cite{AA} and, for  most recent work, papers like \cite{APP,IP,AIpower,AIjoint}, etc.

 The  first   quantity to be studied  historically was the eventual ruin probability
  $$\Rui(x)=P_x[\tz < \I]$$
  for the  {Cram\'{e}r-Lundberg/compound Poisson}  risk model \cite{lundberg1903approximerad,AA}.
  Subsequently, first passage (or exit) problems were studied   in  mathematical finance (barrier options, American options -- see for example \cite{Kyp}),  in
risk \cite{AA}, queueing \cite{APQ}  storage theory \cite{brockwell1982storage,yamazaki2016inventory}, in mathematical biology \cite{ricciardi1999outline}, and in many other applications.
  The typical approach  for a long while consisted in taking Laplace transform of the associated Kolmogorov integro-differential equation involving the generator operator.


  In recent years it became clear
  that most first passage  problems for spectrally negative or spectrally positive L\'evy processes may be reduced to the solution of the two fundamental
  ``two-sided exit'' problems from an interval (TSE), upwards or downwards.
   At their turn, these  can be ergonomically expressed in terms of two scale functions/$q$-harmonic functions $W_{\q}(x), Z_{\q}(x,\th)$. In the case of spectrally negative processes, one ends up with  the following equations:\fn[4]{The first equation generalizes the famous ``gambler's winning" formula for the symmetric random walk $\sRui_0^{b}(x,a)  ={\E}_{x}\left[  \1_{\left\{
\tb<\ta \right\}  }\right]=\fr{x-a}{ b-a}$.}
\begin{align} \la{sRui}
\sRui_{\q}^{b}(x,a)  &  :={\E}_{x}\left[  e^{-\q\tb}\1_{\left\{
\tb<\ta \right\}  }\right]=\fr{W_{\q}(x-a)}{W_{\q}(b-a)}, \ q\geq0, \  a\leq x\leq b, \\
\Rui^{b}_{\q,\th}(x,a)  &  :={\E}_{x}\left[  e^{-\q\ta+\th \left( X_{\ta}-a
\right)}\1_{\left\{  \ta<\tb \right\}  }\right]=Z_{\q}(x-a,\th)-W_{\q}(x-a)\dfrac{Z_{\q}(b-a,\th)}{W_{\q}(b-a)}, \ \th \geq 0.\la{Rui}
\end{align}
We  will  call $\sRui_{\q}^{b}(x,a), \Rui_{\q}^{b}(x,a)$   (killed) {\bf survival} and {\bf ruin} first passage  probabilities, respectively. When $a = 0$, it will be omitted, to simplify the notation.

 \beR Note that the first quotient decomposition above holds true by the absence of positive jumps and by
the strong Markov property, and that this defines $W_\q$ \upm. The second  relation is  equivalent to \eqr{B} below which defines $Z_{\q}$ \upm  ~(see \cite[Thm12]{IP} and Remark \ref{r:proof} below). For many other results in this vein,  see
\cite{Suprun,Ber97,Ber,AKP,Kyp,zhou2007exit,IP,KKR,APP15,LP}, and many other papers listed in the more detailed but still too  succinct  {chronology} in Section \ref{s:ch} below.
\eeR

 \beR \la{r:PK} {\bf The relation between $W(x)$ and $\sRui(x)$}.
When $\q= 0$,  the scale function $W(x):=W_0(x)$ is related to  the eventual ruin $\Rui(x)=P_x[ \tz < \I]$ and ultimate survival probabilities $\sRui(x)=P_x[ \tz = \I]$,  via
\beq \la{rp} \Rui(x)=1-\sRui(x)=1-\k'(0_+) W(x).\eeq
Here $\k$ is the \LE of $X$ given below in \eqr{LE} and the \LT of $W(x)$ is $\H W(s)=\fr{1}{\k(s)}$. Note that above and throughout
the paper we will assume that $\k'(0_+)$ exists, which renders formulas simpler (and is typically satisfied in applications).
\eqr{rp} is related to  the famous Pollaczek-–Khinchine formula for the \LT of the survival function of a spectrally negative L\'evy process
\beq \la{PKs}
 \H {\sRui} (s):=\int_{0}^\I e^{- s x} \;   \sRui(x) \; d x= \fr{\k'(0_+) }{\k(s)}.\eeq

 The scale function $W(x)$ provides an alternative  characterization of
a spectrally negative L\'evy process, which may replace the classic Laplace exponent $\k(s)$. 

\eeR

  \beR
    The eventual ruin and \sur probability    have made the object of numerous numerical studies, for example by inversion of Pad\'e approximations of $\fr{1}{\k(s)}$ \cite{avram2011moments,AAK,avram2018ruin} -- see \cite{AA} for other methods and references. Furthermore, it is easy to adapt  numerical studies of $W$  to yield $ W_{\q}$, by the so called Esscher transform (replacing $\k(s)$ by $\k(s+\q)-\k(\q)$) -- see Remark \ref{r:tr}. Note that
once $W_\q$ and $Z_{\q}$ are computed, we have obtained also the answer to many other problems, thus  removing the need for \LT inversion.  Hence, a cookbook of $W_\q, Z_{\q}$  formulas provides an alternative to the classic Markovian analytic approach.

\eeR

Before continuing, we note that  the last decade has witnessed also very interesting research on {\bf last passage times} -- see \fe \ \cite{baurdoux2009last,paroissin2015first,li2017last,
cai2018occupation}. Since we had  to stop at some point,  these will not be covered in our review.

   {\bf Control of dividends and capital injections}. The next impetus came from  control problems in risk theory which concern versions  of  $X_t$ which are
 reflected/constrained/regulated at first passage times (below or above):
 \begin{eqnarray} \la{Xd}&& X_t^{[a}=X_t  + \; \L_t, \qu  X^{b]}_t=X_t   -  \U_t. \la{ref}\end{eqnarray}
Here, \beq \la{RR} \no && \L_t=\L^{[a }_t=-(\und X_t-a)_-, \quad \und X_t=\inf_{0 \leq s\leq t} X_{s},\\ \no && \U_t=\U^{b]}_t=\le(\ovl X_t -b\ri)_+, \quad \ovl X_t := \sup_{0 \leq s\leq t} X_{s},\eeq
  are  the minimal  ``Skorokhod regulators" constraining $X_t$ to be bigger than $a$, and    smaller than $b$, respectively, and we use the notation $x_+=\max(x,0)$ and $x_-=\min(x,0)$.

  One  problem of historical interest is the de Finetti problem of expected total discounted  dividends until the ruin time $T_{0,-}$, in the presence of a constant (reflecting ) dividend barrier -- see \eqr{div}. Interestingly, its solution
\begin{equation*} V^{b]}(x)=\E_x\le[\int_{[0,T_{0,-}]}e^{-\q t}d   \U_t\ri]=\frac{ W_{\q}(  x)}{W_{\q}^ {\prime}(b)}
\end{equation*}
looks very similar to \eqr{sRui}. Intuitively, this is due to the fact that the two problems differ only
in what happens at the  boundary $b$ (reflection versus absorbtion), which is translated \resp \ into the boundary
conditions $\sRui^b_{\q}(b)=1$, $(V^{b]})'(b)=1$ -- see Remark \ref{divb}. In fact, this is the heart of the $W,Z$ theory: problems which differ only via their boundary behavior have similar answers --
 see  Section \ref{s:10} for further examples.

{\bf Drawdowns and drawups}.  Applications require often the study
of the running maximum and of the process reflected at its maximum/drawdown \beq \la{ddp} Y_t=\ovl{ X}_t-X_t, \quad \ovl X_t= \sup_{0 \leq s\leq t} X_s,\eeq
or that of the running infimum and of the process reflected from below/drawup \beq \la{dup}\und Y_t= X_t-\und X_t, \quad \und X_t=\inf_{0 \leq s\leq t} X_s.
\eeq

 The first passage times of the reflected processes, called
 drawdown/regret time and  drawup time, respectively, are defined for $d>0$ by
\begin{equation} \label{dd}
\begin{aligned}
\t_{d}&:=\inf\{t\geq 0: \ovl X(t)- X(t) > d \},\\
\und \t_{d}&:=\inf\{t\geq 0:  X(t)- \und X(t) > d \}.
\end{aligned}
\end{equation} \iffa
Such times  
turn out to be optimal in several stopping problems, in  statistics \cite{page1954continuous}    in mathematical finance/risk theory  (in problems involving  dividends at a fixed barrier or capital injections) --  see for example \cite{taylor1975stopped,Leh,shepp1993russian,AKP,MP,shiryaevthou,Carr,LLZ17,LLZ17b,baurdoux2017future}
  and in queueing theory (for example when studying idle times until a buffer reaches capacity) -- see for example
  \cite{debicki2012infimum,debicki2015queues}.

{\bf Capital injections/bail-outs}. A second important problem is that of the expected capital injections necessary to maintain a process positive, before reaching an upper barrier; this involves two reflecting  boundaries. Since  problems  with double reflection live on  finite intervals,
    the possibility
  to solve them by  Laplace transforms   seems lost at first; however,  their solutions are also expressible in terms of the fundamental scale functions $W_\q, Z_{\q}$.

For example, the joint \LT of the total regulation/capital injections into a spectrally negative      process \eqref{Xd} reflected at $a$ and of the first up-crossing of a level $b$ is \cite[Thm.~2]{IP}
\begin{equation} \la{B}
\sRui_{\q,\th}^b(x,[a):=\Ea_x \left[e^{-\q T_b^{[a} - \th \L_{T_b^{[a}}}\right]=\Ea_x \left[e^{  - \th \L_{T_b^{[a}}};T_b^{[a} <e_\q\right]
= \dfrac{Z_{\q}(x-a,\th)}{Z_{\q}(b-a,\th)},
\end{equation}
where  $\Ea_x$  denotes the expectation for the process reflected at $a$,  $T_b^{[a}$ denotes the corresponding hitting time \eqr{Tra}, and $e_\q$ denotes an independent exponential random variable of rate $\q$.  This factorization is essentially a direct consequence of  the \str. In our view, it is maybe the most important first passage law -- see Theorem \ref{l:refl}.

{\bf Joint behavior of the process and its drawdown}. The third act in the development of risk theory  was the consideration of the joint behavior of the process and its historical maxima or minima, or, equivalently, of the process and its drawdowns or drawups.
It turns out that this study, just like the previous   problems,  may be reduced to finding the $W_{\q},Z_{\q}$ functions--
see for example Theorem \ref{l:joiDB}.

{\bf Contents}.  We start with a brief review  of   L\'{e}vy processes in Section \ref{s:lev}. 
    Section \ref{s:W0} presents  the function $W_\q$ which is the pillar of this field, and includes three  remarkable results
  in which it appears. Section \ref{s:Z0}  introduces the  $Z_\q$ scale function, and  Section \ref{s:Z} introduces a two variables extension $Z_{\q}(x,\th)$ of  $Z_{\q}(x)$.

We turn next to the extensive and expanding body of knowledge concerning spectrally negative \lev processes. Our $W,Z$ ``cookbook"     collects a list of  some of our  favorite recipes. They come from many recent papers, like  \cite{AKP,Pispot,APP,Iva,IP,
 Ivapot,AIpower,APP15,APY,AIjoint} and other papers cited below, and we apologize for any omission.
 Section \ref{s:10} 
 alone lists
         \nun of the most important first passage laws, dubbed theorems, an eleventh ``meta theorem" including the ``Poissonian/Parisian version" of most of the first \nun theorems is presented in Section \ref{s:Par}, and  other \run results spread throughout the paper are called propositions (this  partition was adopted for the same reasons we organize files in folders).

Section \ref{s:nonh} reviews some $W,Z$ formulas for  smooth \GS functions.
Here the smooth \GS function $Z_{\q}(x,\th)$ which corresponds to the overshoot penalty $e^{\th X_{\tz}}$ is replaced   by a function $G_{w}(x)$ corresponding to an arbitrary penalty function $w( X_{\tz})$.

Section \ref{s:Par}  reviews $W,Z$ formulas for    processes with Poisonian/Parisian observations, and for the more general Omega processes.   The idea, which emerged naturally in the last decade
    in the context of financial modeling, is that ``transgressing boundaries" may pass unnoticed,
    with or without purpose, if observations are not continuous. This gives rise
    to ``soft boundaries", in addition to the traditional reflecting and absorbing ``hard boundaries" from the physics world; it seems \thr an important
      development  in the theory of Markov \procs. This topic is excellently presented
      in the article \cite{AIZ}, but we go beyond that. Quite surprisigly, despite the fact that the methods of proof are different,    we have showed in \cite{AZ}  that several  of the Parisian formulas  coincide with the classic ones, in terms  of two new scale functions (which generalize the classic ones). The same phenomenon was observed in \cite{APY} for processes with Parisian observations within a finite buffer, below which absolute ruin occurs. It is still not understood
      why  the classic and (buffered) Parisian laws look identical, once the appropriate scale functions have been identified. Let us note that due to its theoretical and  applied implications, this topic constitutes  an active field of research, with many open problems, some of which are listed below.

  To illustrate the potential   applicability of $W,Z$ formulas, we have  included in Section \ref{s:div} an important application: the  optimization of dividends, under several objectives. We have chosen this application partly since it is a fundamental brick in the budding
  discipline of risk networks \cite{AM15,AM16,AZ}. We also chose this  to emphasize that the famous and still not completely understood de Finetti optimization problem \cite{deF,gerber1969,AM05,APP,Schmidli,
  Loef08, azcue2014stochastic,APP15} is just one of a family
  of similar  optimization problems which can be tackled via the scale  function methodology, some of which may  be  more  tractable than the original.
    Section \ref{s:ex}   illustrates  the results on examples like Brownian motion \ref{s:BM} and exponential claims  \ref{s:exp}, and Subsection \ref{s:num} illustrates the numerical optimization of dividends for the Azcue-Muler example \cite{AM05}.

    Section \ref{s:dd} reports on  recent results on drawdown  problems. The motivation is
    to explore the idea
that in risk control  (and optimal consumption/harvesting problems) it may be profitable to base decisions both on the position of the underlying process  and on its distance from previous suprema. This suggests basing  decisions on  Azema-Yor/generalized drawdown/trailing stop times,  which involve certain admissible functions of the position and supremum.  This framework provides a natural unification of  drawdown and classic first passage times.

   It was discovered in this context that $W,Z$ formulas continue to hold for spectrally negative \Mar processes \cite{LLZ17b}.  The only difference is that in equations like \eqr{sRui} and \eqr{Rui},  $W_{\q}(x-a)$, $Z_{\q}(x-a,\th)$ must be replaced by  \funs  with one more variable $W_{\q}(x,a)$, $Z_{\q}(x,a,\th)$. Unfortunately, the computation of these
   scale \funs is currently understood in only one particular case outside \lev \procs \and diffusions:
   that of  Ornstein-Uhlenbeck
with phase-type jumps, treated in Jacobsen-Jensen \cite{JacJen}. However, we believe that other diffusions with phase-type jumps will  be treated in the future via variations of this  approach. For that reason,  we decided to present  the last Section \ref{s:dd} in the context  of spectrally negative \Mar processes (note though that this is mostly uncharted territory).

        The paper ends with a short chronology in Section   \ref{s:ch}, and a summary of notations and asymptotic formulas
    in Sections   \ref{s:ch}, \ref{s:n}, \ref{s:asy}.

We hope that our   compilation  may be of help as  a quick introduction to more detailed treatments like \cite{Ber,doney07,Kyp,KKR,KypGS} and also as a cookbook for computing  quantities of interest in applications like risk theory, mathematical finance, inventory and queueing theory, reliability, etc.
We will be forced to make appeal to the literature for many proofs, but some of the most useful methods of attack will be included.


\sec{A glimpse  of L\'evy processes \la{s:lev}}
  A
{\bf L\'evy process} \cite{Ber,Kyp} $X=X_t\in \R, t\geq 0$ may be characterized by its  L\'evy-Khinchine/Laplace exponent/symbol $\k(\th)$,   defined by
\beq \E_0 \left[e^{\th X_t }\right]=e^{t \; \k(\th)},  \la{LE}\eeq
where $\th \in \mD \subset \C$, and $\mD$  includes at least the imaginary axis.

\lev processes and their reflections (drawdowns and drawups)  satisfy a duality result \cite[Prop.~VI.3]{Ber}, \cite[Lem.~3.5]{Kyp}:
\beL \la{p:du} For each fixed $t > 0$, the pairs $(\ovl X_t , X_t -\ovl X_t )$ and $(X_t -\und X_t,\und X_t)$
have the same distribution under $P_0$. \eeL
\beR This result is behind the well-known  duality  between queueing and risk theories, which are concerned with reflected and absorbed processes, respectively. For example,
applying it when $t \to \I$ to the negative of the \CL process $-X$, when $\k'(0_+) >0$, yields the well-known identity between the stationary law of the  M/G/1 workload process and the infimum $\und X_{\I}$ of the \CL risk process -- see \cite{asmussen1992computational,APQ}, and  see  \cite{PDR,boxma2011threshold} for  further applications. \eeR

\beR \la{r:WH} The reflected processes of a L\'evy process are Markov processes 
\cite[Prop.~VI.1]{Ber}; therefore,  nice results on them and first drawdown /drawup passage times  are to be expected.
\eeR

\lev \procs  \  satisfy the \wk \ Wiener Hopf factorization \cite[Prop.~VI.5]{Ber}, a short version of which is:
\beL \la{p:WH}  Let $G_t := \sup \{ 0\leq  s \leq t :
X_s = \ovl X_t\}$ be the last time the process $X$ equals its supremum before or at
time $t$ ($ t-G_t $ is therefore the duration of the last  drawdown  at time $t$). For any \ind \exp \rv $e_\q$ with rate $\q > 0$, the pairs $(\ovl X(e_\q), \ovl G(e_\q))$ and $(X(e_\q)-\ovl X(e_\q), e_\q- \ovl G(e_\q))$ are independent under $P_0$.
\eeL

\ssec{The spectrally negative L\'evy risk model
\la{s:levsn}}

     From now on, $X_t, t\geq 0$ will denote a spectrally negative L\'evy process.
It is natural in applications to restrict to the case   when the \LE  has a L\'evy-Khinchine decomposition of the form
 \beq \la{k} \k(\th) = \frac{\s^2}{2}\th^2 + p \th + \int_{(0,\I)}[e^{- \th y} - 1 + \th y ]\lm(d y), \;  \th \geq 0,\eeq
 with a L\'evy measure $\lm$ of $-X$   satisfying
 \be \int_{(0,\I)} (y\wedge y^2) \lm(dy) < \I \la{lmp} \ee
  (and  $\lm(-\I,0)=0$)\fn[3]{Note that even though $X$ has only negative jumps, for
convenience we work with the L\'evy measure of $-X$.}.
 This implies that
   the growth (or profit) rate \sats $$\E_0 [X(1)]=p=\k'(0_+) \neq \I, $$
    a reasonable assumption in risk theory.

 This assumption excludes L\'evy measures like $\Pi(dx)=x^{-2} dx$ and $\a$-stable processes with $\a \in (0,1)$, but  it allows $\a$-stable processes with $\a \in [1, 2)$  (the \Lm \   is allowed to have infinite mean, as long as $\int_1^\I y \lm_Z(dy)< \I$).

\beR
 $X_t$ is a Markovian process  with infinitesimal generator $\mG$, which acts
on $h\in C^2_0(\R_+)$ by \cite[Thm.~31.5]{Sato}
\beq \label{genrisk}
\mG h(x) = \frac{\s^2}{2}h''(x) + p h'(x) + \int_{(0,\I)}[h(x-y) - h(x) + yh'(x)]\lm(d  y)
\eeq
(where we used \eqr{lmp}).
Incidentally, this may be formally written as $\mG=\k(D)$, where $D$ denotes the differentiation operator.

\eeR
If furthermore  the   jumps of the process have a finite mean
$\int_0^\I y \lm(dy)< \I$ (but not necessarily finite mass, which allows including interesting examples like the Gamma process \cite{DGGamma}), we may rewrite  \eqref{k} as \bea 
\k(\th) = \frac{\s^2}{2}\th^2 + c \th + \int_{(0,\I)}[e^{- \th y} - 1 ]\lm(d y), \qu  \th \geq 0, \qu  c := p+ \int_{(0,\I)} z \Pi (d  z),\eea
which reflects a decomposition into a Brownian motion with parameters $(c,\s)$ and the negative of a
subordinator.
We will call this the {\bf Brownian perturbed finite mean subordinator  risk model}.

A  further particular case to  bear in mind is that
when  the L\'evy measure has finite mass $\lm(0,\I)= \l<\I$. We may write then $\lm(dz)= \l F(dz)$, and rewrite
the process  and its symbol as
  \beq \la{pCL} X_t = x + \s B_t+ c t -\sum_{i=1}^{N^{(\lambda)}_t}  C_i, \qu \k(\th) = \fr{\s^2 \th^2} 2 + c \th + \lambda \H f_C(\th)- \lambda , \eeq
  where $B_t$ is the Wiener process, $C_i, i=1,2,...$ are i.i.d.~{nonnegative} jumps with distribution $F(dz)$,  arriving after
exponentially distributed times with mean $1/\lambda$, and $\H f_C$ denotes the Laplace transform of $C_i$.
This is the {\bf Brownian perturbed compound Poisson}   risk model \cite{DG}. If furthermore $X_t$ has paths of bounded variation, which happens  if and only if
$
\sigma = 0$, we obtain the classic  Cram\'{e}r-Lundberg  risk model \eqref{CL}. The simplicity of this case comes from the fact that its down-ladder times are discrete, which made it a natural favorite in risk theory.

Finally, let us mention the so-called ``\PK" processes which satisfy a generalization of the \PK formula \cite{DG}. The most general version due to \cite{HPV} is obtained by putting together a negative subordinator \satg $\int_{(0,\I)} (y\wedge 1) \lm(dy) < \I$ and an independent  \sn  zero mean  {\bf perturbation}  \satg \eqr{lmp}. The advantage of this class comes from the fact that its {\bf jump down-ladder} times are discrete.

{\bf State dependent \lev \procs}. Nowadays there is also considerable interest in ``\lev processes with state dependent
coefficients".  For example Albrecher and Cani studied the \CL process   with  affine dividends $ X_t= x+ \int_0^t (p- k X_s) \md s - \sum_0^{N_t^{(\l)}} C_i $ \cite{albrecher2017risk}, and  \cite{CPRY} studied a  more general  ``\lev driven Langevin model" $ \md X_t=p (X_t) \md t - \md S_t,  $ where $S_t$ is a \sp \lev \proc.

\sec{The scale function $W_\q$ and its logarithmic derivative $\nu_\q$ \la{s:W0}}

\ssec{Introduction}

 {\bf First passage results  for spectrally negative \lev processes} are remarkably  simpler than in the general case.
 Here everything reduces finally
 to the   determination of  the  ``scale functions'' $ W_{\q}(x):\R_+ \to [0, \infty), \q \geq 0$ defined on the positive half-line  by the Laplace transform \eqr{WLT}, and extended  to be $0$  on $\R_-$.
\beq \label{WLT}
&&\int_0^\infty  \mathrm{e}^{-s x}  W_{\q}(x) d x = \frac {1} {\k(s)-\q}, \qu \for s > \Phi(\q),
 \eeq
where $\Fq$ is the largest \nne root of the \CL equation
  \beq  \Phi(\q):= \sup \{ s \geq 0: \k(s) - \q= 0\}, \qu \q \geq 0. \label{Fq} \eeq
The scale function $ W_{\q}(x)$ is continuous and
increasing on $[0,\I)$ \cite{Bingham}, \cite[Thm.~VII.8]{Ber}, \cite[Thm.~8.1]{Kyp}.

Applying optional stopping at $T_{x,+}$ to the Wald martingale $e^{\Fq X_{t} - \q t} $  yields the fundamental identity
\beq \E_a \left[ e^{- \q T_{x,+}} \right] = e^{-(x-a) \Fq}=P_a [\ovl X(e_\q) >x], \la{e:exp}\eeq
where $e_\q$ is an independent exponential random variable with parameter $\q$ (thus, $T_{x,+}, x \geq 0$ is a subordinator, with
 Laplace exponent $\Fq$   \cite[Thm.~VII.1]{Ber}).

\beR  In the case of general L\'evy processes,  solving first passage problems rests on the Wiener-Hopf factorization of the  {Laplace exponent} with  killing $\k(s)-\q$ \cite[Prop.~VI.5]{Ber} (for meromorphic exponents, this means the identification and separation of the positive and negative roots, see \cite[Sec.~6.5.4]{Kyp} for details).\fn[4]{For a  proof using the Kella-Whitt martingale, see \cite[Thm.~4.8]{Kyp}.}
 The factorization simplifies considerably for  \lev processes which jump in only one direction (as is the case  in  queueing and risk theory), since then one part of the factorization involves only the root $\F(q)$ defined in \eqr{Fq}. Typically, this renders the factorization unnecessary, with most things expressable in terms of the pair of functions $\k,\F$.

  For example, in the spectrally negative case,
  the moment generating function of the drawdown $Y_{e_\q}$ at an exponential time $e_\q$, equal to that of $- \und X_{e_\q}$, satisfies  \cite[Thm.~4.8]{Kyp} \beq \la{gPK} \E_0[e^{-s Y_{e_\q}}]=\fr{ {s}-\F(\q)}{ {\k(s)}-\q} \fr{\q}{\Fq}.\eeq
  {When $\q \to 0$, this becomes the \PK formula}
$$ \E_0[e^{-s Y_\I}]=\fr{\k'(0_+) s}{\k(s)}, $$
which made some authors call \eqr{gPK} the generalized \PK formula.

Another case in which the factorization is easy to compute is that of   two-sided \PH jumps -- see for example \ \cite{asmussen2004russian}.

 \eeR

{\bf The smooth two-sided exit problem}. The most fundamental first passage problem  is the classic {\bf  gambler's winning problem} \cite{gerber1972games}, \cite[Thm.~3]{Suprun}, \cite[(6)]{Ber97}. This is an extension of \eqref{e:exp}, in which one kills the process upon reaching  a lower barrier $a$ which may be taken w.l.o.g. to be $0$.
 \beP \la{p:twos} For any $b > 0$ and $x \in [0, b]$\fn[4]{Note that \eqref{twos} may be obtained by stopping the martingale $ W_{\q}(X_t)$ at $\tb$.},
\beq \label{twos}  &&
		\sRui_{\q}^{b}(x)=\E_x \left[ e^{-\q \tb} \1_{\left\{ \tb <  \tz \right\}}\right]  =\frac { W_{\q}(x)}  { W_{\q}(b)}
 :=e^{-\int_x^b \nu_\q(s) ds}.   		
\eeq
Analytically, $\nu_\q(s)$ is the ``logarithmic derivative of $W_{\q}$ from the right'' \cite[(8.26)]{Kyp}, \beq \la{nu} \nu_\q(s)=\fr{ W_{\q}'(s+)}{ W_{\q}(s)},\eeq
and the ``from the right" will be omitted below since we assume $ W_{\q} \in C^1(0,\infty)$.\fn[5]{Since \eqref{twos} is the Laplace transform of the density of $\tb$, with absorbtion at
$\tz$, a  Laplace inversion will recover the corresponding density.}

\eeP

\beR \la{r:exrate} {\bf Two probabilistic interpretations of  $\nu_q$}. We are trying to avoid as much as possible in our review  the use of excursion theory. However, in preparation for the very important
problem of dividends paid under a constant barrier policy, we will make an exception, and present a ``homemade'' version  of excursion theory, explained in this remark and in Section \ref{s:dd}.
\BEN
\im It has been noted in \cite{albrecher2009tax}  that the last equality in \eqref{twos} may be  interpreted as the probability that no arrival has occurred between times $x$ and $b$, for a nonhomogeneous Poisson process of rate $\nu_\q(s)$.

 This checks with the probabilistic definition of $\nu_q(s)$ provided by \exc \ theory:
 $$\nu_\q(x):= n[\ovl \e >x, s(\e) \leq \kil_\q],$$ where  $n(d \e)$ is the characteristic measure of the Poisson process  of downward excursions $\e$ from a running maximum,  $\ovl \e$ denotes  the height of a downward excursion,
 $s(\e)$ denotes the starting time of an \exc, and $\kil_\q$ is an \ind \exp \rv of rate $\q$ -- see for example  \cite{Ber}, \cite[(12)]{Doney}.

\im  We would prefer to avoid \exc \ theory in our cookbook; however, the concept of \exc \ is too fundamental to be avoided.  We proceed \thr with a ``homemade" version of \exc \ theory for spectrally negative processes,  based on excising the negative excursions of $X_t$.

It has been noted in \cite{AACI,ALL} that differentiating  the last equality in \eqref{twos}  yields
\beq  \fr{d}{ds} \sRui_{\q}^{b}(s)-\nu_\q(s) \sRui_{\q}^{b}(s)=0, \quad \sRui_{\q}^{b}(b)=1. \la{genexc}\eeq
One may recognize here the Kolmogorov equation for the probability that a deterministic process $\T X(s)=s$ starting at $0$,  and also  killed at rate $\nu_\q(s)$ either when a negative \exc \ larger that $s$ occurs, or when an exponential clock of  rate $\q$ ticks,  reaches  $b$ before being killed. ``It turns out"
that $\T X(s)$ may be obtained by taking the running maximum value $s$ as time parameter,  and by
 excising the negative excursions of $ X(t)$ which are larger than $s$.  This interpretation is fundamental, and holds  for spectrally negative Markov processes as well -- see the last section \ref{s:dd}, in particular Remark \ref{r:snM}. $\T X(s)$ will be called from now on ``excised ladder process".

 Note that the quotation marks in ``it turns out" above and below
 mean that the statement can be left as an exercise for the \CL \proc, but needs in general careful treatment, which is beyond the scope of our cookbook.
\EEN

Summarizing this discussion, we retain that $\nu_\q(s)$ represents the    rate of the exponentially distributed period of time the process spends at an upward creeping moment (when $\ovl X_t =X_t$), before a downward excursion bigger than $s$ occurs, and before an exponential clock of  rate $\q$ ticks \cite{Kyp}.

This interpretation of $\nu_\q(s)$   is especially important in the de Finetti problem \eqr{div}, where
we will exploit the fact that the expected dividends $v_{\q}(b)$ paid at a fixed barrier $b$ when starting from $b$ equal the expected discounted time until killing of $\T X$.  This yields finally the simple relation
  \beq
    v_{\q}(b):=  \E_b \left[ \int_{0}^{T_{0,-} } e^{-q  t}dU_t \right]=\nu_{\q}(b)^{-1} \la{v}.\eeq
  This relation can be extended to spectrally negative \Mar processes  with \gen drawdown  \eqr{VLv}.
\eeR

{\bf  The 
smoothness of $W_q$}. 
Regarding the smoothness of the scale function, \ith \  $ W_{\q} \in C^1(0,\infty)$ iff
   the \lev measure has no atoms,  or  $X$ is of unbounded variation.  If  a Gaussian component is present ($\sigma > 0$), then \frt \ $ W_{\q} \in C^2(0,\infty)$.  See \cite{chan2011smoothness,doring2011densities} for further  results on  smoothness, and \cite{Loef08}
  for the case of completely monotone \lev measures\fn[6]{This paper shows that
if the \lev measure has a completely monotone density, $W_{\Fq} \in C^\infty (0,\infty)$, and
$W_{\Fq}'$ is also completely monotone.}. Below,  we will always assume that $ W_{\q}\cd$ is smooth enough
to satisfy the equation $\mG (W_q) (x) = q W_q(x)$ in the classical sense.

{\bf The behavior in the neighborhood of zero of $ W_{\q}$} can be obtained from the behavior of its \LT \eqref{WLT}
at $\I$ \cite[Lem.~ 4.3-4.4]{Kyprianou_Surya}, \cite[Lem.~3.2-3.3]{KKR}:
\begin{equation} \label{asymptotics_zero}
\begin{aligned}
 W_{\q} (0) &= \lim_{s \to \I} \fr{s}{\k(s) -\q}= \begin{cases} \frac 1 {\c}, & \text{if $X$ is of bounded variation/\CL }\\
0, & \text{if $X$ is of unbounded variation }
\end{cases},\\
 W_{\q}' (0_+) &= \lim_{s \to \I}s\left( \fr{s}{\k(s) -\q}- W_{\q} (0) \right)= \begin{cases}  \frac {\q + \Pi(0,\infty)} {\c^2}, & \text{if $X$ is of bounded variation}\\
\frac 2 {\sigma^2}, & \text{if } \sigma > 0, \\
\infty, & \text{if } \sigma = 0 \; \text{and} \; \Pi(0,\infty) = \infty
\end{cases}.
\end{aligned}
\end{equation}
Following the same approach, we may recursively compute  $  W_{\q}''(0)$, etc (these Taylor coefficients may be used in Pad\'e approximations, see \cite{avram2018ruin}). We find, when  the jump distribution has a density  $f$, that
\begin{equation} \label{e:secder}
\begin{aligned}
 W_{\q}'' (0_+) &=\lim_{s \to \I} s \left( s \left(\fr{s}{\k(s) -\q}- W_{\q} (0)\right) -  W_{\q}' (0_+) \right)\\
&= \begin{cases}   \fr 1 c \Big( (\fr {\lambda+ \q }  c )^2 - \fr {\lambda}  c f(0) \Big), & \text{if $X$ is of bounded variation}  \\
-c(\frac 2 {\sigma^2})^2, & \text{if } \sigma > 0
\end{cases},
\end{aligned}
\end{equation}
where the notation for the compound Poisson case is as in \eqref{CL}.
This equation is important in establishing the nonnegativity of the optimal dividends barrier -- see Example \ref{ex:expb}.

We offer now  as appetizer a strikingly beautiful recent application of the scale function due to \cite[(14)]{grahovac2017densities} to the calculation of the {\bf maximal severity
of ruin} \cite{picard1994some} -- see also \cite[Prop XII.2.15]{AA} for the compound Poisson case.
\beP \la{p:Grah} Let $$\eta:=T_{0}=\inf\{t > T_{0,-} :X_t  =  0\}$$
denote  the  {hitting}  time of $0$ (``recovery after  ruin'') -- see also \eqr{fhit}.

The cumulative distribution function of the maximal severity
of ruin $-\und {X}_{\eta}$ (i.e.~the absolute value of the infimum of the process before ``recovery after ruin")
  is  given by
\beq P_x [- \und {X}_{\eta} < u, \; T_{0,-} <\I ]=\fr{W(x+u)-W(x)}{W(u)}. \eeq

\eeP

{\bf Proof:}
By requiring that the first passage time precedes reaching $-u$ and by using
the gambler's winning identity \eqref{twos} one obtains that
\begin{equation}\label{recov:proof:eq}
P_x [- \und {X}_{\eta} < u,\; T_{0,-} <\I ]=\int_0^u P_x [-
{X}_{T_{0,-}} \in dy,\; T_{0,-} <\I ] \fr{W(u-y)}{W(u)}.
\end{equation}
On the other hand, by considering the
event of reaching $0$, but never reaching $-u$ at all we get
\begin{equation*}
\Rui(x)-\Rui(x+u) = \int_0^u P_x [-  {X}_{T_{0,-}} \in dy,\; T_{0,-} <\I
] \sRui(u-y),
\end{equation*}
and by using \eqref{rp} and \eqref{recov:proof:eq} it follows that
\begin{equation*}
W(x+u)-W(x)= \int_0^u P_x [-  {X}_{T_{0,-}} \in
dy,\; T_{0,-} <\I ] W(u-y)=P_x [- \und {X}_{\eta} < u,\; T_{0,-} <\I
	]W(u).
\end{equation*}
\qed

\beR \la{r:twos} We end this subsection by noting  that showing that the function  defined by \eqref{sRui}
has \LT \ \eqref{WLT} (up to a constant), is not trivial.

The first construction via excursion theory is due to \cite[Thm.~VII.8]{Ber}. Other elegant solutions are due to
  \cite{nguyen2005some}, who used a Kennedy  type martingale, and
 to \cite[(3)]{Pispot}, who constructed the scale function as
\beq \la{Wu}  W_{\q}(x)= \Fqp e^{\Fq x} -u_\q(-x)={\Fqp}( e^{ x \Fq}-  \P_x \le[{\tzh  < \kil_\q} \ri]), \; x \geq 0\eeq
where   $u_\q$ is the potential density -- see \eqref{Wsim} below for a proof of the last formula, which can be easily implemented   via Monte Carlo simulation\fn[3]{\eqref{Wu} holds trivially for $ x \in \R_- $ as well, when it reduces to $\P_x \le[{\tzh  < \kil_\q} \ri]=e^{ x \Fq}$, which  may be interpreted as the
value of a payment of $1$ at the  hitting time $\tzh $.}\fn[4]{Noting finally that $u_{\q}(x), x \in \R_+$  is exponential given by $u_{\q}(x)=\Fqp e^{-\Fq x}, \; x \geq 0$ and letting $u_{\q}^+(x)=\Fqp e^{-\Fq x}, x \in \R$ denote the analytic continuation of  $u_{\q}(x), \; x \geq 0$ yields yet another representation $ W_{\q}(x)=   u_{\q}^{+}(-x) -u_\q(-x)$ \cite{avram2002valuation}.}.

The simplest  solution maybe is  to reduce to the case $q=0$ by using the
easily checked Esscher transform relation
\beq  W_{\q}(x) =e^{x \, \Fq} \TW (x). \la{TW}\eeq
Here  $\TW (x)$ denotes the $0$-scale function with respect to the  ``Esscher transformed" measure
$P^{(\Fq)}$ (in general, the transform $P^{(r)}$ of the measure $P$ of a \lev process with \LE $\k(s)$ is the measure of the \lev process with \LE $\k(s+r)-\k(r)$, with $r$ in the domain of $\k(\cdot)$ \cite[Prop. 4.2]{AA}, \cite[3.3 pg.83]{Kyp}).

The advantage of $\TW (x)$ is that this is
 a monotone bounded function, with values in the interval $(\lim_{s\to \I} \fr s{\k(s)}, \fr 1{\k'(\Fq)})$. \Thr for numerical  computation of $ W_{\q}$ it will be  useful to replace it by $\TW(x)$, with \LT
$$
\HW(s)=\fr 1{\k(s+ \Fq)- \q}=\fr 1{\k(s+ \Fq)- \k(\Fq)}:=\fr 1 {\k^{(\Fq)}(s)},
$$
(removing thus the exponential growth). Pad\'e and Laguerre approximations of \eqref{TW} are provided in \cite{AHPS}.

Another probabilistic interpretation of  \eqref{TW}  is
\beq \la{WF} W_{\q}(x)=e^{x \Fq} \T L_\q(x),
\eeq
where
$\T L_\q(b)=
E \Big[\int_0^{\tb} e^{-\q t} d L^0_t\Big]=L^0 _{\tb \wedge \kil_\q}=\Fqp - e^{-\Fq x} u_\q(-x)$
is the expected discounted occupation time at $0$, starting at $0$, before up-crossing the level $b$ \cite[V(18)]{Ber}. This relation extends   to the spectrally negative \MAP (SNMAP) context \cite[(2),(12)]{IP}, has been  used      for computing numerically the SNMAP matrix scale function \cite{Ivapackage}.

\eeR

\beR {\bf $\Fq$ and the  other roots of the \CL equation} $\k(s)-\q= 0$ play a central role in asymptotics computations. Clearly,
$\Fq$ is the asymptotically dominant singularity
$$ W_{\q}(x) \sim  \Fqp e^{x \Fq}=\fr{e^{x \Fq}}{\k'(\Fq)}, \qu x\to \infty.$$

The other poles of the \RHS of \eqr{WLT} (the roots of the \CL equation) intervene, when they exist, in the asymptotics of the eventual ruin probabilities when $\k'(0_+) >0$ and in their numerical  approximations -- see for example \cite{AAK,avram2011moments,avram2012matrix,avram2014survey,avram2018ruin}.

 \eeR
\ssec{Two resolvents in terms of  the $W_\q(x)$ function.\la{s:Wr}}

We will recall here two fundamental resolvent formulas expressed  in terms of $ W_{\q}$.
Resolvents are at a level of  sophistication above the other concepts reviewed in this paper, and
these results will not be proved. However,  once accepted,
they provide us with a convenient point of entrance in our topic.

We introduce first a notation style used throughout the paper.

\beR \la{r:not}  Our cookbook will require notations for several types of boundaries for example absorbing, reflecting, refracting, and  Parisian/Poisonian   stopping or reflecting. To deal with these five cases, it is convenient, following \cite{Ivapot}, to append  the state space to the specification of  a process; the five cases above
  will be denoted below by  $b|, b], b[, b \vdots, b\}$ for an upper boundary, and for a lower boundary by  $|a, [a, ]a, \vdots a , \{ a$ .
 For drawdown  boundaries, the respective notations will be  $\ovl d, \H d, \breve d, \ddot d, \T d$.  Note that the term ``boundary" for the refracting and Parisian cases is meant in the sense of a discontinuous ``regime switching" in the drift and killing parameters of the process, respectively. This convention gives suggestive notations  when composing several mechanisms. For example, for  the ``classic reflection above at $b$, with Parisian reflection below at $b_0$ and absolute ruin at $a< b_0 $" studied in \cite{APY,PYrisks},  the notation for the corresponding state space would be $|a,\{b_0, b_1]$. Such boundaries are useful in optimal control \cite{perez2018optimal}.

 Note  that  absorption delimiters like $|a$ and  $b|$ and  may  and will be often omitted without confusion (so the default for an unspecified end-point is absorbing).  
 \eeR

\beP \la{p:res} Put $W_{\q} (x,a)=W_{\q} (x-a)$ (as a reminder that these formulas hold also for space-inhomogeneous models, like for example for refracted processes  \cite{LZ17})\fn[4]{One of the nice things about the toolkit is that switching to inhomogeneous skip-free processes just requires changing $x-a$ to $x, a$. The only thing specific to L\'evy (and refracted) setting is that $W$ is quasi-explicit.}.

 A) For any bounded interval $[a,b]$ and any Borel set $B\subset [a,b]$, let $$U_\q^{|a,b|}(x,B)=\E_x \left[ \int_0^{ \ta \wedge \tb } e^{-qt} \1_{\left\{X_t  \in B  \right\}} \diff  t \right],$$  denote the {\bf $q$-resolvent of the  spectrally negative  L\'{e}vy process} killed outside the interval $[a,b]$. Then \cite{Suprun}, \cite[Thm.~1]{Ber97}, \cite[Thm.~8.7]{Kyp}, \cite[(14)]{Ivapot}, \cite[Thm.~2.2]{LP}, \cite[Thm.~1]{LZ17},
$U_\q^{|a,b|}(x,B)=\int_{a}^{b} \1_{\{y \in B\}} u_\q^{|a,b|}(x,y) \diff  y$, with resolvent density
\begin{align} \label{res} 
\begin{split}
	u_\q^{|a,b|}(x,y) &= \fr{P_x[X_{\kil_\q} \in dy]}
{\q \; dy}= \frac {W_{\q}(x,a) } {W_{\q}(b,a)} W_{\q} (b,y) -W_{\q} (x,y).
	\end{split}
\end{align}

Note also  the following identities in limiting cases --
see for example \cite[Chapter~8.4]{Kyp}:
%
\begin{eqnarray}
(\q \,d y)^{-1} P\bigl(X_{\kil_\q}\in d y\bigr)&=&\Phi'(\q) \,e^{-\Phi(\q) y}-W_\q(-y),\label{eqpotential}
\\
(\q \,d y)^{-1} P\bigl(X_{\kil_\q}\in d y,\kil_\q<\tb
\bigr)&=&e^{-\Phi(\q)
b}W_\q(b-y)-W_\q(-y),\label{eqpotential1}
\\
(\q \,d y)^{-1} P_b\bigl(X_{\kil_\q}\in d y,\kil_\q<
\tz\bigr)&=&e^{-\Phi(\q) y}W_\q
(b)-W_\q(b-y),\label{eqpotential2}
\end{eqnarray}
where $b>0$ and the killing rate $q\geq0$ is implicit.

B) The {\bf $q$-resolvent of a spectrally negative  L\'{e}vy process absorbed below at $a$  and reflected above at $b$} (see \eqref{ref} for definition of reflection) has the resolvent
density \cite[(21)]{Ivapot}, \cite[Thm.~2.4]{LP}
\begin{align} \label{resb} 
\begin{split}
	u_\q^{|a,b]}(x,y) &=  \frac {W_{\q}(x,a)} {W'_{\q}(b,a)} \le(W'_{\q} (b,y)+  W_{\q}(0) \delta_b(d y)\ri) -W_{\q} (x,y),
	\end{split}
\end{align}
where the derivative is taken with respect to the first variable.
\eeP

 \beR
Letting $b \to \I$ in \eqref{res} we find the resolvent on intervals bounded only below for any Borel set $B\subset [a,\I)$, which is closely related to Dickson's formula in the actuarial literature
\begin{align} \label{resi}
\begin{split}
U_\q^{|a}(x,B) = \int_{a}^{\I} \1_{\{y \in B\}} u_\q^{|a}(x,y) \diff  y, \qu u_\q^{|a}(x,y) &=  W_{\q}(x-a) e^{-\F(\q) (y-a)}  -W_{\q} (x-y).
	\end{split}
\end{align}

\eeR
\beR
For other resolvent laws involving all possible combinations of boundary conditions (reflection or/and absorbtion), see \cite{Kyp,Ivapot,LP}. Note that the   proofs use typically \exc \ theory. One exception is \cite[Thm.~4.1]{perez2018optimal}, who compute the resolvent density $u_{\q,\r}^{]0}(x,y)$
with Parisian reflection at  Poisson observation times of intensity $\r$. The proof uses the \Mar \  property in the  bounded variation case,
and a \LT \ approach in the unbounded variation case. \eeR

\sec{Obtaining the  $Z_\q(x)$ function in terms of $W_\q(x)$ by using the resolvent  \la{s:Z0}}

  The first resolvent formula will now be used to introduce
  the second pillar of this theory, the  scale function $ Z_{\q}$, which intervenes in the ``non-smooth-exit law" below. Using this together with the ``smooth-exit law" \eqr{twos}
  will be essential in deriving the other recipes offered below.
  \beP
   \la{p:Z} A) {\bf The \LT \ of the time until the lower
  boundary $0$, if this precedes an upper boundary $b>0$},  is given by \cite[(10)]{AKP}

\beq \label{RZ} &&
\Rui_\q^b(x):=\E_x\left[e^{- \q \tz} ; \tz <\tb \right] = Z_{\q}(x ) -   \fr { W_{\q}(x)}{ W_{\q}(b)} Z_{\q}(b),
\eeq
where $ Z_{\q}(x )=1+ \q \ovl  W_{\q}(x)$, $\ovl  W_{\q}(x)=\int_0^x   W_{\q}(u) d u$.

B) {\bf The \LT \ of the time until the lower
  boundary $0$ in the presence of reflection at an upper boundary} $b \geq 0$ is
  \cite[Prop.~5.5]{APP15}, \cite[Thm.~6]{IP}
\beq\label{bdruinrefl} &&
\Rui_\q^{b]}(x):=\Eb_x\left[e^{- \q \tzb } \right]= Z_{\q}(x ) -   \fr { W_{\q}(x)}{ W'_\q(b)} Z'_\q(b ),
\eeq
where $\Eb$ denotes expectation for the process reflected from above at $b$ and
\beq \la{Trb} \tzb =\tz \; \1_{\{\tz< \tb\}}+ \t_b \;  \1_{\{\tb<\tz\}} \eeq
denotes the first passage  below $0$ under this measure (recall that $\t_b$ is a drawdown  time \eqref{dd}, or, equivalently, the time when the process starting at $b$ and Skorokohod reflected at $b$ is ruined\fn[4]{When $x=b$, \eqr{Trb} simplifies to $T_0^{b]} = \t_b $.}).

\eeP

Here is a proof of Proposition  \ref{p:Z}, borrowed from \cite{LZ17} (who consider the more general case of Omega models).

{\bf Proof}.  A): Put $T =\min(\tz,\tb)$, and consider the elementary
identity:
\beq \la{resost}\int_0^T \q e^{-\q t} d t=1-e^{-\q T}.\eeq

By denoting
\begin{equation}\label{Wovl}
\ovl W_{\q}(x ):=\int_0^x W_{\q}(y)dy,
\end{equation}
taking expectation and using the resolvent formula \eqref{res}, we get
\bea \q \int_0^b u^{|0,b|}_\q(x,y) d y&=&1-\fr{ W_{\q}(x)}{ W_{\q}(b)}- \Rui_\q^b(x) \Eq \\\q \le( \fr{ W_{\q}(x)}{ W_{\q}(b)} \int_0^b
 W_{\q}(b-y )d y-\int_0^x
 W_{\q}(x-y )d y\ri) &=&1-\fr{ W_{\q}(x)}{ W_{\q}(b)}- \Rui_\q^b(x) \Lra\\\Rui_\q^b(x)=1-\fr{ W_{\q}(x)}{ W_{\q}(b)}-\q \le( \fr{ W_{\q}(x)}{ W_{\q}(b)} \ovl
 W_{\q}(b)-\ovl
 W_{\q}(x )\ri)&=&1+ \q \ovl
 W_{\q}(x ) -\fr{ W_{\q}(x)}{ W_{\q}(b)}\le(1+ \q \ovl
 W_{\q}(b )\ri).\eea
Putting now $ Z_{\q}(x )=1+ \q \ovl  W_{\q}(x)$ yields the result.

B): Applying the same steps to $\tzb$, we find
\bea \E_x \Big[\int_0^{\tzb} \q e^{-\q t} d t\Big]=\E_x \Big[1-e^{-\q \tzb}\Big]=1- \Rui_\q^{b]}(x) &=& \\\q \le( \fr{ W_{\q}(x)}{ W_{\q}'(b)} \le(\int_0^b
 W_{\q}'(b-y )d y+  W_{\q}(0)\ri)-\int_0^b
 W_{\q}(x-y )d y\ri) & \Lra&\\\Rui_\q^{b]}(x)=1-\q \le( \fr{ W_{\q}(x)}{ W_{\q}'(b)}
 W_{\q}(b)-\ovl
 W_{\q}(x )\ri)&=& Z_{\q}(x ) -\fr{ W_{\q}(x)}{ W_{\q}'(b)} Z_{\q}'(b).\eea

 \beR These two proofs illustrate the very important  method  of   integrating  resolvent densities -- see \cite{Ivapot} for a compendium of resolvent formulas. For a  direct proof not using resolvents, in the case of Brownian motion, see  \cite[Thm 1.1]{Ebe}.  \eeR

 \beR Note the similar structure of \eqref{RZ} and \eqref{bdruinrefl} (a phenomenon which will keep  recurring below).   Formally, switching from absorption at $b$ to  the measure $\Eb$ involving reflection at $b$ only requires switching  the respective boundary conditions $ \Rui_\q^{b}(b)=0,\le(\Rui_\q^{b]}\ri)'(b)=0.$
  Now the first boundary condition is obvious, like any absorbtion boundary condition, but not the second.

Let us examine now a ``failed direct approach" to establish \begin{equation} \le(\Rui_\q^{b]}\ri)'(b)=0 \Eq \Rui_\q^{b]}(b-\e)-\Rui_\q^{b]}(b)=o(\e). \la{bcr} \end{equation}
 Using now the decomposition \eqr{Trb} yields
 \bea &&
\Rui_\q^{b]}(b) -\Rui_\q^{b]}(b-\e)=
 \Rui_\q^{b]}(b)-\le(Z_{\q}(b-\e) -   \fr { W_{\q}(b-\e)}{ W_\q(b)} Z_\q(b ) +   \fr { W_{\q}(b-\e)}{ W_\q(b)} \Rui_\q^{b]}(b)\ri)\\&&=\Rui_\q^{b]}(b) \left(1-  \fr { W_{\q}(b-\e)}{ W_\q(b)}\right)-\le(Z_{\q}(b-\e) -   \fr { W_{\q}(b-\e)}{ W_\q(b)} Z_\q(b ) \ri)\\&&=\e \le[\Rui_\q^{b]}(b) \fr { W_{\q}'(b)}{ W_\q(b)}+\le(Z_{\q}'(b) -   \fr { W_{\q}'(b)}{ W_\q(b)} Z_\q(b ) \ri)\ri]+o(\e)=\e \le[Z_{\q}'(b) +\fr { W_{\q}'(b)}{ W_\q(b)} \le(\Rui_\q^{b]}(b)-   Z_\q(b ) \ri)\ri]+o(\e).
\eea
The \bco  on the derivative is equivalent thus to the \bco  on the function
$\Rui_\q^{b]}(b)=\Eb_b\left[e^{- \q \tzb } \right]= Z_{\q}(b ) -   \fr { W_{\q}(b)}{ W'_\q(b)} Z'_\q(b)$, which we wanted to avoid establishing. A more sophisticated approach is thus needed.
For the \CL model,  the  \bco \ \eqr{bcr} on the derivative has been established in  \cite{lin2003classical}, using the regenerative property of the Poisson process at claim instants (their proof is quite ingenious).  For spectrally negative \lev processes,  the use of \exc \ theory   seems unavoidable.
\eeR

\beR The Propositions \ref{p:twos}-\ref{p:Z} and most of the results  in this review  may be modified    to apply formally to the context of spectrally negative and \sp \Mar processes,  which include for example the  {\bf  continuous state-space branching processes} (CSBP) -- see for example \cite[Ch. 12]{Kyp} (in particular Thm.~12.8), and the  continuous-state branching processes
with immigration (CBI)
  introduced by Kawazu and
Watanabe \cite{Kawazu}, which may characterized in terms
 of two \LEs $\psi$, $\k$ of \sp \lev processes. However, while $W,Z$ exist (as functions of two \vars),  no straightforward method for their  computation 
 is available.\fn[4]{Recall that CSBPs are characterized by generators of the form $x \psi(D)$, where $\psi(D)$ is the generator of a \sp \lev process, and that they may be obtained from \sp \lev process by a time-change called the Lamperti transformation -- see \cite{CLUB}. This acts on the Skorokhod
space $\mathbb{D}$ of c\`adl\`ag trajectories with values in $E=[0,\I]$, as follows:
for any $f \in \mathbb{D}$,  introduce the additive functional $I$ and its inverse $\ola I$, given by
$I_t=I_t(f):=\int_0^t f(s) d s \in[0,\I], \qu \ola I_t=\ola I_t(g):=\inf\{s \geq 0: I_s(g) >t\}=I_t(\fr 1 g) \in[0,\I].$
The Lamperti transformation $L : \mathbb{D} \to \mathbb{D}$
is defined by
$L(f) = f \circ \ola I$
(note that $L(f) (t) = f(\I)$ if $\ola F_t = \I$, so  that $0,\I$ indeed are  absorbing for $L(f)$).
It may be checked that $L$ is a bijection of
$\mathbb{D}$, with  inverse given by $L^{-1}(g)= g \circ I (g) $. An extension to the CBI case is offered in \cite{caballero2013lamperti}.
 However, the Lamperti transformation seems too complicated to yield a method for the  computation of $W,Z$     in terms of the \lev Laplace exponents. It is intriguing to investigate whether simple formulas for $W,Z$ are available in these cases at all.}
\eeR
\beR \la{r:tmin}
  Adding \eqref{twos} and \eqref{RZ},  we find that for $T =\min(\tz,\tb)$
  \begin{equation} \la{tau} \E_x \left[ e^{- \q T} \right]=P_x[ T \leq \kil_\q]= Z_{\q}(x ) -   \fr { W_{\q}(x)}{  W_{\q}(b)}( Z_{\q}(b )-1)=1-\q  \left(\fr{ W_{\q}(x)}{  W_{\q}(b)} \ovl W_{\q}(b)-\ovl W_{\q}(x)\right),\end{equation}   which recovers \cite[Cor.~1]{Ber97}  (up to the omission of $\q$ there). Since this must be less than $1$, it follows that the function $\fr{ \ovl W_{\q}(x)}{ W_{\q}(x)} $ is increasing, or, equivalently, that $\ovl  W_{\q}(x)$ is log-concave, and
  \beq \fr{ W_{\q}'(x) \ovl  W_{\q}(x)}{ W_{\q}^2(x)} <1 \la{posW} \eeq
  see also \cite{loeffen2010finetti}.

   For a second probabilistic proof of \eqr{posW}, consider the  time from $b$
to $0$ of a reflected process \eqref{bdruinrefl},
   which is equal in law to the drawdown  time  $\t_b$\fn[3]{That is easily understood by fixing the maximum at $b$, which changes the negative of the drawdown  into the Skorokhod reflected process.}.
Choosing  $x=b$  in \eqref{bdruinrefl} yields
   \begin{equation} \la{d}
{\de}_{\q}(b):=\Rui_\q^{b]}(b)=\E_0\left[e^{- \q  \t_b } \right]= Z_{\q}(b) -\fr{ W_{\q}(b) Z_{\q}'(b)}{ W_{\q}'(b)}=1-  \q \left(\fr { W_{\q}^2(b) }{ W'_\q(b)}-\ovl  W_{\q}(b ) \right).
\end{equation} Since this must be less than $1$, the \nny of the term in parenthesis follows.

\eeR

 {\bf Reduction of first passage problems to the computation of the solutions $ W_{\q}$ and $ Z_{\q}$ of TSE}. It turns out that
 the solutions of a great variety of first passage problems reduce ultimately to the solutions of the two-sided smooth and non-smooth first passage problems of exit from a bounded  interval (TSE). Thus, they may be expressed in terms of  $ W_{\q}$ \cite{Ber97}, and further simplified by the introduction of the  second scale function $ Z_{\q}$ \cite{AKP}.
Many calculations and inversions of \LTs may be  replaced for spectrally negative \lev processes by
  the computation of  the $W$ and $Z$ scale functions -- see \cite{Pisexit,Pispot,Pisexc,APP,IP}, to cite only a few papers.  Furthermore, the formulas reviewed  hold  as well for spectrally negative Markov additive processes,  where the appropriate matrix scale functions were identified in \cite{KP,Iva,IP}, for random walks (the compound binomial risk model) \cite{AV}, and for   positive \ssr \ Markov processes with one-sided jumps \cite{vidmar2018temporal,vidmar2018exit}.

Somewhat surprisingly, it appeared recently that
  the recipes reviewed below apply equally to spectrally negative \lev processes with (exponential) Parisian absorbtion or  reflection below \cite{landriault2014insurance,
 AIZ,albrecher2017strikingly,BPPR,APY}, with  the appropriate scale functions $W,Z$ identified in \cite{APY,AZ}.
 This  mystery was explained in \cite{LP,LZ17,vidmar2018first}, who showed that
 the $W,Z$ recipes appropriately extended apply  to the general class of Omega models, of which  Parisian Poissonian models are a particular case.  In fact, the second paper considers
 even more general models with refraction \cite{KL,KPP}.

\sec{The three variables $Z_{\q}(x,\th)$ scale function/Dickson-Hipp operator  applied to $W_\q(\cdot)$ \la{s:Z}}
Let $\H {_x W_{\q}}(\th)$ denote the Laplace transform  of the shifted scale function $_x W_{\q}(y):=W_{\q}(x+y)$ (the composition of shift with \LT \ is also called
   Dickson-Hipp operator).

When the \LT \; $\E_x [ e^{\th X_{\tz} }]$ of the first  position of the process after exiting $[0, \I)$ is of interest, 
one ends up  working with the
 {\bf two variables $ Z_{\q}$ scale function} \cite{AKP,IP},\cite[Cor.~5.9]{APP15}, defined for $\theta\in \C$ such that the real part $\Re(\th)> \F(\q)$ (to ensure integrability) by:
\begin{equation}\la{Zt}
\begin{aligned}
Z_{\q}(x,\th)&=  \le(\k(\th)-\q\ri)\int_0^\I e^{- \th y}  W_{\q}(x+y) dy :=\le(\k(\th)-\q\ri) \H {_x W_{\q}}(\th)\\&= \fr{\k(\th)-\q}{\th-\Fq}  W_{\q}(x)+\E_x \le[e^{-\q \tz +\th X_{\tz }} \; \1_{\{\tz < \I\}} \ri]\\&=\fr{\k(\th)-\q}{\th-\Fq}  W_{\q}(x)+ \Rui_{\q,\th}(x), \quad \Re(\th)> \F(\q).
\end{aligned}
\end{equation}
(see Corollary \ref{r:asycreep} A) for the proof of the last decomposition.)
Thus, up to a constant,    $Z_{\q}(x,\th)$  is    the Laplace transform $\H {_x W_{\q}}(\th)$ of the shifted scale function $_x W_{\q}(y):=W_{\q}(x+y)$,  and  the normalization ensures that  $Z_{\q}(0,\th)=1$.

 \beR The first term in the decomposition above is asymptotically dominant for $\q>0$.  The second term simplifies in the \CL case  when $x=\q=0$ to
 \bea \E_0 \le[e^{\th X_{\tz}} \; \1_{\{\tz < \I\}} \ri]= 1 -   \fr {\k(\th)}{ c \th}= \fr{\H{\ovl \lm}(\th)}c, \ \for \th >0,\eea
 identifying the well-known \LT of the deficit at ruin starting from $0$ for the \CL process, where $\H{\ovl \lm}$ denotes the Laplace transform of the tail of the L\'evy measure $\ovl \lm(y) = \lm(y,\infty)$. \eeR

  The analytic continuation of \eqref{Zt} is
 \begin{equation}Z_{\q}(x,\th)= e^{\th x} + \big(\q- \k(\th)\big) \int_0^x e^{ \th (x-y)} W_{\q}(y) dy, \; \th \in \C. \la{Z} 
 \end{equation}
This implies that \begin{equation} \la{Zexp}\bc Z_{\q}(x,\th)= e^{\th x}, &x\leq 0\\   Z_{\q}(x,\Fq)= e^{ x \Fq}, &x \in \R \ec.\end{equation} 
\beR
We can also identify $ Z_{\q}(x,\th)$ via its Laplace transform in $x$:
\bea &&\H{Z_{\q}}(s,\th)=({s-\th})^{-1}(\k_\q(\th)^{-1}-\k_\q(s)^{-1})  \k_\q(\th) =
\fr{\k(s)-\k(\th)}{s-\th} \fr 1 {\k(s)-\q}, \; \k_\q(s):=\k(s)-\q\\&& \Lra \H Z_{\q}(s)=s^{-1} \k(s) \k_\q(s)^{-1} .\eea
\eeR

We list now   some useful easy to check formulas involving $Z_{\q}(x), Z_{\q}(x,\th)$:

\beq \la{ovl}
 Z_{\q}(x)& =&  1 + \q \overline{W}_\q(x)= c  W_{\q}(x) + \fr{\s^2}2 W'_\q(x)  -\int_0^x  W_{\q}(x-y)  \ovl \lm(y) dy, \la{oldZ} \\
\overline{Z}_\q(x)&:=& \int_0^{x}  Z_{\q} (z) d z = x + \q \int_0^{x} \int_0^z  W_{\q} (w) d w d z \la{ZbarD}\\Z_{\q}^{(1)}(x) &=& \fr{\partial Z_{\q}(x,\th)}{\partial \th}_{\th=0}= \ovl   Z_{\q}(x) -\k'(0_+) \ovl   W_{\q}(x), \la{Z1}
\\ Z_{\q}'(x,\th)&=&\th Z_{\q}(x,\th) +(\q- \k(\th)) W_{\q}(x),\la{Zder}
\eeq
where $'$ denotes here and  below derivative with respect to $x$ and $\ovl \lm(y) = \lm(y,\infty)$. The second formula for $ Z_{\q}(x)$ is a particular case of \eqr{Gdec}. Let us check it now when  $\s>0$:
\bea
&&1 + \q \overline{W}_\q(x)=1+\int_0^x \mG \le({W}_\q\ri)(y)dy\\
&&= 1 + \frac{\sigma^2}{2}  (W'_\q(x)-W'_\q(0_+)) + c  (W_{\q}(x)-W_\q(0_+)) + \int_0^\infty \left( \int_0^x W_q(y-z) dy - \int_0^x W_q(y) dy  \right) \Pi(dz)\\
&&= 1+ \frac{\sigma^2}{2}  (W'_\q(x)-W'_\q(0_+)) + c  W_{\q}(x) + \int_0^\infty \left( - \int_{x-z}^x W_q(y) dy  \right) \Pi(dz)-\1_{\{\int  \Pi(dz)<\infty \& \s=0\}}\\
&&=  \frac{\sigma^2}{2}  (W'_\q(x)-W'_\q(0_+)) + c  W_{\q}(x) -\int_0^x  W_{\q}(x-y)  \ovl \lm(y) dy +\1_{\{ \s>0\}}\\&&=c  W_{\q}(x) + \fr{\s^2}2 W'_\q(x)  -\int_0^x  W_{\q}(x-y)  \ovl \lm(y) dy,
\eea
where we integrated $\q {W}_\q(x)=\mG \le({W}_\q\ri)(x)$ with $\mG$ given in \eqref{genrisk}, and used Fubini, integration by parts, and \eqr{asymptotics_zero}.

\beR For Brownian motion, \eqr{oldZ} yields
$$Z_{\q}(x) = c  W_{\q}(x) + \fr{\s^2}2 W'_\q(x) =c  W_{\q}(x) + \fr{W'_\q(x)}{W'_\q(0)}  .$$
\eeR

\beR \la{r:tr} Note that  for  $ x \leq 0$, it holds that
$\overline{W}_\q(x) = 0, \;  Z_{\q}(x) = 1, \; \overline{Z}_\q(x) = x$,
and  that $Z_{\q}(x,\th)$ is proportional to an Esscher transform; indeed, it is easy to check that
$ W_{\q- \k(\th)}^{(\th)}(x)=e^{-\th x}  W_{\q}(x), \; Z_{\q- \k(\th)}^{(\th)}(x)=e^{-\th x} Z_{\q}(x,\th)$. Recall that the Esscher transform refers to an exponential change of measure using the martingale $e^{\th X_t - \kappa(\theta) t}$, $t \geq 0$. For each $\theta $  in the domain of $\k(\cdot)$, the process $X$ remains in the class of spectrally negative L\'evy processes,  is characterized by the \LE $\kappa(\cdot +\theta)- \kappa(\theta),$ and $ W_{\q}^{(\th)}$, $ Z_{\q}^{(\th)}$ denote the scale functions of $X$ under this change of measure.\fn[5]{Before the introduction of the notation $Z_{\q}(x,\th)$ in \cite{IP,APP15}, results were  expressed in terms of  Esscher transformed scale functions.}
 \eeR

\beR It is easy to check by taking \LT \ \cite{Pisexit,LRZ} that the convolution of two $W$ scale functions satisfies the equation
\be W_\q*W_\l(x)=\fr{W_\q(x)-W_\l(x)}{\q-\l}.  \ee

The analogue formula for the convolution of two $Z$ scale functions is more complicated. When $\s=0$, \ith
\bea (Z_\q*Z_\l)(x)=(\fr{Z_\l-Z_\x}{\q-\l})*(\ol \lm)(x).  \eea

However the  convolution of $W$ and  $Z$ is again simple \cite[Lem. 4.1]{AIZ}
\be W_\q*Z_\l(x,\th)=\fr{Z_\q(x,\th)-Z_\l(x,\th)}{\q-\l} . \ee

\eeR

{\bf The history of Z}.
   The second  scale function  $ Z_{\q}(x)$  was introduced in the thesis of M. Pistorius (which the first author codirected with A. Kyprianou),  as a means of  expressing in a simpler  way  both the results of \cite{Suprun,Ber97} and some new results involving   reflected processes and drawdown stopping (used ``Russian options''). See \cite[(6)]{AKP} for the first published reference. Its importance became clearer after its  further  use in \cite{Pisexit,Pispot,kyprianou2005martingale,nguyen2005some,
   Doney,Pisexc}.

   By some  historical error, all these papers, as well as the textbook \cite{Kyp}, omitted the information that
   the  "birth certificate" of the function $Z$ was signed in the thesis of Pistorius and in  \cite{AKP}. Instead, reference was made to the pioneering work \cite{Ber97}, which however contains no $Z$ function.

   The   three variables extension $Z_{\q}(x,\th)$   was introduced essentially in \cite{AKP} as an Esscher transform of $ Z_{\q}(x)$ -- see Remark \ref{r:tr}.
   Then, the simultaneous papers \cite{IP} and  \cite[Cor.~5.9]{APP15} (first submitted  in 2011, ArXiv 1110.4965)  proposed the direct definition \eqref{Z}, without the Esscher transform from previous papers.

   Subsequently,    $Z_{\q}(x,\th)$ was shown in \cite
  [Thm.~5.3, Cor.~5.9]{APP15} to be  a particular case  of a ``smooth Gerber-Shiu function'' \cite[Def.~5.2]{APP15} associated to an {exponential payoff} $e^{\th x}$. More precisely, $Z_{\q}(x,\th)$ is  the unique ``smooth'' solution of
\beq \la{hext}\bc (\mG-\q I) Z_{\q}(x,\th)=0, & x \geq 0\\ Z(x,\th)=e^{\th x},& x \leq 0\ec,\eeq
where $\mG$ is the Markovian generator \eqref{genrisk} of the process $X_t$ -- see \cite[(1.12), (5.23), Sec.~5]{APP15} and Section \ref{s:nonh}.

   $Z_{\q}(x,\th)$ was used first as generating function for the smooth Gerber-Shiu functions associated to  power  rewards $1, x, x^2$, which were denoted respectively by $Z_q, Z_q^{(1)}, Z_q^{(2)}, \dots$. Subsequently,
   it started  being used intensively in exponential Parisian ruin problems following the work of \cite{AIZ}.

    As of recently, several papers \cite{APP,KL,Iva,IP,Ivapot,AIZ,AIpower,APY,AZ}
    showed that  L\'evy formulas expressed in terms of  $ W_{\q}(x)$ and  $ Z_{\q}(x)$ or $Z_{\q}(x,\th)$   hold also for doubly reflected processes\fn[4]{for the construction of these, one may use a recursive approach, or the recent paper \cite{kruk2007explicit}},  refracted processes, spectrally negative \MAP,  processes with Parisian absorption or reflection, and combinations of these features.
More precisely, formulas which hold for the L\'evy model continue to hold for the others,  once  appropriate
{(matrix) scale}  functions are identified.

  We will call this body of related first passage formulas the scale  functions kit or cookbook.  Its availability  means that the analytic work required to solve a first passage problem  may often be replaced by looking up in the  cookbook.
   The next section contains  \nun of our favorite recipes. 

\sec{\Nun first passage laws \la{s:10}}

We will start with the  easiest problem, which involves only $W_q\cd$.

\ssec{Expected discounted dividends}

We  review now expected discounted dividends $U$ 
under both  reflection and absorbtion regimes. These are especially important in the control of reserves processes -- see Section \ref{s:div}.

\beT \la{l:div} A) The {\bf expected total discounted dividends up to $\tzb$} are given by
\begin{equation} V^{b]}(x):=\Eazb_x\le[\int_{[0,\tzb]}e^{-\q t}d   \U_t\ri]=\frac{ W_{\q}(  x)}{W_{\q}^ {\prime}(b)}, \la{div} 
\end{equation}
where  $\Eazb$ denotes the law of the process  reflected from above at $b$,
and absorbed at $0$ and below.

B) The {\bf expected total discounted dividends over an infinite horizon} for the doubly reflected process,  with expectation denoted  $\Ezb$, are given by   \cite[(4.3)]{APP}
\begin{equation}\label{divSLG}
V^{[0,b]}(x):=\Ezb_x\le[\int_0^{\I} e^{-\q t}d  \U_t\ri] =   \frac{Z_{\q}(x)}{  Z'_{\q}(b)}.
\end{equation}
\eeT

{\bf Proof}. A) Since $V^{b]}(x)= \fr { W_{\q}(x)}{ W_{\q}(b)} V^{b]}(b)$ by the smooth-exit law  \eqr{twos}, the essential part is proving the result for $x=b$, i.e. that $V^{b]}(b)= \fr { W_{\q}(b)}{ W_{\q}'(b)}=\nu_\q(b)^{-1}$, where the latter (\exc \ theoretic)  quantity has already
been introduced in
Remark \ref{r:exrate}. For the \CL case, a direct computation of $V^{b]}(b)$ is provided
in \cite[Lem 6.4]{KypGS}; for the \sn case, a
  generalization to all moments of the discounted dividends (using \exc \ theory) may be found  in \cite[Thm 10.3]{Kyp}.

To see the idea behind the  \exc \ theory proof, note, following \cite{albrecher2018linking}, that
\beq \la{iddiv} \Eb_x \le[ \int_0^{\tzb} e^{-\q t} d \U_t\ri] =\Eb_x \le[ \int_0^{\tzb \wedge e_\q}   d \U_t\ri] =\Eb_x \le[ \U_{\tzb \wedge e_\q}  \ri].\eeq
Finally, the law of variable $\U_{\tzb \wedge e_\q} |x=b$   is exponential with parameter
$\nu_\q(b)$, cf. Remark \ref{r:exrate} (see also
Theorem \ref{l:DP} A) below for a generalization).

B) Again, it is enough to prove the result for $x=b$, since
\begin{align*} V^{[0,b]}(x) &= \Ezb_x\le[\int_0^{\I} e^{-\q  t}d  \U_t\ri]=
\Ezb_x\le[\int_{\tbz}^{\I} e^{-\q  t}d  \U_t\ri]=\Ezb_x\le[ e^{-\q  \tbz} \int_0^{\I} e^{-\q  t}d  \U_t\ri]\\
&=\fr{ Z_{\q}(x)} { Z_{\q}(b)} \Eb_x \le[ \int_0^{ e_\q}   d \U_t\ri]=\fr{ Z_{\q}(x)} { Z_{\q}(b)} \Eb_x \le[ \U_{e_\q} \ri].
\end{align*}
It turns out that for $x=b$, the variable $\U_{e_\q} $ under the measure $ \Ezb_x$ is exponential with parameter
$\fr{ Z_{\q}'(b)} { Z_{\q}(b)}$, yielding the result (see Theorem \ref{l:joiDB} and Remark \ref{r:ZpZ} below for a generalization and further references). \qed

\beR \la{divb} Since  the \bco  $V^{b]}(b)=\fr { W_{\q}(b)}{ W_{\q}'(b)}$ in  A) requires \exc \ theory, one might  try to  establish instead the simpler \boun \con on the derivative
\begin{equation} (V^{b]})'(b)=1,\end{equation}
which says roughly that $$V^{b]}(b)-V^{b]}(b-\e)\sim \e.$$
\BEN \im
Let us start with the \CL model, and  follow the  derivation   suggested  in \cite{gerber2006note}, which  note that when starting from $b-\e$, no dividends are gained during a period of $\fr \e c$, while  when starting from $b$, dividends roughly equal to $c \fr \e c=\e$ are gained during this period.

\Mp construct  the processes starting from $b$ and $b- \e$ on the same probability space, and let $A$ denote the event
 that there is no jump in the interval $[0,\e/c]$. Over this event, the processes
 are coupled at time $\fr \e c$ and the only difference between the dividends comes from the interval $[0,\e/c]$.
 Putting now together the contribution over $A$ and over its complement yields:
 \begin{align*}&
 \fr{V^{b]}(b) - V^{b]}(b-\e)}{\e}= \\&\fr c{\e} \int_0^{\fr{\e}c} e^{- q s} ds + {\e}^{-1}\int_0^{\e/c} \l e^{-\l s-\q s} \int_0^b (V(b+ c s-x) -V(b-\e+ c s -x)) f(x) dx ds\\
 &\leq  \fr c{\e} \int_0^{\fr{\e}c} e^{- q s} ds + {\e}^{-1} \int_0^{\e/c} \l e^{-\l s-\q s} \int_0^b \fr{\l}c V(b+ c s-x) {\e} f(x) dx \to 1,
 \end{align*}
 where we used the  increasingness and locally Lifschitz  property of the value function \cite{Schmidli}, \cite[1.3, Prop.~1.3, p.9]{azcue2014stochastic},  in the \CL case.

\im We turn now to the spectrally negative Levy model. Armed with our two exit laws, we find:
\bea &&  V^{b]}(b-\e) =  \sme{b-\e}{b} V^{b]}(b) + \left( \nsme{b-\e}{b} \right)\times 0
\Eq \\&& \frac{V^{b]}(b) - V^{b]}(b-\e)}{\e}  = \frac{W_q(b)-W_q(b-\e)}{\e W_q(b)} V^{b]}(b) \Eq \\&&  (V^{b]})'(b) =  \fr{W_\q'(b)}{W_\q(b)} V^{b]}(b)\eea
and we fall back on the \prob of  tackling  $V^{b]}(b)$, suggesting  that the boundary condition is not trivial and that the use of \exc \ theory (see \cite{Ber}) is unavoidable in general. Note \how ~that  the perturbed \CL model was solved  in \cite{li2006distribution}, via a perturbation  approach.
\EEN

\eeR

\ssec{The total discounted capital injections/bailout  law, with non-smooth regulation}
The next  result \cite{Pisexit,IP}   shows the importance of $Z$  for reflected spectrally negative \lev processes.
It also provides a \gene of the fundamental survival \pro formula \eqr{sRui}.

\beT \la{l:refl} {\bf The \LT of the discounted capital injections/bailouts for the process reflected below}.
Let $X^{[0}_t$  denote the process reflected at $0$ \eqref{Xd} with regulator $\L_t=-\und X_t $, let $\E^{[0}_x$  denote expectation for this process and let
  \beq \la {Tra}  T_b^{[0}=\tb \; \1_{\{ \tb  <\tz \}}+ \und
  \t_b \;  \1_{ \{\tz <\tb\} }\eeq denote the first passage  to $b$ of $X^{[0}_t$ , to be called ``reflected up time''. The total capital injections into the    process reflected at $0$, until the first up-crossing of a level $b$ satisfy
\cite[Thm.~2]{IP}:
\begin{equation} \label{refbailout} 
\sRui_{\q,\th}^b(x,[0):=\Ez_x \left[e^{-\q T_b^{[0} - \th \L_{T_b^{[0}}}\right]=\Ez_x \left[e^{  - \th \L_{T_b^{[0}}};T_b^{[0} <e_\q\right]
=\bc \dfrac{Z_{\q}(x,\th)}{ Z_{\q}(b,\th)} & \th <\I\\\E_x \left[e^{-\q \tb} \1_{ \{ \tb < \tz \}}\right]=
\dfrac{W_{\q}(x)}{ W_{\q}(b)} &  \th =\I\ec.\fn[4]{The  result  \eqref{refbailout} above may be viewed as the fundamental law of spectrally negative \lev processes, since it implies the fundamental smooth two-sided exit formula \eqref{twos}. Note also that formally, replacing  absorption at the boundary $0$ by reflection leads
to replacing $W$ by $Z$;  this will be further confirmed in   several of the results below.}\end{equation}
\eeT

%

\beR  \la{r:proof}
Theorem \ref{l:refl} was first proved in  \cite[Thm.~2]{IP} as a consequence of a more general result \cite[Thm.~13]{IP}, but we prefer to use the observation   that it is essentially equivalent to  \eqref{sevruin}  \cite{IP}. Indeed,  \eqref{Tra} implies:
  \begin{equation} \la{proof} \Ez_x \left[e^{-\q T_b^{[0}-\th \L_{T_b^{[0}}}\right]=
\E_x \left[e^{-\q \tz +\th X_{\tz}} ;  \tz < \tb \right] \; \Ez_0 \left[e^{-\q T_b^{[0}-\th \L_{T_b^{[0}}} \right]+  W_{\q}(x) W_{\q}(b)^{-1}. \end{equation}
If the first term is known one gets an equation for the deficit at ruin
\begin{eqnarray*}&& Z_{\q}(x,\th) Z_{\q}(b,\th)^{-1}=  W_{\q}(x) W_{\q}(b)^{-1} + \E_x \left[e^{-\q \tz +\th X_{\tz}} ;  \tz <\tb \right] Z_{\q}(b,\th)^{-1},  \end{eqnarray*}
with the known solution
$
\E_x\left[e^{-\q\tz + \th X_{\tz}};\tz <\tb\right]=Z_{\q}(x,\th) -  W_{\q}(x) W_{\q}(b)^{-1} {Z_{\q}(b,\th)} $. And  if the deficit at ruin is known, one may use \eqref{proof} with $x=0$ to solve for $\Ez_0 [e^{-\q \tb-\th \L(\tb)}]$, provided that $ W_{\q}(0) \neq 0$. When $ W_{\q}(0) = 0$, one must start with a ``perturbation (approximation) approach'', letting $x \to 0$  \cite{zhou2007exit}-- see also Section \ref{s:pr}, where this result is proved directly, in the more general context of Parisian ruin.
\eeR


\ssec{Deficit at ruin}
We turn now to problems of deficit at ruin. We will present here a generalization of the ``non-smooth-exit law'', featuring
the $Z_{\q}(x,\th)$ function.
\beT \la{l:s}{\bf Deficit at ruin for a process absorbed or reflected at $b>0$}.

 A)
 The joint Laplace transform of the first passage time of $0$ and the undershoot for a process absorbed  at $b>0$
 is given by \cite[Prop.~5.5]{APP15}, \cite[Cor.~3]{IP}, \cite[(5)]{AIZ}
 \beq \label{sevruin} \Rui_{\q,\th}^b(x):=\E_x\left[e^{- \q \tz+ \th X_{\tz}} \1_{\{\tz<\tb\}} \right]=Z_{\q}(x,\th) -   \fr { W_{\q}(x)}{  W_{\q}(b)}Z_{\q}(b,\th), \; x \geq 0.
\eeq

   B) The joint Laplace transform of the first passage time at $0$ (``reflected ruin time'', see \eqref{Trb}) and the undershoot  in the presence of reflection at a barrier $b \geq 0$ is
  \cite[Prop.~5.5]{APP15}, \cite[Thm.~6]{IP}
\begin{equation}\label{sevruinrefl}
\Rui_{\q,\th}^{b]}(x):=\Eb_x\left[e^{- \q \tzb + \th X_{\tzb}} \right]=Z_{\q}(x,\th) -   \fr { W_{\q}(x)}{ W'_\q(b)} Z'_\q(b,\th ), \; x \geq 0.
\end{equation}
\eeT

{\bf Proof sketch:}  A)  is a  consequence of the harmonicity/$\q$-martingale property of $ Z_{\q}(X_t,\th)$, and of the boundary condition it satisfies \eqr{hext}. Indeed, stopping the martingale $e^{- \q t}  Z_{\q}(X_t,\th)$ at
$\min(\tb,\tz)$ yields
\begin{align*}  Z_{\q}(x,\th)&=\E_x \left[e^{- \q \tb} Z_{q}(b,\th) \1_{\{ \tb <\tz \}}\right]+
\E_x \left[e^{- \q \tz}  Z_{\q}(X_{\tz},\th) \1_{\{ \tz <\tb \} }\right]\\
&=\fr { W_{\q}(x)}{  W_{\q}(b)} Z_{\q}(b,\th)+
\E_x \left[e^{- \q \tz + \th X_{\tz}} \1_{\{\tz <\tb \}}\right] =\fr { W_{\q}(x)}{  W_{\q}(b)} Z_{\q}(b,\th)+\Rui_{\q,\th}^b(x).
\end{align*}
Note also that using another  (less smooth)   harmonic function with the same boundary condition, necessarily of the form $Z_{\q}(x,\th)+ k  W_{\q}(x)$, $k \neq 0$ would not change anything, since
$ W_{\q}(x) $ would cancel in the final result.\fn[5]{A direct proof using the resolvent formula \eqref{res} and \eqref{resost}
is also possible.}

B) Conditioning at $\min(\tb,\tz)$ shows that $\Rui_\q^{b]}(x,\th)$ is also of the form $ Z_{\q}(x,\th)-k  W_{\q}(x)$.  To determine $k$, we need to use  either the (non-trivial) boundary condition $\le(\Rui^{b]}_{q,\th}\ri)'(b)=0$ or the final value
$$\Rui_{\q,\th}^{b]}(b)=Z_{\q}(b,\th) -   \fr { W_{\q}(b)}{ W'_\q(b)} Z'_\q(b,\th ).$$ The latter has been
established in the related drawdown  literature -- see \eqr{de} and Theorem \ref{l:dd} for a generalization and further references. \qed

 \beC  \la{r:asycreep} \begin{enumerate}[A)] \im By using $ \lim_{b\to\infty}\frac{Z_{\q}(b,\th)}{W_{\q}(b)}=
 \frac{\k(\th)-\q} {\th-\Phi(\q)}$ (see \eqr{ZtoW} below) in \eqref{sevruin}, we recover \cite[(7)]{AIZ}
\beq \la{nob} \E_x \le[e^{-\q \tz +\th X_{\tz}} \; \1_{\{\tz < \I\}} \ri]= Z_{\q}(x,\th) -  W_{\q}(x)  \fr {\k(\th)-\q}{ \th-\Phi(\q)}, \quad \th > \F (\q). \eeq

\im
The relation \eqref{nob} holds as well for $\th=0,$ by analytic continuation, recovering the classic ruin time transform \cite[(10)]{AKP}
\beq \E_x \le[e^{-\q \tz } \; \1_{\{\tz < \I\}} \ri]=  Z_{\q}(x) -  W_{\q}(x)  \fr {\q}{\Phi(\q)}. \la{sevruin0}\eeq

\im The limit of \eqref{nob} when $\th \to \I$, which is the second term in the asymptotic expansion \eqref{Zt}, is
\beq &&\lim_{\th \to\I}
 \E_x \le[e^{-\q \tz + \th X_{\tz}}\ri]
 =\lim_{\th \to\I} \le(Z_{\q}(x,\th)-  \fr{\k(\th)-\q}{\th-\Fq}  W_{\q}(x)\ri)=\no\\&&\E_x \le[e^{-\q \tz } ; X_{\tzh}=0\ri]=\fr{\s^2}2 \le( W_{\q}'(x)-\Fq  W_{\q}(x)\ri). \la{creep}
 \eeq
 The last equality is the so-called ``creeping law'' \cite[Cor. 2]{Pispot}, \cite[(2.30)]{KKR}.
\im
A similar result for the  {hitting}  time of $0$ (``recovery after  ruin'') may be obtained by letting first $\th \to \Fq$ in \eqref{nob}.

Indeed, using $\k'(\Fq)=\fr {1}{\Fqp}$ and \eqr{Zexp}, we find
\begin{equation*} \E_x \le[e^{-\q \tz +\Fq X_{\tz}} \; \1_{\{\tz  < \I\}} \ri]=  Z_{\q}(x,\F (\q)) - \fr {1}{\Fqp} W_{\q}(x)  =e^{ x \Fq} - \fr {1}{\Fqp} W_{\q}(x). \la{fhit}\end{equation*}

Turning now to the \LT  of the  {hitting}  time of $0$, we find that for $x\geq 0,$
\beq && \E_x \le[e^{-\q \tzh } \; \1_{\{\tzh  < \I\}} \ri]=
\E_x \Big[\E_x \le[e^{-\q \tz -\q (\tzh - \tz)} \; \1_{\{\tz  < \I\}}| X_{\tz} \ri] \Big] = \\&&
  \E_x \le[e^{-\q \tz +\Fq X_{\tz}} \; \1_{\{\tz  < \I\}} \ri] =e^{ x \Fq} -  \fr {1}{\Fqp} W_{\q}(x),  \no \la{fhit1}\eeq  (alternatively, this formula may be obtained by a martingale stopping argument,   and holds for $ x \in \R $ as well).
This yields the representation of $ W_{\q}$ announced in \eqr{Wu}:
\beq \la{Wsim}  e^{ x \Fq}-\fr{ W_{\q}(x)}{\Fqp} =   \P_x \le[{\tzh  < \kil_\q} \ri].\eeq

 \end{enumerate}
\eeC

\ssec{From drawdowns to the   dividends-penalty  law}

This section and the following ones will exploit the connection between drawdown s  and  dividends. Namely,  the law  of the drawdown  triple and that of the dividend triple
\beq
\Big(\t_b , \ovl X_{\t_b} -X_0,  Y_{\t_b}-b, \Big)|\{ X_0=b\}, \quad  \Big(\tzb, \U_{\tzb}, - X_{T_0^{b]}})\Big)
  \la{dddiv}
\eeq
 coincide. See Figure \ref{fig:triples} below, where the paths of the process $X^{b]}$  are obtained from the paths of the process $X$ on the right by Skorokhod reflection  at $b$. For the picture of $X$, we may assume that $X_0=b$ for simplicity, but that is not necessary.    Now note that: a) the times $\tzb$ and $\t_b$ coincide; b) the total regulation equals the sum of the projections
 on the $X$ axis of the segments when $X$ is at a running maximum; c) the last drop must be the same on both pictures, since no reflection occurs during the last drop.  Thus $b- X_{\tzb}= Y_{\t_b}$.\fn[3]{To understand Skorokhod reflection informally, imagine the process $X$ arrives to $b$ from below, and encounters a barrier. If the barrier is fixed, it is forced to stick to the barrier until the first impulse downwards. If the barrier is movable, it is just raised during running maximum periods. In physics, under these two hypotheses, $b- X_t$  represents the distance to $b$  with respect to a fixed and  moving frame, respectively.}

 \begin{figure}[h]
 \begin{center}
 \includegraphics[width=0.48\textwidth]{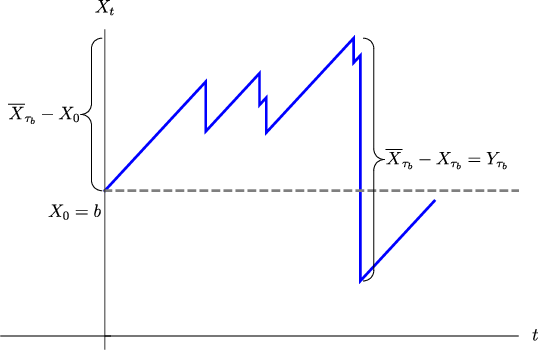}\quad
 \includegraphics[width=0.48\textwidth]{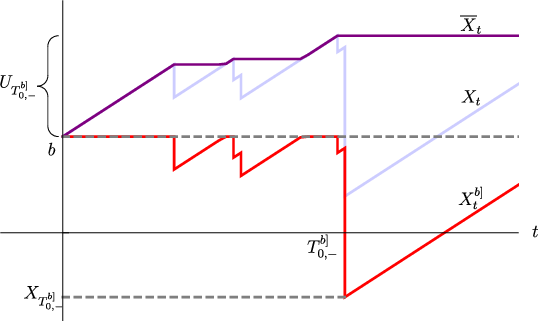}
 \end{center}
 \caption{Drawdown and dividend triples \eqref{dddiv}}
 \label{fig:triples}
 \end{figure}

This section  reviews first the independence of the law of the supremum  $\ovl X_{\td_d}-x$  of the law of  the (killed) drawdown
achieved on the last downwards excursion.
The former law is exponential with parameter $ \nu_\q(d)=\fr{ W_{\q}'(d)}{ W_{\q}(d)}$ (recall this follows intuitively from the fact that the upward ladder process with downward  excursions excised is
a drift killed at rate $\nu_\q(d)$). The independence is due intuitively to the fact that
each time the  upward ladder process  reaches a new point, the search for the killing excursion larger than $d$
starts again.

Equivalently, by \eqr{dddiv}, the independence of the dividends until ruin and of the final deficit
when starting from $b$ follows. When starting from $x<b$, one gets the famous Dividends-Penalty identity
first obtained in \cite{lin2003classical}.
\beT \la{l:dd}
 {\bf The deficit at   drawdown} \cite{MP}, \cite[Thm.~3.1]{LLZ17}, \cite[Prop.~3.1, 3.2]{LVZ}\fn[5]{We have re-expressed the result using the transformations in Remark \ref{r:tr}.} satisfies:
\begin{eqnarray} \la{funL} && {\de}_{\q,\th}(d,x,s):=\E_x \left[
e^{ -\q \td_d - \th (Y_{\td_d}-d)}; \ovl X_{\td_d} \in d s
\right]=
\Big(  \nu_\q(d)\; e^{-  \nu_\q(d) (s-x)_+}\; ds \Big)  \T \de_{q,\th}(d) \\&& \Eq
\E_x \left[
e^{ -\q \td_d - \th (Y_{\td_d}-d)- \vt (\ovl X_{\td_d} -x)}
\right]=
\fr{  \nu_\q(d)}{\vt+ \nu_\q(d) } \T \de_{q,\th}(d),  \no
\end{eqnarray}
where \be\la{de} \T \de_{q,\th}(d)=E_x \left[
e^{ -\q \td_d - \th (Y_{\td_d}-d)}\right]=Z_{\q}(d,\th)-W_\q(d)\fr{Z_{\q}'(d,\th)}{W_\q'(d)}.\end{equation}
\eeT

Using now the alternative interpretation furnished by \eqr{dddiv}
 yields  a powerful generalization of the deficit at ruin with reflection, Theorem \ref{l:s} B):
\beT \la{l:DP}   Let
\begin{equation*}
 DP_{\q,\th,\vt}^{b]}(x):=\Eb_x\left[e^{ -\q \tzb  + \th X_{\tzb}-\vt  \U_{\tzb}} \right]
\end{equation*}
denote the {\bf  dividends-penalty Laplace transform}\fn[4]{On an arbitrary interval
$|a,b]$, we will use the notation $\Rui_{\q,\th,\vt}^{b]}(x,a)$}.

A) When $x=b$, \ith
\begin{equation} \la{DPb}
 DP_{\q,\th,\vt}^{b]}(b)=
\fr{  \nu_\q(b)}{\vt+ \nu_\q(b) } \T \de_{\q,\th}(b).\end{equation} {Thus, when starting from $x=b$, the dividends $ \U_{\tzb \wedge \kil_\q}$  and the deficit at ruin $X_{\tzb \wedge \kil_\q}$ are independent, with the first  variable
having an exponential distribution
 \cite{Kyp}.}

B) Furthermore
\cite[Thm.~6]{IP}:
\begin{eqnarray} \la{DP}
 && DP_{\q,\th,\vt}^{b]}(x)=
Z_{\q}(x,\th) -  W_{\q}(x)  H_{DP}(b),
\\&& \la{HDP}
H_{DP}(b)=  \fr{  Z_{\q}'(b,\th)+ \vartheta Z_{\q}(b,\th)  }{ W_{\q}'(b)+ \vartheta  W_{\q}(b)}.
\end{eqnarray}

\eeT

{\bf Proof}: A)  When  starting at $x=b$ one may apply  Theorem \ref{l:dd} from the drawdown  literature.\fn[5]{Putting $l(x):=DP_{\q,\th,\vt}^{b]}(x)$, the (mixed) boundary condition at $x=b$  is now  $l'(b)+ \vartheta l(b)=0$; this offers another line of attack, at least in the \CL case.}

 B) Stopping at
$\tz \wedge \tb$ yields that $l(x)$ satisfies:
\begin{equation*}
 l(x)=Z_{\q}(x,\th) -     \fr{ W_{\q}(x)   }{  W_{\q}(b)} Z_{\q}(b,\th)+ \fr{  W_{\q}(x)   }{ W_{\q}(b)} l(b)=Z_{\q}(x,\th)+  W_{\q}(x)\fr{ l(b)- Z_{\q}(b,\th)   }{  W_{\q}(b)}
\end{equation*}
and the result follows from part A) by easy algebra.
 \qed

\beR It is easy to check that when $x=b$,  the transform \eqref{DP} factorizes and we recover \eqr{DPb}:

 \begin{align*}
 DP_{\q,\th,\vt}^{b]}(b)&=\fr{Z_{\q}(b,\th){ W_{\q}'(b)}- Z_{\q}'(b,\th){ W_{\q}(b)}} {{ W_{\q}'(b)}+\vt { W_{\q}(b)}}\\
 &=\fr{\fr{ W_{\q}'(b)}{ W_{\q}(b)}}  {\fr{ W_{\q}'(b)}{ W_{\q}(b)}+\vt} \le(Z_{\q}(b,\th)-Z_{\q}'(b,\th)\fr{ W_{\q}(b)}{ W_{\q}'(b)}\ri)=
  \fr{\fr{ W_{\q}'(b)}{ W_{\q}(b)}}
 {\fr{ W_{\q}'(b)}{ W_{\q}(b)}+\vt}
  \T \de_{\q,\th}(b).
\end{align*}
\eeR

\beR \la{r:DZW}Setting $\vt=0$ in $DP_{\q,\th,\vt}^{b]}(b)$
yields \begin{equation} \la{DelZW}
DP_{\q,\th,0}^{b]}(b)={ \de}_{\q,\th}(b)=\Rui_{\q,\th}^{b]}(b)=\fr{\D_{\q,\th}^{(ZW)}(b)}{ W_{\q}'(b)} 
\in (0,1),
\end{equation}
where we denoted
\begin{equation}\label{deltaZW}
\D_{\q,\th}^{(ZW)}(x,b):=Z_{\q}(x,\th)  W_{\q}'(b)-Z_{\q}'(b,\th) W_{\q}(x), \; \D_{\q,\th}^{(ZW)}(b)=\D_{\q,\th}^{(ZW)}(b,b).
\end{equation}

 The obvious nonnegativity of $\D_{\q,\th}^{(ZW)}$ implies that the function
$\fr{Z_{\q}(b,\th)}{  W_{\q}(b)}$ is decreasing (other papers refer to this as the log-convexity of  $Z_\q(x)$). It also
implies an
  upper bound for the  Wronskian
  $$0\leq \D^{(\ovl W W)}:= W_\q^2 (b)-\ovl W_\q (b) W_\q (b)) \leq \q^{-1} W_\q' (b).$$ \eeR
  The \nne \ of the  Wronskian

\ssec{From bailouts to the joint dividends-bailouts law}
After dividends, we now turn to bailouts as defined by $\L_t=-\min (\und X_t, 0)$,  and finally to their joint law.

\beT \la{l:be} {\bf Bailouts  until an exponential time}.
\bea &\text{A)   }& \E_x^{[0} \left[e^{- \th \L_{e_\q} }; e_\q < \tbz
\right]= 1- Z_{\q}(x)- Z_{\q}(x,\th)\fr{1- Z_{\q}(b)}{ Z_{\q}(b,\th)}\\
& \text{B)   }& \E_x^{[0} \left[e^{- \th \L_{e_\q \wedge \tbz}}
\right] = 1- Z_{\q}(x)+ Z_{\q}(x,\th)\fr{ Z_{\q}(b)}{ Z_{\q}(b,\th)}\\
&\text{C)   }& \E_x^{[0,b]} \left[e^{- \th \L_{e_\q} } \right]= 1- Z_{\q}(x)+ Z_{\q}(x,\th)\fr{ Z_{\q}'(b)}{ Z_{\q}'(b,\th)}. \eea\eeT

{\bf Proof}. A) Decompose $l(x):=\E_x^{[0} \left[e^{- \th \L_{e_\q} }; e_\q < \tbz\right]$ as
\begin{align*}
l(x)&=\E_x^{[0} \left[e^{- \th \L_{e_\q} }; e_\q < \tz \wedge \tb\right] +\E_x^{[0} \left[e^{- \th \L_{e_\q} }; \tz \leq  e_\q < \tb)\right]\\
&=P_x^{[0} \left[ e_\q < \tz \wedge \tb\right] + \E_x^{[0} \left[e^{ \th X_{\tz} }; \tz \leq  e_\q \wedge \tb\right] \; \E_0^{[0} \left[e^{ -\th \L_{e_\q} };   e_\q < \tbz\right]\\
&=\left(1- Z_{\q}(x)+ \fr{ W_{\q}(x)}{ W_{\q}(b)}\le({ Z_{\q}(b)}-1\ri)\right) + \le( Z_{\q}(x,\th)-\fr{ W_{\q}(x)}{ W_{\q}(b)}{ Z_{\q}(b,\th)}\ri) l(0),
\end{align*}
where we used the minimum law \eqref{tau} and the deficit law \eqref{sevruin}. In the \CL case when $ W_{\q}(0)\neq 0$  we may plug $x=0$ and conclude that
$$l(0)=\fr{\q \ovl W(b)}{ Z_{\q}(b,\th)}.$$
The same may be shown in the general case by a perturbation argument. Plugging now $l(0)$ yields the result A).

B) follows by adding \eqref{refbailout}.

C) follows by conditioning at time $e_\q \wedge \tbz$, where $h(x):=\E_x^{[0,b]} \left[e^{- \th \L_{e_\q} } \right]$. Indeed,
\bea h(x) = \Big(1- Z_{\q}(x)+ Z_{\q}(x,\th)\fr{ Z_{\q}(b)-1}{ Z_{\q}(b,\th)}\Big)+\fr{ Z_{\q}(x,\th)}{ Z_{\q}(b,\th)} h(b) \eea
\bea \Lra \fr{h(x) + Z_{\q}(x)-1}{ Z_{\q}(x,\th)}=\fr{h(b) + Z_{\q}(b)-1}{ Z_{\q}(b,\th)}=\fr{ Z_{\q}'(b)}{ Z_{\q}'(b,\th)}, \, \for x, \eea
where for the last equality we have used $h'(b)=0$ and the fact that for two functions $f$ and $g$, $f(x)/g(x)=c$ implies $f'(x)/g'(x)=c$. \hfill\qed

\beR  By letting $b \to \I$ in B) we recover \cite[Lem.~3.1]{AIjoint}. \eeR

\beT \la{l:joiDB} {\bf The joint dividends-bailouts law for a  process doubly reflected at $0$ and $b$, over an exponential horizon}.

The dividends-bailouts function is given by
\beq \la{DB}
 &&DB_{\q}^{[0,b]}(x,\th,\vt):=\Ezb_x\left[e^{  -\vt \U_{e_\q}- \th \L_{e_\q} }\right]= {1-Z_{\q}(x)}+
Z_{\q}(x,\th) DB_{\q}^{[0,b]}(0,\th,\vt),  \\&& \no DB_{\q}^{[0,b]}(0,\th,\vt)= \fr{ Z_{\q}'(b)+ \vartheta (Z_{\q}(b)-1) }
{Z_{\q}'(b,\th)+ \vartheta Z_{\q}(b,\th))}= \q \fr{W_{\q}(b)+ \vartheta \ovl W_{\q}(b) }
{Z_{\q}'(b,\th)+ \vartheta Z_{\q}(b,\th))}:=H_{DB}(b).
\eeq
\eeT

{\bf Proof}:
  Conditioning at $e_\q \wedge \tbz$ and using Theorem \ref{l:be} A) and Theorem \ref{l:refl}  we find
\begin{align*} l(x)&:=DB_{\q}^{[0,b]}(x,\th,\vt)=\Ezb_x\left[e^{  - \th \L_{e_\q} }; e_\q < \tbz\right] +\Ezb_x\left[e^{  - \th \L_{\tbz} };   \tbz \leq e_\q\right] l(b)\\&=1- Z_{\q}(x)- \fr{Z_{\q}(x,\th)}{Z_{\q}(b,\th)}(1- Z_{\q}(b))+
\fr{Z_{\q}(x,\th)}{Z_{\q}(b,\th)} l(b)\\&\Lra \fr{l(x)-1+ Z_{\q}(x)}{Z_{\q}(x,\th)}= \fr{l(b)-1+ Z_{\q}(b)}{Z_{\q}(b,\th)}=l(0).\end{align*}

The value of  $l(b)$
\be  l(b)=\Ezb_b \left[e^{  -\vt \U_{e_\q}- \th \L_{e_\q} }\right]=
\fr {Z_{\q}'(b,\th) + \Big(Z_{\q}(b,\th)   Z_{\q}'(b)-Z_{\q}'(b,\th)   Z_{\q}(b)\Big)}
{Z_{\q}'(b,\th)+ \vt Z_{\q}(b,\th)} \la{ZpZ}\ee
was obtained in \cite[Thm.~1]{AIjoint},
via  \exc \ theoretic arguments.
 \hfill\qed

\beR \la{r:ZpZ}  When $\th=0$, \eqref{ZpZ}  shows that discounted dividends starting from $b$ over an exponential horizon, with double reflection,  have an exponential law with parameter $\fr { Z_{\q}'(b)}{ Z_{\q}(b)}$, a surprising result which seems to have gone unnoticed.
Also,
 $ \Ezb_x [e^{-\vt \U_{e_\q}}]= 1-\fr{\vt Z_{\q}(x)}
 {Z_{\q}'(b)+ \vt Z_{\q}(b)}$,  recovering
 $ \Ezb_x [\U_{e_\q}]=\fr{Z_{\q}(x)}{Z_{\q}'(b)}$ \cite[(4.3)]{APP}.

 Putting $\vt=0$ in \eqref{DB} yields Theorem \ref{l:be} C), and differentiating recovers \cite[(4.4)]{APP}
\begin{align*} \Ezb_x [\L_{e_\q}] &= \fr{1}{
Z_{\q}'(b)}\Big[ Z_{\q}(x)\Big(Z_{\q}(b)-\k'(0_+)  W_{\q}(b)\Big)-\Big(\ovl   Z_{\q}(x) -\k'(0_+) \ovl   W_{\q}(x) \Big) \q W_{\q}(b)\Big]\\
&=\fr{ Z_{\q}(x) Z_{\q}(b)-\ovl   Z_{\q}(x)   Z_{\q}'(b)-\k'(0_+)  W_{\q}(b)}{
Z_{\q}'(b)}=\fr{ Z_{\q}(x) Z_{\q}(b)    }{
Z_{\q}'(b)}-\ovl   Z_{\q}(x)-\fr{\k'(0_+)}{\q},\end{align*}
where $\ovl   Z_{\q}(x)$ is defined in \eqref{ZbarD}.

\eeR

\ssec{Expected discounted bailouts}
We recall now results on  expected discounted bailouts until $\tb$ and over an infinite horizon, which may be obtained simply by differentiating  the corresponding moment generating functions in Theorem \ref{l:be} B), C).
\beT \la{l:VL} Put  \begin{equation} \la{GB} G_{\q}^B(x)=Z_{\q}^{(1)}(x)=\fr{\partial Z_{\q}(x,\th)}{\partial \th}_{\th=0}=\ovl  Z_{\q}(x) -\k'(0_+) \ovl   W_{\q}(x).\end{equation}
A) The {\bf expectation of the total discounted bailouts  up to $\tb$} for $0\leq x\leq b$ is \cite[Cor.~3.2 (ii)]{APY}:
\begin{equation} \label{VL}
	B^{b}(x):=\Ez_x\left[\int_0^{T_b^{[0}} e^{-\q t}\diff  \L_t\right]=
\Ez_x\left[ \L_{T_b^{[0} \wedge e_\q} \right]=\frac{Z_{\q}( x)}{Z_{\q}( b)} G_{\q}^B(b)-G_{\q}^B(x).
\end{equation}

B) The {\bf expected total discounted bailouts  over an infinite horizon}, with reflection at $b$ are \cite[(4.4)]{APP}:
\begin{equation}\label{VLSLG}
 B^{[0,b]}(x)=\Ezb_x\le[\int_0^{\I} e^{-\q t}\diff  \L_t\ri] =\Ezb_x\le[ \L_{e_\q}\ri]= \frac{ Z_{\q}( x)}{Z'_{\q}( b)} (G_{\q}^B)'(b)-G_{\q}^B(x).
\end{equation}
 $G_{\q}^B$ may also be taken to be   \beq \la{G} G_{\q}^B(x)= \ovl{ Z}_\q(x)+\frac{\k'(0_+)}{\q},\eeq
in both results.

\eeT

\beR \la{GSnu} As may be easily checked, the first expression for  $G_{\q}^B$, i.e. $Z_{\q}^{(1)}(x)$, is the smooth \GS function (see \cite{APP15} and next section), fitting  the value of $w(x)=x$ at $0$, and also its derivative in the non-\cP case. Without smoothness, the \GS function is  unique only up to adding a multiple of the \corr scale function, and   simpler expressions  like \eqr{G} may be available.\eeR

\beR \la{r:W2z}
 Note that  several relations for the  process  reflected below like  \eqref{VLSLG}, and  the relation $ \E_x^{[0} [e^{- \q \tb}]= \frac{Z_{\q}(x)}{Z_{\q}(b)}$ \cite{AKP} 
 may be obtained formally from analog relations  for the process absorbed at $0$,  by substituting the second scale function $Z_{\q}$ instead of the first scale function $ W_{\q}$. \eeR

\ssec{Results obtained by differentiating the moment generating functions \la{s:diff}}
We turn now to obtain the expectations of the ruin time, exit time from an interval, reflected ruin time, reflected up time and  recovery after  ruin time, obtained by differentiating  the  respective moment generating functions \eqref{sevruin0}, \eqref{tau}, \eqref{sevruinrefl}, \eqr{refbailout}, \eqref{fhit} with respect to $\q$
(making  use of the  analyticity of $W_\q$ in $\q$ \cite[Lem.~8.3]{Kyp}), and  putting $\q=0$. In the proof of B) below, we additionally use  the fact that
when some function $f$ is differentiable at 0, it holds that
$\fr{\partial \big[\q f(\q)\big]}{\partial \q}_{\q= 0}= f(0)$.

\beT \la{p:time} A) When $\k'(0_+) <  0 \Lra \F(0) >0$, it holds that \bea &&\E_x
\le[ \tz   \ri]=\frac{W(x)}{\F(0)}- \ovl W(x).\eea

 When $\k'(0_+) >  0 \Lra \F(0) =0$, it holds that
\begin{align*} \E_x
\le[ \tz  \; \1_{\{\tz < \I\}} \ri] &= W(x) \lim_{\q \to 0} \fr{\Fq -\q \F'(q)}{\Fq^2}+ \k'(0_+) W^{*2}(x)- \ovl W(x) \\
&=-\k'(0_+)^2 \fr{\F''(0_+)}2  W(x)  + \k'(0_+) W^{*2}(x) - \ovl W(x)\\
&= \fr{\k''(0_+)}{2\k'(0_+)}  W(x) + \k'(0_+) W^{*2}(x)- \ovl W(x),
\end{align*}
where we used
\beq \F''(x)=-\fr{\k''(x)}{(\k'(x))^3}\eeq
and the series expansion \cite[(8.29)]{Kyp}
\beq \la{qexp}  W_{\q}(x)=\sum_{k=0}^\I \q^k W^{*,k+1}(x),\eeq
with $W^{*,k}(x)$ denoting convolution.

B) Put $T=\tz \wedge \tb$. Then\fn[4]{This provides  a third proof of the monotonicity of
   $\fr{\ovl W(b)}{ W(b)} $ (see Remark \ref{r:tmin}).}
\beq \E_x [T]=\fr{W(x)}{ W(b)} \ovl W(b)-\ovl W(x).\eeq

C) $$\Eb_x
\le[\tzb   \ri]= W(x) \fr{W(b)}{W'(b)}- \ovl W(x) \Lra \E
\le[\t_{b}   \ri]=  \fr{W(b)^2}{W'(b)}- \ovl W(b).$$

D) $$E
\le[T^{[0}_{b}   \ri]=   \ovl W(b).$$

E)  \beq &&\E_x
\le[ \tzh ; \tzh < \I   \ri]=\k'(\F(0)) W^{*,2}(x) + \fr{\k''(\F(0))}{\k'(\F(0))}  W(x)- \fr { x e^{ x \,\Phi(0) } }{\k'(\F(0))}.\eeq

  When $\k'(0_+) >  0 \Lra \F(0) =0$, this simplifies to
\beq &&\E_x
\le[ \tzh ; \tzh < \I   \ri]=\k'(0_+) W^{*,2}(x) + \fr{\k''(0_+)}{\k'(0_+)}  W(x)- \fr x {\k'(0_+)}.\eeq

\eeT

\beR In the particular \cP case, A) reduces, using $W(x)=\fr{\sRui(x)}{\k'(0_+)}$ and $\k''(0)=\lambda \E[C_i^2]$ to \cite[(11.3.26)]{RSS}
\bea &&\E_x
\le[ \tz  \; \1_{\{\tz < \I\}} \ri]= \fr{\k''(0)}{2 \k'(0_+)^2}  \sRui(x) - \frac{1}{\k'(0_+)} \int_0^x \sRui(y) \Rui(x-y) d y.\eea
Our examples show that the expected time to ruin conditioning on ruin happening is unimodular, with a unique maximum. This maximum could be viewed as a  reasonable lower bound for the initial
reserve, which postpones ruin as much as possible (in the worst case).
\eeR
\beR

To show  the \nny of C),  it suffices to take $x=b$, where the \nny holds by the log-concavity of $\overline{W}^{(q)}$, proved in Remark \ref{r:tmin}.

When $b \to \I$ and $\k'(0_+) <  0$, C) converges to A).

 When $x=0$, C) yields  the ``0-cycle  law'' \cite[Prop.~3.2(i)]{starreveld2016occupation}
\beq \Eb_0
\le[\tzb  \ri]= W(0) \fr{W(b)}{W'(b)}.\eeq
\eeR

To give an idea of very recent developments in the $W,Z$ theory, we end this section with a  hitting time result which holds
for certain Omega spectrally negative Markov processes  as well \cite[Cor.~1]{LZ17} (the proof is quite elegant).
\beT For $x,i \in (a,b)$,  it holds that
$$\E_x\le[e^{- \int_0^{T_{\{i\}}} \q  d s}; T_{\{i\}} \leq \ta \wedge \tb\ri]=\fr{ W_{\q}(x-a)}{ W_{\q}(i-a)}- \fr{ W_{\q}(x-i)}{ W_{\q}(b-i)}\fr{ W_{\q}(b-a)}{ W_{\q}(i-a)}.$$ \eeT
For the general  result with Omega non-constant killing, it suffices to replace $\int_0^{T_{\{i\}}} \q  d s$ by $\int_0^{T_{\{i\}}} \w(X_s)  d s,$  where  $\omega:\mathbb{R}\to\mathbb{R}_{+}$ is an arbitrary locally bounded nonnegative measurable state dependent discounting, to replace $b-a$ by $b,a$,..., etc.,   and to identify the scale function $W_\w$ \cite{LP,LZ17} -- see also Section \ref{s:ome}.

\sec{Smooth Gerber-Shiu functions: $Z_{\q}(x,\th)$  is replaced   by the smooth Gerber-Shiu function $G_{w}(x)$ \la{s:nonh}}
    When $e^{\th X_{\tz}}$ is replaced  in the previous formulas  \eqref{sevruin}, \eqref{sevruinrefl} by an arbitrary penalty function $w( X_{\tz}), w:(-\infty,0]\to\mbb R$, extensions of these formulas still hold for $$\mc V^b (x):=\E_x\le[e^{-\q\tz}w(X_{\tz}) \1_{\{\tz<\tb\}}\ri],$$ if one   replaces $Z_{\q}(x,\th)$
by an {\em infinite horizon Gerber-Shiu penalty function}
$$\mc V (x): = \E_x\le[e^{-\q\tz}w(X_{\tz})\ri].
$$
Indeed, applying the strong Markov property at $\tb$ immediately yields
$$\mc V (x)=\mc V^b (x)  + \fr { W_{\q}(x)}{  W_{\q}(b)}\mc V (b) \Lra \mc V^b (x) =\mc V (x)- \fr { W_{\q}(x)}{  W_{\q}(b)}\mc V (b).$$
Note that $\mc V (x)$ is not unique: it  may be replaced in the identity above by adding to it any multiple of $ W_{\q}(x)$ \cite[Prop.~5.4]{APP15}. 

For this reason, \cite[Thm.~5.3]{APP15} identify the unique ``smooth Gerber-Shiu function'' ${G} $ \cite[Def.~5.2]{APP15},  which exists if $w$ satisfies some minimal  integrability conditions. Under these,
given  $0<b<\infty$,  $x\in(0,b)$,
there exists a unique smooth  function ${G} =G_\q$ so that the following hold:
\beq&&
\mc V^b (x)=\E_x\le[e^{-q
\tz}w\le(X_{\tz}\ri)\mbf 1_{\{\tz<\tb\}}\ri]  = {G} (x ) - \frac{ W_{\q}(x ) }{ W_{\q}(b )} {G} (b ), \la{serui} \\&&
\mc V^{b]} (x )=\Eb_x\le[e^{-q
\tzb}w\le(X_{\tzb}\ri)\ri]  = {G} (x )-
 \frac{ W_{\q}(x ) }{ W_{\q}'(b )} {G} '(b ).\la{rserui}
\eeq
Stated informally,  both   problems above admit decompositions involving the same ``non-homogeneous  solution'' ${G} $.

The  ``smoothness'' required is:
 \begin{eqnarray}\label{dFw0}
\begin{cases}
{G} (0) = w(0),\\
{G} '(0_+) = w'(0_-), & \text{in the case $\sigma^2>0$ or $\lm([0,1])=\infty$}.
\end{cases}
\end{eqnarray}
Under these conditions, the function ${G} $ is unique. Furthermore, it  may be represented as \cite[(5.13) Lem.~5.6]{APP15}:
\begin{align}
{G} (x)&=   w(0)  Z_{\q}(x) +  w'(0-) \fr{\s^2}2    W_{\q}(x) +  \int_0^x  W_{\q}(x-y) \int_{z=y}^{\I}  [w(0)-w(y-z)] \lm(dz)  dy \no\\
&=  w(0) \le(\fr{\s^2}2   W'_\q(x) + c  W_{\q}(x)\ri)+ w'(0-) \fr{\s^2}2    W_{\q}(x) -   \int_0^{x}  W_{\q}(x-y)  w^{(\lm)}(y) dy \la{Gdec} ,
\end{align}
where $w^{(\lm)}(y)=\int_{z=y}^\I  [w(y-z)] \lm(dz)$ is the
 expected liquidation cost conditioned on a pre-ruin position of $y$, with ruin causing jump bigger than $y$. The second equality follows by using \eqref{oldZ}.

 \beR The last term in the second equality in \eqref{Gdec} fits the  ``non-local'' part of $w$, and
 the first two terms  may be viewed as boundary  fitting terms. Indeed, this holds  since
 $\fr{\s^2}2   W'_\q(0_+) + c  W_{\q}(0_+) =1,   \fr{\s^2}2    W_{\q}(0_+)=0$, and $\fr{\s^2}2   W''_\q(0_+) + c W'_\q(0_+) =0,   \fr{\s^2}2  W'_\q(0_+)=1$.  \eeR

\beP \la{p:GSexp} For  $w(x)=e^{\th x}$, the Gerber-Shiu function is $Z_{\q}(x,\th)$ and the decomposition \eqref{Gdec} becomes:
\bea Z_{\q}(x,\th)=  Z_{\q}(x) + \th \fr{\s^2}2  W_{\q}(x)+\int_0^x  W_{\q}(y) \int_{x-y}^{\I}  [1-e^{\th(x-y-z)}] \lm(dz)  dy.\eea
\eeP

This may be easily checked by taking Laplace transforms, since
$$\H  W_{\q}(s)\fr{\k(s)-\k(\th)}{s-\th}=\H  W_{\q}(s)\Big(\fr{\k(s)}{s}+ \th \fr{\s^2}2+ \fr{\H \pi(s)-\H \pi(\th)}{s-\th}- \fr{\H \pi(s)-\H \pi(0)}{s}\Big).$$

\sec{Poissonian/Parisian
detection of bankruptcy/insolvency, and occupation times \la{s:Par}}

A useful type of models developed recently \cite{AIZ,albrecher2017strikingly,APY} assume  that insolvency is only {\bf  observed periodically},
  at an  increasing sequence of \emph{Poisson observation times} $\mT_\r =\{t_i,i=1,2,...\}$, the arrival times   of an independent Poisson process of rate $\r $, with  $\r  > 0$ fixed\fn[4]{The concept of periodic observation may be   extended to the  Sparre Andersen  (non L\'evy) case,  using geometrically distributed  intervention times at the times of claims. This  deserves further investigation.}.
  The analog concepts for first passage times are the stopping times
\begin{eqnarray}T_{b,+}=T_{b,+}^\r = \inf\{t_i:\; X_{t_i}  > b\}, \quad T_{a,-} =T_{a,-}^\r  = \inf\{t_i > 0:\; X_{t_i} < {a}\} \label{def_tau_a_plus_minus}
 \end{eqnarray}
  Under
 {Parisian observation times}, first passage is recorded
 only when the most recent excursion  below $a$/above $b$
 has exceeded an exponential random variable $e_\r$  of rate $\r $.
 We  use here the same notation as for classic first passage times (which correspond to the case
 $\r  =\I$).

\beR We will refer to stopping at $T_{0,-}$ as (exponential) \textbf{Parisian absorption}.
 A {\bf spectrally negative  \lev processes with (exponential) Parisian reflection} below $0$ may be defined  by   pushing the process up to $0$ each time  it is below $0$ at an observation time $t_i$. In both cases, this will not be made explicit in the notation; classic and Parisian absorbtion and reflection will be denoted in the same way.

 Note that the  case $\r  \to 0 $ corresponds to complete leniency; default
 is never observed. We see thus that Parisian inspection is an  intermediate situation  between continuous inspection and no inspection, and can help to render modelling more realistic.  \eeR

 It was recently observed   that  the classic first passage laws
listed above hold with a ``Parisianly observed" lower boundary,
once $ W_{\q},  Z_{\q}$ are replaced by  appropriate generalizations, defined by \cite{APY,AZ}:
\begin{eqnarray} \la{Z2}
 && Z_{\q,\r}(x,\th):=\fr{\r} {\q+\r- \k(\th)} Z_{\q}(x,\th)
 +\fr{\q-\k(\th)}{\q+\r- \k(\th)} Z_{\q}(x,\Fqr )\\&&=\fr{\r} {\q+\r- \k(\th)} (Z_{\q}(x,\th)-Z_{\q}(x,\Fqr))
 + Z_{\q}(x,\Fqr), \no\\ && {W_{\q,\r} (x)}:= \fr{ \Fqr -\Fq}{\r} Z_{\q}(x,\Fqr ),
 \end{eqnarray}
with  the value for $\th=\Fqr $ being interpreted in the limiting sense.\fn[4]{When $\r  \to \I$, the Parisian results
reduce to  the classic ones, since $Z_{\q,\r}(x,\th), W_{\q, \r }(x )$ are asymptotically equivalent to $ Z_{\q}(x,\th ),  W_{\q}(x)$. The first assertion is trivial, for the second  see \eqr{ZWF}. The notation $W_{\q,\r} (x):= \fr{ \Fqr -\Fq}{\r } Z_{\q}(x, \Fqr )$  has been chosen to emphasize that this replaces, for processes with Parisian ruin, the $W_\q$ scale function in  the classic ``gambler's winning''  problem, and also to ensure a convenient \asy \ behavior.}

\beR    Exponential Parisian detection below $0$ is related to  the \LT of the total ``occupation time spent in the red''
$$ T^{<0}:= \int_0^\I \1_{\{X_t <0\}} d t, $$
 a fundamental risk measure studied by
 {\cite{picard1994some,zhang2002total,loisel2005differentiation}.

  Indeed, the probability of Parisian ruin not being observed (and of recovering  without  bailout) when $\kappa'(0_+) >0$ is   \cite[Cor.~1,Thm.~1]{LaRZ}, \cite[(11)]{AIZ}
   \be  P_x[T_{0,-}=\I]=P_x[ T^{<0}< e_\r ]
  = \E_x \left[e^{- \r  T^{<0}}\right]= \kappa'(0_+) \fr{\F(\r )}  {\r }  Z(x,\F(\r )). \la{redext}\ee}
 When $x=0,$
 this reduces to
 \beq \la{SAid}  P_0[T_{0,-}=\I]=P_0[ T^{<0}< e_\r ]
  = \E_0 \left[e^{- \r  T^{<0}}\right]= \kappa'(0_+) \fr{\F(\r )}  {\r } ,\eeq
 a quantity which
 could be viewed as  a {\bf model dependent} extension of the profit parameter $\k'(0_+)$,  measuring the profitability of a risk process.

 Note  that $\kappa'(0_+) \fr{\F(\r )}  {\r } $ furnishes also the \LT of six other remarkable random variables besides $T^{<0}$, by the  ``\SA identities'' due to \cite[Prop.~1.1,(2)]{ivanovs2016sparre}. Differentiating \eqref{SAid} with respect to $\r $  when $\kappa'(0_+)>0$ shows that the Sparre-Andersen-Ivanovs variables  have all expectation $-\fr{\F''(0)}2=\fr{\k''(0)}{(\k'(0))^3}$, a quantity which appeared already in several previous computations.
\eeR

 The following proposition lists  some basic first passage results for processes with Parisian detection of ruin, reflected or absorbed, following \cite{AIZ,BPPR,APY}. Note that these results coincide with the ones with classic, ``hard'' detection of ruin, and imply them when $\r   \to \I$. 
\beT \la{twostep} {\bf First passage results for processes with classic detection at a smooth boundary $b >0$ and Parisian detection beolw $0$, followed by stopping or by reflection}.
Let $X$ be a spectrally negative \lev process with Parisian detection below $0$, and fix $b>0$.   Assuming $x\in[0,b]$ and $\q, \r >0,0\leq \theta < \I$, using the notation of Remark \ref{r:not} and letting  $W_{\q,\r }(x)$ and $Z_{\q,\r }(x,\th)$ be defined   by \eqref{Z2},  the  following hold:

\BEN 
\im The ``gamblers survival formula takes the form \cite[(12)]{AIZ}
\be \E_x^{\{0} \Big[e^{- \q \tb} \1_{\{\tb < T_{0,-}\}}\Big]=\frac{W_{\q,\r } (x)}{  W_{\q,\r } (b)} .\ee
\im A)  The {\bf expected discounted dividends (upper regulation at  $b$)
 until $\tzb$}    are \cite[(27)]{AIZ}:
\begin{equation} \la{divPar}
 V^{\vdots 0,b]}(x)=\Eb_x\le[\int_0^{T_{0,-}} e^{-\q  t}d  \U_t\ri] =   \frac{W_{\q,\r } (x)}{  W'_{\q,\r } (b)}=\frac{Z_{\q}(x, \Fqr )}{  Z_{\q}'(b, \Fqr )}
\end{equation}

B)  The {\bf expected discounted dividends
  with reflection at $0$ at Parisian times, until the total bail-outs  surpass an exponential variable $e_\th$}  \cite[(15)]{AIpower} are
 \begin{eqnarray}V_{\U}^{\{ 0,b]}(x,\th)=
 \E^{\{ 0,b]}_x\le[\int_0^\I e^{-\q s} \1_{ \{\L(s) < e_\theta\}} d  \U(s) \ri]=\frac{Z_{\q,\r}(x,\th)}
{Z_{\q,\r}'(b,\th)} \label{e:divskillbo}\end{eqnarray}

\beR When $\theta=0$, this becomes \cite[Cor.~3.3]{APY}:
\begin{eqnarray}\label{eq:DRLa2}
 V_{\U}^{\{0,b]}(x)=\E^{\{ 0,b]}_x\le[\int_0^{\I} e^{-\q  t}d  \U_t\ri] &=&    \frac{Z_{\q,\r } (x)}{  Z_{\q,\r}' (b)}.
\end{eqnarray}
\eeR

\im The {\bf   capital injections/bailouts law for a process with Parisian reflection at $0$, until $\tb$} \cite[Cor.~3.1 ii)]{APY}.
Let $L_t$ denote the regulator for the process with Parisian reflection at  $0$ and $\E^{\{ 0}_x$ the expectation for such process.
Then:
 {
\begin{equation} \la{Parisbailouts}  
\sRui_{\q,\th,\r} ^{\{ 0,b|}(x):=\E^{\{ 0}_x [e^{-\q \tb - \theta \L_{\tb}}]=
\bc \frac{Z_{\q,\r}(x,\th) }{Z_{\q,\r}(b,\th)}& \theta <\I\\ \E_x^{\{0}[e^{-q\tb}; \tb<T_{0,-}]=
\frac{W_{\q,\r } (x) }{W_{\q,\r } (b)}&  \theta =\I\ec.\end{equation}}

\im {\bf Deficit at ruin for a process absorbed or reflected at $b>0$}.

 A)
 The joint Laplace transform of the Parisian first passage time of $0$ and the undershoot for a process absorbed  at $\tb$
 is given by    \cite[(15)]{AIZ}:
 \fn[4]{the second expression in \eqr{sevruinP} uses a simpler, non-smooth \GS function --see Remark \eqr{GSnu}.}

\begin{equation} \label{sevruinP}
\Rui_{\q,\th,\r} ^{\vdots 0,b|}(x):=\E_x\left[e^{\theta X_{T_{0,-}} } \1_{\{T_{0,-} < \tb \wedge e_q\}}\right]=
Z_{\q,\r}(x,\th)-W_{\q,\r } (x)
{W_{\q,\r } (b)}^{-1}{Z_{\q,\r}(b,\th)}
\end{equation} \bea =\fr{\r  } {\q + \r   -\k(\th)}\le(Z_{\q}(x,\th)-W_{\q,\r } (x)
{W_{\q,\r } (b)}^{-1} Z_{\q}(b,\th)\ri)\eea

   B) The joint Laplace transform of the first passage time at $0$ and the undershoot  in the presence of reflection at a barrier $b \geq 0$ is
\begin{equation}\label{sevruinreflP}
\Rui_{\q,\th,\r} ^{\vdots 0,b]}(x):=\Eb_x\left[e^{- \q T_0^{b]} + \th X_{T_0^{b]}}} \right]=Z_{\q,\r}(x,\th) -   \fr {W_{\q,\r } (x)}{ W'_{\q,\r } (b)} Z_{\q,\r}'(b,\th), \; x \geq 0.
\end{equation}

\im
Let $U_{\q,\r } ^{\vdots a,b|}(x,B)=\E_x \left[ \int_0^{ \ta \wedge \tb } e^{-qt} \1_{\left\{X_t  \in B  \right\}} \diff  t \right],$  denote the {\bf $q$-resolvent of a doubly absorbed spectrally negative  L\'{e}vy process with Parisian ruin}, for any Borel set $B\subset [a,b]$. Then  \cite[Thm.~2]{BPPR}
\begin{align} \label{resolvent_density} 
\begin{split}
	U_{\q,\r } ^{\vdots a,b|}(x,B) &= \int_{a}^{b} \1_{\{y \in B\}}\left( \frac {W_{\q,\r } (x-a) W_{\q,\r }  (b-y)} {W_{\q,\r } (b-a)} -W_{\q,\r }  (x-y) \right) \diff  y, \quad a < x < b.
	\end{split}
\end{align}

\im
 The {\bf  dividends-penalty law for a  process reflected at $b$, with Parisian ruin} 
 is:
\begin{equation} \la{DP2}
 DP_{\q,\theta,\vt}^{\vdots 0,b]}(x):=\Eb_x\left[e^{  -\vt \U_{T_{0,-}}+ \theta X_{T_{0,-}}} ; T_{0,-} <e_q\right]=Z_{\q,\r}(b,\th) - W_{\q,\r } (b)   \fr{ Z_{\q,\r}'(b,\th)+ \vartheta Z_{\q,\r}(b,\th)  }{W_{\q,\r } '(b)+ \vartheta W_{\q,\r } (b)} \end{equation}
 \beq =
\left(Z_{\q}(x,\th) -Z_{\q,\F(\q+\r ) }(x) {H_{\F(\q+\r ) }(b)}^{-1}{H_{\theta}(b)}\right)\r  ({\q+\r- \k(\theta)})^{-1}, \la{mis}
\eeq
where $H_\th (b)= \vt Z_{\q}(b,\th)+ Z_{\q}'(b,\th)=(\theta+ \vt) Z_{\q}(b,\th)+(\q-\k(\theta)) W_\q(b).$\fn[5]{The structure of this formula reflects the fact that $\Fqr $ is a removable singularity.}  We included the second, rather  complicated    formula,  to allow comparison with the original formula in \cite[(23)]{AIZ}.

\beR
When $x=b$, we may factorize the transform  $\Eb_b\left[e^{ \theta X_{T_{0}} -\vt \U_{T_{0}}} ; T_{0} <e_q\right]$ \eqref{mis} as:
\begin{equation}
\fr{\nu_{\q,\r } }{\nu_{\q,\r } + \vt}
\Bigl(Z_{\q}(b,\th)-{\nu_{\q,\r } } ^{-1}\Big(\theta Z_{\q}(b,\th) +\le(q-\k(\theta)\ri) W_q(b)\Big)\Bigg) \fr{\r  } {\r +\q- \k(\theta)}, \la{fac}\end{equation}
where $ \nu_{\q,\r } =V^{b]}(b)^{-1}=
{W_{\q,\r } '(b)}{W_{\q,\r } (b)}^{-1}=
{Z_{\q,\F_{\q+\r } }'(b)}{Z_{\q,\F_{\q+\r } }(b)}^{-1}$. 
Indeed,
\bea  &&Z_{\q}(b,\th) -Z_{\q,\F_{\q+\r } }(b)
\le( (\F_{\q+\r } + \vt) Z_{\q,\F_{\q+\r } }(b)-\r   W_\q(b)\ri)^{-1} H_\th(b)
\\&&=Z_{\q}(b,\th) -\le( \vt + \F_{\q+\r }  - \r   W_\q(b) Z_{\q,\F_{\q+\r } }(b)^{-1}\ri)^{-1}H_\th(b)\\&&=Z_{\q}(b,\th) -\le( \vt + \nu_{\q,\r } \ri)^{-1} H_\th(b), \eea
 and
\eqref{fac} follows by simple algebra.
By \eqref{fac},  $ \U_{T_{0}}$ and $X_{T_{0}}$ are independent when starting from $b$, and the former
has an exponential distribution with parameter $\nu_{\q,\r } $ \cite[(23),(26)]{AIZ}.

 When $\vt=0$, this result reduces to \eqref{sevruinreflP}.
\eeR

\im
A) The {\bf expected total discounted bailouts at Parisian times up to $\tb$} are given for $0\leq x\leq b$ and $q>0$ by
 \cite[Cor.~3.2 ii)]{APY}:
\begin{align}\label{ParBail}
	B^{\{ 0,b|}(x):=\Ez_x&\left[\int_0^{\tb}e^{-qt}\diff  \L_t\right]=\frac{Z_{\q,\r } ( x)}{  Z_{\q,\r } ( b)} G_{\q,\r}^B(b)-G_{\q,\r}^B(x)
\end{align}
where \beq \la{ParBailG} G_{\q,\r}^B (x)
=\frac{\r  } {\q+\r }   \Big(\ovl{Z}_{q}(x)+\frac{\k'(0_+)}{\q}
\Big)
=\frac{\r  } {\q+\r }G_{\q}^B (x) .\eeq

B)  The {\bf expected total discounted bailouts at Parisian times over an infinite horizon}, with reflection at $b$ are \cite[Cor.~3.4]{APY} (see also \cite[Thm.~3.2]{zhao2017optimal}, where $Z_{\q,\r } ( x) $ is denoted by $B_2(x)$\fn[4]{Our sign of $\fr p q$ in formula \eqr{ParBailG} for $G_B(x)$ is opposite to that
  in formulas (3.26) and (3.30) of \cite{zhao2017optimal}, since they consider \sp \procs.}):
\begin{eqnarray}\label{eq:localt0}
V^{\{ 0,b]}(x)=\Ezb_x\le[\int_0^{\I} e^{-q t}\diff  \L_t\ri] &=& \frac{Z_{\q,\r } ( x) }{ Z_{\q,\r}' ( b)} (G_{\q,\r}^B)'(b)-G_{\q,\r}^B(x).
\end{eqnarray}

\EEN
\eeT
\beR
Note that each result from Theorem \ref{twostep} has its analog in classical detection of ruin. Indeed,
\begin{itemize}
\item (2) corresponds to the dividends Theorem \ref{l:div};
\item (3)
is the Parisian analog of the bail-outs Theorem \ref{l:refl} (\cite[Thm.~2]{IP})  ;
\item (4) A) and B) are Parisian analogues of Theorem \ref{l:s} A) and B) (\cite[Prop.~5.5]{APP15});
\item (5) corresponds to the resolvent formula \eqref{res};   it is natural to conjecture that the  resolvents for (partly) reflected processes will also be of  the same form as the classic ones \cite[Thm.~1]{PDR}, \cite[Thm.~2, Cor.~2]{Ivapot};
\item (6) is the Parisian analog of the dividends-penalty Theorem \ref{l:DP};

\item (7) corresponds to the expected total discounted bailouts Theorem \ref{l:VL}. One may check that
    \begin{eqnarray}\label{eq:localt1}
V_{\q,\th,\r}^{\{ 0,b]}(x)=\Ezb_x\le[\int_0^{\I} e^{-q t}\1_{ \{\L(s) < e_\theta\}} \diff  \L_t\ri] &=&
\frac{Z_{\q,\r}(x,\th) }{ Z_{\q,\r}'(b,\th)} (G_{\q,\r}^B)'(b)-G_{\q,\r}^B(x).
\end{eqnarray}
\end{itemize}

\eeR

\beQ
It is natural to conjecture that the outstanding  results which have not yet been extended from the classic to the Parisian case, like  Theorem \ref{l:joiDB} on the joint distribution of dividends and bailouts,   the  optimality of barrier policies with fixed  final penalty \eqr{Kb} \cite[Prop.~4.3]{hernandez2016time},  the  optimality of barrier policies for the \SLG objective \cite[Lem.~2]{APP},  etc, hold in the  Parisian case as well. \eeQ

\beQ
The fact that the results for the Parisian case coincide with the classical ones suggest that the known first passage results with hard ruin for SNMAPs \cite{KP,Iva,IP,AIrisk} might generalize to the Parisian case, provided that  properly defined scale matrix functions are introduced, and   multiplied in  correct order.   To facilitate further work, we provide non-Parisian SNMAP references for the corresponding results of Theorem \ref{twostep}:  for (2) A) and B) see \cite[Cor.~3]{IP} and \cite[Thm.~6]{IP} respectively; for (3) see \cite[Thm.~2]{IP}; for (4) see \cite[Thm.~2, Cor.~2]{Ivapot}; for (5) see \cite[Thm.~6]{IP}.

Most interesting is the problem of resolvents. One case already resolved is the resolvent density $u_{\q,\r}^{\{0}(x,y)$
with Parisian reflection at  Poisson observation times of intensity $\r $, obtained in
\cite[Thm.~4.1]{perez2018optimal}. It is not easy to prove that their result  converges when $\r   \to \I$ to the classic one in \cite[(22), Cor.~2]{Ivapot}.
\eeQ

\beQ It would be interesting to generalize the $W,Z$ formalism in a way which applies also to the case
of periodic observations of the smooth boundary. \eeQ
\beR
Some of the results above have been extended   to processes $X_\de^{[0[}(t) $ with classic reflection at $0$ and refraction at the maximum \cite[(3),Thm.~3.1]{AIpower}, and to processes $X_\de^{b[}(t) $ with $\de$-refraction at a fixed point $b$ \cite{KL,Kyp,KPP,Ren,PY}.

 Thus, \eqref{Parisbailouts} holds with $Z_{\q}(x,\th)$ replaced by $Z_{\q}^{\fr 1{1-\de}}(x,\th)$ \cite[Thm.~3.1]{AIpower}.  The proof uses the  probabilistic interpretation
$ \Ez_x [e^{-\q \tb - \theta \L_{\tb}}]=P[\tb < e_\q
\wedge K_{\theta}],$
where $K_{\theta}$ is the first time when the total bail-out exceeds
an independent exponential random variable $e_\theta$.
Finally,  \cite[(22)]{AIZ}  extend this to the case when $\tb$ is replaced by its Parisian version.

Similar results hold also
for  processes $X_\q^{b[}(t) $ with $\de$-refraction at a fixed point $b$ \cite{KL,KPP,Ren,PY}. The
scale functions are:
\begin{align}
\wrb_q(x)&=W_q(x)+\delta\int_b^x\mathbb{W}_q(x-y)W_q'(y)dy, \\
\zrb_{q}(x,\theta)&=Z_{q}(x,\th)+
\delta\int_b^x\mathbb{W}_q(x-y)Z_{q}'(y,\th)dy,
\end{align}
where $\mathbb{W}_q$ is the scale function of $X_t-\delta t$.

 For example, by \cite[Cor. 2]{KPP}, it holds that
\beq \la{nbr}
&&\E_x \Big[e^{- \r   T^{<0}}\Big]=P_x[T_{0,-}=\I]= (\k'(0_+) -\q)\fr{\F(\r  )}  {\r  - \q \F(\r  )}  \zrb_\q(x,\F(\r  )), \quad  0\leq \q \leq \k'(0_+). \eeq
\eeR

\ssec{Elements of proof for Theorem \ref{twostep}} \la{s:pr}

In the following, we provide some proofs for Theorem \ref{twostep}.
Before   that, let us record some useful preliminaries.

 \beP \la{p:st} For $z\leq0$, \ith

 \begin{enumerate}[A)] \im   the ``recovery before Parisian ruin'' probability is
 \bea &&P _z[T_{\{0\}}<e_\r ]=\E[ e^{-\r   T_{\{0\}}}]= e^{\Phi(\r  ) z}\\&& \E _z[e^{-\q T_{\{0\}}}; T_{\{0\}}<e_\r ]=\E[ e^{-(\r  +\q) T_{\{0\}}}]= e^{\Phi(\r  +\q) z}.\no \eea

\im  $$\E_z \left[e^{-\q e_\r  + \th X_{e_\r }};
T_{\{0\}}<e_\r \right]=e^{\Fqr z}\E_0 \left[e^{-\q e_\r  +\th X_{e_\r }}\right]=e^{\Fqr z}\frac{\r  }{\r  +\q -\k(\th
)}, \for \th \neq \F(\r  ).$$

\im
\begin{equation*} \label{eqexptime} \E_z\left[e^{-\q e_\r  +\th X_{e_\r }}; e_\r <
T_{\{0\}}\right]=\E_z \left[e^{-\q e_\r  + \th X_{e_\r }}
\right]-e^{\Phi(\r  ) z}\E \left[e^{\th
X_{e_\r }} \right]=\frac{\r  }{\r  + \q -\k(\th
)}
\bigl(e^{\th z}-e^{\Phi(\r  + \q) z}\bigr), \quad \th\geq0.
\end{equation*} %
\end{enumerate}
\eeP

\prf A) The second equation follows from the first, which is just the \fund identity \eqref{e:exp} (or  set $z \leq 0$ in \eqr{fhit}). B) follows  by the strong Markov property at $T_{\{0\}}$, and C) follows from B).

{\bf Proof of Theorem \ref{twostep}.2}
By the strong Markov property,  we may decompose $l(x,b):=\E_x [e^{-q\tb-\theta \L_{\tb}}] , \theta>\r  + \q,$ in three parts:
\begin{equation*}
\begin{split}
&l(x,b)
=\E_x[e^{-q\tb };\tb <\tz]
+\E_x\left[e^{-q\tz} \E_{X_{\tz}}[e^{-qT_{\{0\}}}; T_{\{0\}}<e_\r ];\tz<\tb  \right]l(0,b)  \\ &+\E_x\left[e^{-q\tz}\E_{X_{\tz}}[e^{-qe_\r +\theta X_{e_\r }}; e_\r <T_{\{0\}}];\tz<\tb   \right]l(0,b)=\frac{W_q(x)}{W_q(b)}+\\
&l(0,b)\left[\E_x[e^{-q\tz +  \Fqr  X_{\tz}}; \tz<\tb]+C\right]=\frac{W_q(x)}{W_q(b)}+l(0,b)\left[{Z_{\q}(x, \Fqr )}-
W_q(x)\frac{Z_{\q}(b, \Fqr )}{W_q(b)} +C\right]
\end{split}
\end{equation*}
where we have used Proposition \ref{p:st} A).

For the third part we  use Proposition \ref{p:st} C).
 We find
 \bea &&C=\E_x\left[e^{-q\tz}\E_{X_{\tz}}[e^{-qe_\r +\theta X_{e_\r }}; e_\r <T_{\{0\}}];\tz<\tb   \right]\\&& =\frac{\r  }{\r  + \q -\k(\th
)}\E_x\left[e^{-q\tz}\bigl(e^{\th {X_{\tz}}}-e^{\Phi(\r  + \q) {X_{\tz}}}\bigr);\tz<\tb   \right].\eea
Finally
\begin{equation*}
\begin{split}
&l(x,b)=\left\{\frac\r {\r  + \q -\k(\theta)}\left(Z_{q}(x,\th)- Z_{\q}(x, \Fqr )-W_q(x)\frac{Z_{q}(b,\th)-Z_{\q}(b, \Fqr )}{W_q(b)}\right)     +Z_{\q}(x, \Fqr ) \right.\\
&\left. -W_q(x)\frac{
Z_{\q}(b, \Fqr )}{W_q(b)}\right\}l(0,b) +\frac{W_q(x)}{W_q(b)}=\le\{Z_{\q,\r}(x,\th)-W_q(x)
\frac{Z_{\q,\r}(b,\th)}{W_q(b)}\right\}l(0,b) +\frac{W_q(x)}{W_q(b)}.\\
\end{split}
\end{equation*}

Now in the finite variation case we may substitute $x=0$, and, using $W_q(0) >0,$ conclude that $l(0,b)=\frac 1{Z_{\q,\r}(b,\th)}$, which yields the result.

In the infinite variation case, we may use a  perturbation approach.
For $b>x>0$,  we have
\begin{equation}\label{g}
\begin{split}
l(0,b)
&=\E [e^{-q\tau^+_x}; \tau^+_x<e_\r ]l(x,b)+\E[e^{-qe_\r +\theta X_{e_\r }}; e_\r <\tau^+_x, X_{e_\r }<0]l(0,b) \\
&+\int_0^x \E[e^{-qe_\r }; e_\r <\tau^+_x, X_{e_\r }\in dy]l(y,b)dy=e^{-  \Fqr  x}l(x,b) +I_2(x)l(0,b)+I_3(x),
\end{split}
\end{equation}
\begin{equation*}
\begin{split}
I_2(x)&=\r \int_{-\infty}^0\left(e^{-  \Fqr  x}W_{\r  + \q}(x-y)-W_{\r  + \q}(-y)\right)e^{\theta y}dy\\
&=\r \int_0^\infty e^{-  \Fqr  x-\theta y}W_{\r  + \q}(x+y) dy-\frac\r {\k(\theta)-q-\r  }\\
&=\r \int_x^\infty e^{-  \Fqr  x-\theta (z-x)}W_{\q +\r }(z) dz-\frac\r  {\k(\theta)-q-\r   }\\
&=\frac\r  {\k(\theta)-q-\r   }(e^{-  \Fqr  x+\theta x}-1)- \r  \int_0^x e^{-  \Fqr  x-\theta (z-x)}W_{\q +\r }(z) dz   \\
&=\frac\r  {\k(\theta)-q-\r   }(e^{-  \Fqr  x+\theta x}-1)+o(W_q(x)).
\end{split}
\end{equation*}

We can check that
\begin{equation*}
\begin{split}
&e^{-  \Fqr  x}(Z_{\q}(x, \Fqr )-Z_{q}(x,\th))\\
&=e^{-  \Fqr  x}\left[e^{  \Fqr  x}(1-\r   \int_0^x e^{-  \Fqr  y}W_r(y) dy)-e^{\theta x}(1-\r   \int_0^x e^{-\theta y}W_r(y) dy)  \right]   \\
&=1-e^{-  \Fqr  x+\theta x}+o(W_q(x)),
\end{split}
\end{equation*}
\[Z_{\q}(x, \Fqr )=e^{  \Fqr  x}\left(1-q\int_0^x e^{-  \Fqr  y}W_q(y)dy \right)  =e^{  \Fqr  x}+o(W_q(x)), \text{   and}\]
\[ I_3(x)\leq \int_0^x E[e^{-qe_\r }; e_\r <\tau^+_x, X_{e_\r }\in dy]dy=\r  \int_0^x e^{-  \Fqr  x}W_{\q +\r }(x-y)dy=o(W_q(x)). \]

Solving now (\ref{g}) for $l(0,b)$ and letting $x\goto 0+$, we find again
\begin{equation*}
\begin{split}
l(0,b)&=\lim_{x\goto 0+}\frac{e^{-  \Fqr  x} \frac{W_q(x)}{W_q(b)}}
{e^{-  \Fqr  x} W_q(x)\frac
{Z_{\q}(b, \Fqr )}{W_q(b)}+\r   e^{-  \Fqr  x} W_q(x)\frac{Z_{\q}(b, \Fqr )-Z_{q}(b,\th)}{(\k(\theta)-q-\r   )W_q(b)}+o(W_q(x))}  \\
&=\frac{\k(\theta)-q-\r   }{(\k(\theta)-q)
Z_{\q}(b, \Fqr )-\r   Z_{q}(b,\th) }=\frac 1{Z_{\q,\r}(b,\th)}.
\end{split}
\end{equation*}

\ssec{Spectrally negative Omega Processes\la{s:ome}}

Recently, it was discovered that the classic exponential Parisian formulas may be  further extended to Omega models,  \cite{albrecher2011optimal,gerber2012omega,LP,LZ17}, in which a state-dependent
  rate of  killing (or observation) rate $\w(x)$ is used, where $\omega:\mathbb{R}\to\mathbb{R}_{+}$ is  an arbitrary locally bounded nonnegative measurable function. Exponential Parisian models are just the particular case when $\w(x)$ is a step  function with two values. 

  Analogs of Propositions \ref{p:twos}, \ref{p:res} and of Theorems \ref{l:s}, \ref{l:refl} are provided in  \cite[Thm.~2.1-2.4]{LP}, who  showed that the first passage theory of Omega models rests on two functions $\{\mW_\w(x), x\in \mathbb{R}\}$ and $\{\mZ_\w(x), x\in \mathbb{R}\}$
 called $\omega$-scale functions,  which are defined uniquely as the solutions of the  renewal equations:
\begin{align}
\mW_\w(x)=&\ W(x)+ \int_0^{x} W(x-y) \omega(y) \mW_\w(y)\,dy, \label{eqn:Ww}\\
\mZ_\w(x)=&\ 1+ \int_0^{x} W(x-y) \omega(y) \mZ_\w(y)\,dy,\label{eqn:Zw}
\end{align}
 where $W(x)$ is the classical zero scale function.

Furthermore, \eqref{eqn:Ww}, \eqref{eqn:Zw}   may be  generalized to nonhomogeneous models \cite[Lem.~3]{LZ17}:
\begin{align}
\mW_{ \T \w}(x,a)=&\mW_\w(x,a)+ \int_0^{x} \mW_{\w}(x,y) \le(\T \w(y)-\w(y)\ri) \mW_{ \T \w}(y,a)\,dy, \label{eqn:Wwn}\\
\mZ_{\T \w}(x,a)=&\ Z_\w(x,a)+ \int_0^{x} W_{ \w}(x,y) \le(\T \w(y)-\w(y) \ri) \mZ_{ \T \w}(y,a)\,dy. \label{eqn:Zwn}
\end{align}

Note that in the case of constant $\omega(x)=\q$, these  reduce
 \begin{equation}\label{falka}
 W_{\q}-W= \q W_{\q}* W \quad\mbox{and}\quad  Z_{\q}-Z= \q  W_{\q}*Z,
\end{equation}
which can be easily checked by taking the Laplace transforms of their both sides and by using the expansion \eqref{qexp}.

 \ssec{Occupation times}
Here is an elegant  result \cite[Thm.~3.1]{li2015two}
  on the joint law of the occupation times above and below $0$ of a spectrally negative L\'evy process.

 \beP  Introduce the auxiliary function \cite[(1)]{li2015two} (a slight modification of which had essentially  appeared  already in \cite[6]{LRZ}),
defined  for all $x \in \R$ and $\r  ,\q \geq 0$ by:
\begin{equation} \label{W_a_def}
\mW^a_{\r  ,\q}(x) :=  \bc W_{\r } (x),& 0\leq x \leq a\\
W_{\r } (x)+(\q-\l)\int_{a}^{x}  W_{\q}(x-y)W_{\r } (y)\diff  y =W_{\q}(x)+(\r -\q) \int_0^{a}W_{\q}(x-y)W_{\r } (y) \diff  y, & 0\leq a \leq x \\ W_{\q}(x), & a\leq 0\ec
\end{equation} where the second equalities hold by the convolution identity
$W_\r  \, * W_{\q}(x)=\fr{W_\r  (x)- W_{\q}(x)}{\r -\q}$ \cite[(5)]{LRZ}.\fn[4]{Note that
these functions satisfy \cite[(2.18)]{APY}
	$
	\lim_{a\downarrow -\infty} \frac{\mW_{\r  ,\q}^a(x)}{W_{\r } ({a})}
	= Z_{\q}(x,\Phi(\r ) ) 	\quad
\textrm{and} \quad \lim_{x\uparrow \infty}\frac{\mW_{\r  ,\q}^a(x)}{W_{\q}(x)}
	=Z_{\r } ({a},\Phi(\q))
	$.}
Let $L^-_t=\int_0^t \1_{(-\I,0)}(X_s) ds, L^+_t=\int_0^t \1_{(0,\I)}(X_s) ds$ denote the occupation times  below and above $0$. Then, $\for \r  _-, \r  _+ >0$ and $\for x,y \in \R$ it holds that
\bea && \int_0^\I e^{-\q t} \E_x \left[e^{-\r  _- L^-_t- \r  _+
L^+_t},X_t  \in d y\right] d t\\&&=\left(\fr{\F(\q + \r  _+)-\F(\q + \r  _-)}{\r  _+-\r  _-} Z_{\q + \r  _+}(x, \F(\q + \r  _-)) Z_{\q + \r  _-}(-y, \F(\q + \r  _+))- \mW_{\q+\r _-,\q+\r _+}^{-y}(x-y)\right) dy.\eea
\eeP
\beR Starting from $x=0,$ the result loses its symmetry, and simplifies to  \cite[Thm.~3.1, Rem.~3.2]{li2015two}
\begin{equation*}
\begin{aligned}
(d y)^{-1} \int_0^\I e^{-\q t} \E_0 \left[e^{-\r  _- L^-_t- \r  _+ L^+_t},X_t  \in d y\right] d t&=\fr{\F(\q + \r  _+)-\F(\q + \r  _-)}{\r  _+-\r  _-} Z_{\q + \r  _-}(-y, \F(\q + \r  _+)-W_{\q+\r _-}(-y)\\&=\fr{\F(\q + \r  _+)-\F(\q + \r  _-)}{\r  _+-\r  _-} \E_{-y} \left[e^{-(\q + \r  _-) \tz + \F(\q + \r  _+) X_{\tz}}\right].
\end{aligned}
\end{equation*}

Integrating the final position yields \cite[Cor.~3.1]{li2015two}
\bea && \int_0^\I e^{-\q t} \E_0 \left[e^{-\r  _- L^-_t- \r  _+ L^+_t}\right] d t=\fr{\F(\q + \r  _-)}{(\q + \r  _-)\F(\q + \r  _+) }.\eea

This implies  \cite[Rem.~4.1]{LaRZ}, \cite[Cor.~3.2]{starreveld2016occupation} \beq \la{SM} && \int_0^\I e^{-\q t} \E_0 \left[e^{- \r   L^+_t}\right] d t=\fr{\F(\q)}{\q \F(\q + \r  ) }.\eeq
\eeR

\beR \textbf{Asymptotics of occupation times for a reflected process}.
A  general result for the time $L^{[0,b]}_t=\int_0^t \1_{[0,b]}(X_s) d s$ spent in $[0,b]$   by a process with positive drift (and thus with $\F(0)=0$) reflected at $b$ is provided in
\cite[Thm.~3.4]{starreveld2016occupation}:
\beq && \int_0^\I e^{-\q t} \E_0 \left[e^{- \r   L^{[0,b]}_t}\right] d t=\fr{\F(\q)}{\q}
 \fr{Z_{\l}(b,\Fq)}{\r   W_\r  (b)+ \F(q) Z_{\l}(b,\Fq)}, \la{roct}\eeq
which recovers the previous result \eqref{SM} by using $\lim_{b \to \I}
\fr{Z_{\l}(b,\Fq)}{ W_\r  (b)}=\fr{\r  }{\F(\q+\r )-\F(\q)}$.

The  large deviations rate  for $L^{[0,b]}_t$ has been obtained in  \cite[Thm.~3.3]{starreveld2016occupation}, as  a direct consequence of the G\"artner-Ellis theorem, which states that this   is  the Legendre transform of
\begin{equation} \r (r):= \lim_{t \to \I} \fr1{t}
\log\Big[\E[e^{-r  L^{[0,b]}_t}] \Big]=\lim_{\q \to 0} \fr{\Fq}{\q} \fr{Z_{\l}(b,\Fq)}{\r   W_\r  (b)+ \Fq Z_{\l}(b,\Fq)}=\fr{1}{p} \fr{Z_\r  (b)}{\r   W_\r  (b)}.\end{equation}
\eeR

 \sec{Optimization of dividends  \la{s:div}}

 Risk theory  initially revolved around minimizing the probability of
ruin.  However,  insurance companies are realistically more interested in maximizing company  value than minimizing risk and an alternative approach is therefore
to study optimal dividend policies, in the sense of maximizing the expected value of the
sum of discounted future dividend payments until the time of ruin, as suggested by De Finetti in the 1950 \cite{deF}-- se also  Miller and Modigliani   \cite{miller1961dividend}.

A second interesting objective to maximize introduced by Shreve, Lehoczky and Gaver (1984) \cite{SLG}, is the expected discounted cumulative dividends for the reflected process  obtained by redressing the reserves
by   capital injections, at a proportional cost, each time this becomes necessary.

These two objectives and certain generalizations are easily expressed for spectrally negative \lev processes  in terms of  the  scale functions $W,Z$ (at least when restricting to barrier policies).

\ssec{The de Finetti objective with Dickson-Waters modification for spectrally negative processes\la{s:deF}}
This objective proposed by de Finetti (1957) \cite{deF} is to maximize expected discounted dividends until the  ruin time. It makes sense  to include   a    penalization for  the final deficit \cite{dickson2004some}, arriving at:
   \beq\label{deFmod} &&V_{w}(x)=\sup_{\pi} V_{w}^{\pi}(x),\\&& V_{w}^{\pi}(x)=\E_x \le[ \int_0^{\tz} e^{- \q t} d \U^\pi_t +   e^{- \q \tz} w(X_{\tz})\ri]\no  := V^{\pi}(x)+   \Rui_{\q,w}^{\pi}(x). \eeq
   Here  $\U^\pi_t$ is an ``admissible'' dividend paying policy, and $w(x)$ is a bail-out penalty function\fn[4]{The value function must satisfy in a viscosity sense
     the HJB equation \cite[(1.21)]{azcue2014stochastic}:
   $ \mG ( V)(x):=\max[\mG_\q V(x),1-V'(x),V(x)-w(x)]=0,$
   where $\mG_\q V(x)$ denotes the discounted infinitesimal generator of the uncontrolled
surplus process, associated to the policy of continuing without paying dividends.
The second operator $1-V'(x)$ is associated to the possibility of modifying the surplus by a lump payment, and the third to bankruptcy.}.

The most important class of policies is that of constant barrier policies $\pi_b$, which modify  the surplus  only when  $X_t>b$,
  by a lump payment bringing the surplus at $b$, and then  keep it there by Skorokhod reflection, until the next negative jump\fn[3]{In the absence of a Brownian component, this  amounts to paying all the income while at $b$}, until the next claim.

Under   a reflecting barrier strategy $\pi_b$, the dividend part of the de Finetti objective
   has a simple expression \eqr{div} in terms of the $W$ scale function :
   \begin{equation*}  V^{b]}(x)=\Eazb_x\le[\int_{[0,T_0^{b]}]}e^{-\q t}d   \U_t \ri]= \frac{W_\q(  x)}{W_{\q}^ {\prime}(b)}, 
\end{equation*}
where  $\Eazb$ denotes the law of the process  reflected from above at $b$,
and absorbed at $0$ and below. This formula reflects the representation $$V^{b]}(x)=\E[ e^{-q \tb} ; \tb < \tz] \Eazb_b\le[\int_{[0,T_0^{b]}]}e^{-\q t}d   \U_t\ri]=\E[ e^{-q \tb} ; \tb < \tz]\;  \Eazb_b\le[   \U_{T_0^{b]} \wedge \kil_q}\ri],$$ and the fact that the local time $\U_t$ at $b$ with reflection at $b$ is an exponential random variable.

The ``barrier function'' \beq H_D(b):=\fr {1}{ W'_\q(b)}, \; b \geq 0, \la{GDeF}\eeq
  plays
a central  role  in the solution of the problem, and   the optimal dividend policy is  often a barrier strategy at its   maximum. In particular, when the barrier function is differentiable and has a unique local maximum $b^*>0 \Lra W''_\q(b^*)=0$, this $b^*$ yields the optimal dividend policy. Furthermore,   the value function
  \beq \la{dbpr} V(x)= \sup_{b \geq 0} V^{b]}(x)=V^{b^*]}(x)\eeq
  is  then the largest concave minorant  of  $W_\q(x)$.  In the presence of several inflection points, \how \ the optimal policy is multiband \cite{AM05,Schmidli,Loef,APP15}.

  The first numerical examples of multiband policies  were produced in \cite{AM05,Loef}, by Cram\'{e}r-Lundberg model \eqref{CL} with  Erlang claims $E_{2,1}$. However, it was shown in \cite{Loef} that multibands cannot occur when $W'_\q(x)$ is increasing after its  last global minimum  $b^*$ (i.e. when no local minima are allowed after the global minimum).\fn[5]{One instance when that happens is when the \Lm \ is completely monotone. Then, \eqr{Wu} may be written as $  W_{\q}(x)= \Fqp e^{\Fq x} -{\Fqp} \int_0^\I e^{-x t} \mu_\q(dt), \; x \geq 0,$
for some finite measure $\mu_\q$. This implies $W_\q'''(x) \geq 0, x \geq 0$, and implies finally that  $W_\q'(x) $ is convex, with a unique minimum.}

  \cite{Loef} further made the interesting  observation that in the Brownian perturbed Cram\'{e}r-Lundberg model \eqref{pCL} with Erlang claims $E_{2,1}$  (which are non-monotone), multiband policies may occur  for $\s$ smaller than a threshold value, but  barrier polices (with non-concave value function!) will occur when $\s$ is big  enough.

  Figure \ref{f:pd} displays the first derivative $W_\q'(x),$ for $\s^2/2 \in \{ {\color{blue} 1/2}, {\color{yellow}1}, {\color{green} 3/2}, {\color{red} 2} \} $. The last two values yield barrier polices with non-concave value function, due to the presence of an inflection point in the interior of the interval $[0,b^*]$.
   \figu{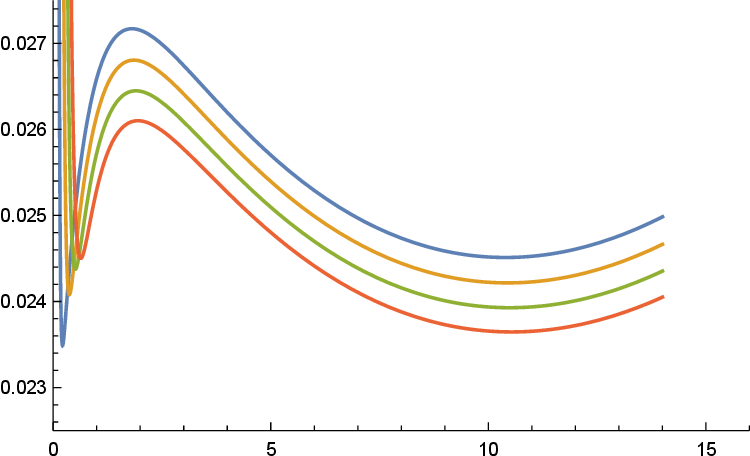}{Graphs of the Loeffen example for $\k(s)= \frac{\s^2 s^2}{2}+c \; s+\l
   \left(\frac{1}{(s+1)^2}-1\right)-\q, c=\frac{107 }{5}, \l=10, q=\frac{1}{10} $, $\s^2/2 \in \{ \color{blue} 1/2, \color{yellow}1, \color{green} 3/2, \color{red} 2 \} $.}{0.7}

Even when barrier strategies do not  achieve the  optimum,  and  multi-band policies must be used instead, constructing the solution must start by
  determining the global maximum of the barrier function   \cite{AM05,Schmidli,APP15}.
  We will only consider barrier strategies in this review.

The penalty part of the objective \eqref{deFmod} for a barrier strategy $\pi_b$ can be expressed  as $\Rui^b_{\q,w}(x)=
{G}_w(x)-W_\q(x)\fr {{G}_w'(b)}{ W'_\q(b)}$ \eqref{rserui}, where ${G}_w(x)$ is  the smooth \GS function associated to the penalty $w$ (see Section \ref{s:nonh});
finally, the modified de Finetti value function is:
\beq
V_{w}^{b]}(x)=\bc {G}_w(x) +   W_\q(x)\fr {1-{G}_w'(b)}{ W'_\q(b)}& x \leq b\\x-b + V_{w}^{b]}(b) &x \geq b  \ec  \la{VDeF}
\eeq

The corresponding barrier function is \beq H_w(b):=\fr {1-{G}_w'(b)}{ W'_\q(b)}, \; b \geq 0. \la{GDeFD}\eeq

 The most important cases of bail-out costs ${w}(x)$ are
\BEN \im exponential $w(x)=e^{\th x}$, when ${G}_w(x)=Z_{\q}(x,\th)$ (Proposition \ref{p:GSexp}), and \im  linear  $w(x)= k x -K$. For $x < 0$, the constants $k >0$ and  $K\in \R$ may be  viewed as  proportional and  fixed bail-out costs, respectively.\fn[6]{The cases  $k \in (0, 1]$ and $k >1$   correspond  to  management being held responsible
for only part of the deficit at ruin, and to having to pay extra costs  at liquidation, respectively. When $K < 0$, early liquidation is rewarded; when $ K> 0$,  late ruin  is rewarded.}
In this case as well, ${G}_w(x)$ may be obtained by using
$Z_{\q}(x,\th)$ as generating function in $\th$, i.e.~the coefficients of $K,k$ in ${G}_w(x)$ are found by differentiating with respect to $\th$ the $Z_{\q}(x,\th)$ scale function  
$0$ and $1$ times respectively, and taking $\th=0$. This yields
\beq {G}_w(x) =k Z_{\q}^{(1)}(x) - K Z_\q(x ),\la{e:linct}\eeq
where $Z_{\q}^{(1)}(x)$ is given by \eqref{Z1}. 
In the simple, but important particular case $w(x)=-K$, the modified de Finetti value function  and barrier function are
\resp
  \beq \la{HK} && V_{K}^{b]}(x)=-K Z_\q(x) +   W_\q(x)\fr {1+K Z_\q'(b)}{ W'_\q(b)},\\&&H_K(b):=\fr {1}{ W'_\q(b)}+ K\fr { Z_\q'(b)}{ W'_\q(b)}=\fr {1+K \q W_\q(b)}{ W'_\q(b)}. \no
\eeq

\EEN

\beR
Optimality largely rests on the sign of the numerator $$H_w'(b)=\fr{-W_\q''(b)+(G_w' W_\q''-W_\q' G_w'')(b)}{(W_\q')^2(b)}.$$

For \eqref{HK} for example,
\beq H_K'(b)=\fr{ K \q \D_q^{(W)}(b)-W_\q''(b)}{(W_\q')^2(b)},\eeq
where
\be \la{DW} \D_q^{(W)}(b):=\le((W_\q')^2 -W_\q W_\q''\ri)(b)=(W_\q')^2(b) \fr{d}{d b}  \left( \fr{W_\q}{W_\q'} \right) (b).\ee

Since the excursion rate $\nu(b)=\fr{W_\q'}{W_\q}(b) $ is by definition decreasing (see Remark \ref{r:exrate}),
it follows that  $\D_q^{(W)}(b) \geq 0$.\fn[3]{incidentally, when $\s >0,$  this  is also implied by  the {\bf   creeping drawdown  law} \cite{MP}, \cite[(2.5)]{LLL}:
\beq \E_x \Big[ e^{- \q \t_a} ; Y_{\t_a} =a \Big] = \fr{\s^2} 2 \fr{\D_q^{(W)}(a)}{W_\q'(a)}, \; \for x. \eeq}

Let $b_0$ denote the last maximum of   the unconstrained  $H_D(b)$,   and, $\for b \geq b_0$,  let
\beq K(b)= \fr{W_\q''(b)}{q \D_q^{(W)}(b)} \geq 0,  \la{Kb} \eeq
denote the unique $K \geq 0$ satisfying $H'_K(b)=0$.

Then, assuming {\bf  \cmy of the \Lm}, \cite[Prop.~4.5, Thm.~4.4]{hernandez2016time} show that for every $b \geq b_0$ $K(b)$  is
strictly increasing. \Thr barrier policies are optimal and $b$ yields the optimal barrier  for the cost $ K(b)$ (in their paper,    the parameter $K$  intervenes as a Lagrange multiplier associated to a time constraint).
\eeR

 \ssec{Optimal de Finetti dividends barrier until Parisian ruin}
Differentiating \eqr{divPar} and using twice \eqr{Zder}, we find that the optimal de Finetti dividends barrier $b$ until Parisian ruin
must satisfy
\be {\th}(\fr{\th}{\l} Z_{\q}(b,\th) -  W_q(b)) ={ W_\q'(b)}, \quad \th = \Fqr  \la{bFP}\ee
(note that
the same  equation was obtained in \cite{noba2018optimal}
in the context of a different, but equivalent problem
involving running costs).

When $\r \to \I,$ the LHS of \eqr{bFP}
converges to $W_\q''(b)+W_\q'(b)$ by \eqr{ZWF}. Thus, $\lim_{\r \to \I} b_\r^\ast  = b^\ast$, recovering the classic optimality
equation.

An important case is that when the optimal  dividends barrier is $0$; this
may be viewed as a measure of the process involved corresponding  to an "efficient company" (ready to pay dividends) -- see  \cite{AM16}. The "efficiency"
 condition here is
 $$\Fqr(\frac{\Fqr}{\r} -W_\q(0)) \geq W_\q'(0)$$
 see also \cite{renaud2019finetti}.

\ssec{The Shreve-Lehoczky-Gaver infinite horizon objective, with linear penalties} \la{s:SLG}

We turn now to an objective  which was first  considered in a diffusion setting by Shreve, Lehoczky,  and Gaver (SLG) \cite{SLG} -- see also  \cite{boguslavskaya2003optimization,lokka2008optimal} -- to be called  SLG objective.

Suppose  a subsidiary must be bailed out each time its surplus is negative,
and  assume the penalty costs are linear  $w(x)=k x $. The optimization objective   of interest  combines   discounted dividends $\U_t$, and cumulative bailouts $\L_t$ 
  \beq\label{e:Sobj} &&V_{S,k}(x)=\sup_{\pi} V_{S,k}^{\pi}(x), \no\\&&
  V_{S,k}^{\pi}(x)=\E^{\pi}_x \le[ \int_0^\infty e^{- \q t} d \U^\pi_t - k \int_0^\infty e^{- \q t} d \L_t^\pi\ri]
  \eeq
  where $\pi$ is a dividend/bailout 
  policy, and $k \geq 1$. 

  Importantly, for \lev processes the optimal dividend/bailout policy $\pi$  is always of {\bf constant barrier} type \cite{APP}
    , and the objective for fixed $b$ has the simple expressions
 provided in   \cite[(4.3),(4.4)]{APP}
(and included above as \eqref{divSLG}, Theorem \ref{l:div} and \eqref{VLSLG}, Theorem \ref{l:VL}), resulting in\fn[4]{As already noted in Remark \ref{r:W2z}, this has the same form as the de Finetti objective 
\eqref{e:linct} with
$Z$ replacing $W$.}:
\begin{eqnarray}\label{value-SLG}
&&V_{S,k}^{[0,b]}(x)=V^{[0,b]}(x)- k B^{[0,b]}(x)=\fr{Z_{\q}(x)}  {Z'_{\q}(b)} +k \le(Z_{\q}^{(1)}(x) -\fr{Z_{\q}(x)}  {Z'_{\q}(b)} (Z_{\q}^{(1)})'(x) \ri) \no \\&& =k \left( \overline{Z}_\q(x) + \frac{\k'(0_+)}{q} \right) +Z_\q(x) H_k^{SLG}(b),
\end{eqnarray}
with barrier function
\beq H_{k}^{SLG}(b)=\fr{1 -k Z_{\q}(b)}
{\q W_{\q}(b)}
 \la{HS}\eeq
 -- see also \cite[Prop. 3.1]{wang2018dividend} for a generalization
 involving fixed dividend costs $K$. This impulse control problem involves replacing   the reflection barrier by a $b_1,b_2$ band. It turns out that the \valf \ is of the same form, but the barrier function changes, to
 $$H_{k,K}^{SLG}(b)=\fr{b_1 - b_2 -K - k \le(\ol Z_\q(b_2)-\ol Z_{\q}(b_1)\ri)}
{  Z_\q(b_2)- Z_{\q}(b_1)}.$$
Note that the derivation becomes simpler  than in the reflection case.

The next proposition merges new results from \cite[Prop. 1]{AGR} with previously known results from \cite[Lem.~2]{APP}. The main object is the function $k_f \colon [0,\I) \to [k_0,\I)$ defined by
\begin{align} & k_f(b):=\fr{W_{\q}'(b)}{Z_{\q}(b) W_{\q}'(b)- \q  W_{\q}^2(b)} , \quad b > 0 ,\\& k_0:=k_f(0_+)=\fr{W_{\q}'(0_+)}{W_{\q}'(0_+)- \q  W_{\q}^2(0_+))} =
\bc
1 , & \text{if $X$ is of unbounded variation,}\\
1 + \fr{\q}{\lm(0,\I)} , & \text{if $X$ is of bounded variation.}
\ec\la{kf} \end{align}

This function \red{is increasing, by the well known identity  \cite[Thm 1]{AKP}\fn[4]{Some papers refer to this as the log-convexity of  $Z_\q(x)$.}
$$
\E_x \left[ \mathrm{e}^{-q \tau^b} \right] = Z_q(x) - q \frac{W_q(b)}{W_q^\prime (b)} W_q (x) \; \Lra k_f(b)=\fr 1{\E_b \left[ \mathrm{e}^{-q \tau^b} \right]},
$$
and since the map $b \mapsto \E_b \left[ \mathrm{e}^{-q \tau^b} \right]$ is decreasing.}

The monotonicity  allows us to re-parametrize the problem in terms of the optimal barrier $b_k$ associated to a fixed cost $k$.

\beP  \la{p:APP} Assume $X$ is a SNLP and $K=0$. We have the following results:
\BEN

\im For fixed $x,b$, the function $k \mapsto V_k^{0,b}(x)$ defined in~\eqref{value-SLG} is non-increasing.

\im For $k=k_f(b)$, the value function defined in~\eqref{value-SLG} can be written as follows:
\begin{eqnarray} \la{Vkf}
&&V_{k_f(b)}^{0,b}(x) = k_f(b) \left[ \ol Z_{\q}(x) +  \frac{p}{q} -Z_\q(x) V^{b]}(b)  \right]= k_f(b) \left[ Z_{\q}^{(1)}(x) + Z_\q(x) \left( \frac{p}{q} - \fr{W_\q(b)}{W_\q'(b)} \right) \right] ,
\end{eqnarray}
where
$
Z_\q^{(1)} (x)$ is defined in \eqr{Z1}, and $V^{b]}(b)$ is the \deF objective when starting at the barrier.

\im For fixed $k$, the barrier function $H_k^{SLG}$ defined in~\eqref{HS} is an increasing-decreasing function with a unique maximum $b_k \geq 0$.   Moreover, if $b_k>0$, then $k_f(b_k)=k$.

\EEN
\eeP

\begin{proof}

\BEN

\im This is obvious since the \SLG  value function \eqref{value-SLG} is decreasing in $k$, and
the value function $V_k^{0,b}(x)$ can be seen as the maximum of $\E_x \left[ \int_0^{\infty} \mathrm{e}^{-q t} \left( \mathrm d D_t - k \mathrm d C_t \right) \right]$ over control couples $(C,D)$ keeping the surplus in $[0,b]$.  Since the cost functional is non-increasing in $k$, our assertion follows.
\im Recalling \eqr{value-SLG}, we need to show that
\be -H_{k_f(b)}^{SLG}(b)=k_f(b) V^b(b).\ee
Indeed, it is easy to check that the equality
$$\frac{k Z_q(b)-1}{q W_q(b)}=k \fr{W_q(b)}{W_q^\prime(b)}$$
holds for $k=k_f(b)$.

\im  For the sake of completeness, let us reproduce this proof from \cite[Lem. 2]{APP}. The  derivative of the barrier function \eqref{HS} \sats
\be
q \la{f} \fr{H' W_q^2}{W_q'}(b)=f(b):=k \; \fr{\D_{\q}^{(ZW)}(b)}{ W_{\q}'(b)}-1= k \; \E_b [e^{-\q \t^{b}} ] -1 =\fr k { k_f(b)}  -1 ,
\ee
where $\D_{\q,0}^{(ZW)}=Z^{(\q)}(b)  W_{\q}'(b)-\left(Z^{(\q)}\right)'(b) W_{\q}(b)$ (see  \eqref{deltaZW}).
The sign of the derivative of the barrier function \eqref{HS} coincides \thr with that of   $ f(b)= k \; \Eb_b [e^{-\q T_0^{b]}} ] -1 $. Clearly the latter function $f$ is  decreasing in $b$ from $\lim_{b \to 0} f(b)=\fr{k}{k_0}-1$ to $-1$.


\EEN

\end{proof}

\beR \la{b**}
We may conclude therefore that  if $$k \leq k_0 \Eq f(0) \leq 0 \Eq  \fr{W_{\q}'(0)}{W_{\q}'(0)- \q  W_{\q}^2(0))} \geq k,$$ then $b_k=0$ is the optimal barrier, and otherwise there is a unique global and local maximum
\satg $$\fr{W_{\q}'(b_k)}{Z_{\q}(b_k) W_{\q}'(b_k)- \q  W_{\q}^2(b_k))}=k={\T \de}_\q  ^{-1} (b_k), b_k>0.$$
\eeR

\beR
 The last identity in Proposition \ref{p:APP} turns out useful  in establishing the so called \LZ for  Brownian motion with drift -- see\cite{lokka2008optimal}, \cite{lindensjo2019optimal} --  and for  the \CL model with exponential jumps \cite{AGR}. These results state that, depending on the size of transaction costs, one of the following strategies is optimal:
\begin{enumerate}
\item if the cost $k$ of capital injections is below a  critical  point   $k_c$,   then it is optimal to pay dividends and to inject capital, according to a double-barrier strategy, meaning that ruin never occurs;
\item if the cost of capital injections is above the  critical  point   $k_c$,  it is optimal to use a  single-barrier strategy and declare bankruptcy at the first passage below $0$.

\end{enumerate}

The crucial point in these two cases is that  a further identity  holds which allows expressing  the RHS of \eqr{Vkf} in terms of the $W$ scale function, and implies \be \la{kid}V_{k_f(b^\ast)}^{0,b^\ast}(x)=V^{dF}(x),\ee
where $b^{\ast}$ denotes the optimal barrier level in de Finetti's problem.

\Mp in the \BM \ case,  note the easily checked identities
$$ Z_q^{(1)}(x) + Z_q(x) \left( \frac{p}{q} - V^b(b) \right)=Z_q^{(1)}(x)=\fr{\s^2}2 W_\q(b) \Lra V^{0,b^\ast}_{k_f({b^\ast})}(x)=V^{dF}(x),$$
and use then the monotonicity of $V^{SLG}_{k}(x)$ in $k$.

Similar computations  establish the \LZ in the \CL case with \expoc ~\cite{AGR}.
\eeR

\ssec{The dividends and penalty objective, with exponential utility}
  Given $\delta, \theta, \vartheta>0$, one may consider the barrier strategy obtained by minimizing the objective \eqref{DP}. Such an objective is based on exponential utility that rewards late ruin and cumulative dividends while penalizing deficit at ruin. Recall that the barrier function of \eqref{DP} is
  \bea  H_{DP}(b)=\fr{  Z_{\q}'(b,\th)+  \vartheta  Z_{\q}(b,\th)  }{ W_{\q}'(b)+ \vartheta W_{\q}(b)}.\eea

 For  $\th=\vt=0$, this reduces to $\q \fr{  W_{\q}(b)  }{ W_{\q}'(b)}$, which is  clearly an  increasing function.
  For $\th=0$, \eqref{DP} reduces to a dividends and time objective, with \Bf
 \beq H_{DT}(b)=\fr{  Z_{\q}'(b)+ \vt  Z_{\q}(b)  }{ W_{\q}'(b)+ \vt  W_{\q}(b)}. \la{th0}\eeq
 This  bounded function, with values in between $H_{DT}(0)=\fr{\q  W_{\q}(0) + \vt}{ W_{\q}'(0)+ \vt  W_{\q}(0)}$,
 and $H_{DT}(\I)=\fr{\q  + \vt \fr{\q}{\Fq}}{\Fq+ \vt }$, is the \Bf of the objective
 \begin{equation} \la{e:DT}
 DT^b(x,\vt):=\Eb_x\left[e^{ -\q T_0^{b]}  -\vt  \U_{T_0^{b]}}} \right]=\Eb_x\left[e^{  -\vt  \U_{T_0^{b]}}} ; T_0^{b]} < e_\q \right]=
 Z_{\q}(x) -  W_{\q}(x)   \fr{  Z_{\q}'(b)+ \vartheta  Z_{\q}(b)  }{ W_{\q}'(b)+ \vartheta  W_{\q}(b)}.
\end{equation}
\beR
Note that this objective encourages taking dividends soon; in fact, everything is lost at $e_\q$, which must be interpreted as a catastrophic event. An alternative would be to minimize $\Eb_x\left[e^{  -\vt  \U_{T_0^{b]} \wedge e_\q}} \ri]$,
which would also encourage taking dividends soon, but with less urgency. The optimal barrier for this last objective should increase with respect to that of \eqref{e:DT}.
\eeR

 \beR The  sign of the derivative of the barrier function \eqref{th0} of the exponentiated dividends and time objective  \eqref{e:DT} is determined by
\bea &&\le(Z''_{\q}(x) + \vt Z'_{\q}(x) \ri)
\le( W'_\q(x)+ \vt  W_{\q}(x) \ri) -\le( W''_\q(x)+ \vt W'_\q(x) \ri) \le( Z'_{\q}(x)+ \vt Z_{\q}(x)\ri).\eea\fn[4]{Even after simplification
$$ \q \Big(\vt^2 ( W_{\q}(x)^2-W'_\q(x) \ovl  W_{\q}(x))+ \vt (W'_\q(x) W_{\q}(x)-W''_\q(x) \ovl  W_{\q}(x))+ W'_\q(x)^2-W''_\q(x)  W_{\q}(x)\Big) -\vt \le( W''_\q(x)+\vt W'_\q(x)\ri), $$
this seems hard to analyze.}

Some numerical results involving the exponential utility barrier functions \eqref{HDP}, \eqref{th0} and their critical points are presented in Section \ref{s:num}.
We have never found multi-modal instances, suggesting that the optimal policy is simpler to implement than that for the de Finetti objective. 

   \eeR

\beR
 For comparison with \eqref{e:DT},    consider also the linearized value function (see Theorem \ref{p:time} C) and Theorem \ref{l:div} A))
   \begin{align*}
   \Eb_x \left[\q T_0^{b]} + \vt \U_{T_0^{b]}} \right] &= \q\le(W(x) \fr{W(b)}{W'(b)}- \int_0^x W(y)dy \ri) + \vt \fr {W(x)}{ W'(b)}\\
   &=- \q \int_0^x W(y)dy +W(x)  \fr{ \q W(b)+ \vt}{W'(b)}   
   \fr{\q W(b)+ \vt}{W'(b)},
   \end{align*}
    which needs to be maximized.

    The optimization \eqref{e:DT} may then be viewed
as a risk sensitive optimization with   {\bf exponential utility} $e^{- x}$, applied to the random variable $\q T_0^{b]} + \vt \U_{T_0^{b]}}$.
\eeR

 \ssec{Optimization of   dividends for   spectrally positive processes \la{s:divpos}}
The dividends of a spectrally positive process $X_t$ are the bailouts of its dual $-X_t$. Furthermore, for a fixed upper barrier $b$, the argument $x$ of the scale functions must be replaced by $b-x$. The end result for the de Finetti problem is
\cite[Lem. 2.1]{bayraktar2013optimal}
\beq \la{duDeF} V(x)= Z_\q(b-x) \fr{ G_{\q}^B(b-a)}{ Z_\q(b-a)}- G_{\q}^B(b-x), \quad G_{\q}^B(x)=\ovl Z_\q(x) + \fr{\k'(0_+)}{\q}, \; \q >0, x \leq b.\eeq
 Barrier policies $b^*$ are always optimal,  and smooth fit yields that $\ovl Z_\q(x)=\fr{p_+}{\q}$ \cite[Thm. 2.1]{bayraktar2013optimal}.

 Since stopping happens now without overshoot, the only relevant penalty of ruin is $w(x)=-K$,
 and \eqr{duDeF} still holds, with $G_{\q}^B(x)$ replaced by $G_{\q}^B(x)-K$ \cite[Thm. 3.1]{yin2013optimal}.

For Parisian observation of \deF dividends and a final ruin penalty $K$, the  value function
is given by \eqr{eq:localt0}, applied to $b-x$, and the optimal barrier must satisfy the \eq \cite[(3.40),Lem. 3.6]{zhao2017optimal}, \cite[Lem. 4.2]{perez2017optimality}
$$\fr{\r }{\q +\r}\le(\ol Z_\q(b)- \fr p q \ri)+  \fr{Z_{\q,\r}(b)}{\Fqr}+K=0.$$
This has a unique positive root if and only
if $\fr{\r }{\q +\r} \fr p q > \fr 1{\Fqr}+K.$

For \SLG dividends with costs $k x +K$ for a capital injection of $x$, and with Parisian observation,   the  value function $V(x)$ \cite[Thm. 4.1]{zhao2017optimal}
is obtained by choosing a level $a$ for \ci \  and a barrier $b$, such that $V(a)=V(0)+ k a + K, V'(a) =k, V'(b)=1$.  This yields \cite[(4.10)]{zhao2017optimal}
$$\bc Z_{\q,\r}(b-a)=k\\\fr{\r}{\q +\r}(\ol Z_\q(b)- \ol Z_\q(b)) + \fr {Z_{\q,\r}(b)-Z_{\q,\r}(b-a)}{\Fqr}=k a +K\ec.$$

\section{Examples 
}\label{s:ex}
\ssec{
Brownian motion with drift\label{s:BM}}
For Brownian motion
with drift  $X_t = \s B_t + \mu  t$, $\mu\neq 0$ (a possible model for small claims), $\k(\th) =\mu \th + \fr{\s^2}2 \th^2$ and let $\g=\fr{2 \mu}{\s^2}$ be the adjustment coefficient. The roots of $\k(\th) - \q=0$ are $\rho_{1}=(-\mu+\Dis)/\s^2=\F(\q)$ and $\rho_{2}=(-\mu-\Dis)/\s^2$ where $\Dis = \sqrt{\mu^2 +
2\q\s^2}$. The $W$ scale function is
\begin{equation}\label{BMWsf}
W_\q(x) =\frac{1}{\Dis} [
e^{\rho_{1} x}-e^{\rho_{2} x}]=\frac{1}{\Dis} [
e^{(-\mu+\Dis)x/\s^2}-e^{-(\mu+\Dis) x/\s^2}]= \fr{2 e^{- \mu x/\s^2}}{\Dis} \sinh( x \Dis/\s^2)
\end{equation}
and
$$\ovl W_\q(x)=\bc \frac{1}{\Dis} [
\fr{e^{\r_{1} x}}{\r_{1}} -\fr{e^{\r_{2} x}}{\r_{2}}-\fr{\Dis}{ \q}],& \q >0\\
\frac{1}{\mu} [x -\fr{1-e^{-\g x}}{\g}], &\q= 0
\ec.$$

The second scale function for $x \geq 0$ is:
$$Z_{\q}(x,\th)= Z_\q(x)+ \th \fr{\s^2}2 W_\q(x)=\frac{\q- \k(\th)}{\Dis}\left[\frac{e^{ \rho_{1} x}}{
\rho_{1}-\th}-\frac  {e^{\rho_{2} x}}{\rho_{2}-\th}\right].$$
One may check that for every $q$
\bea \D_{\q,\th}^{(ZW)}(x,x)=\fr 2{\s^2} e^{- \g x}, \;
\D_q^{(W)}(x)=(W_\q')^2(x)- W_\q(x)W_\q''(x)=\fr 4{\s^4} e^{- \g x},\; \Lambda_0(x):=\fr{W_0''(x)}{\D_0^{(W)}(x)}=-\mu.
\eea

Finally, the general result for reflected stopping times
\eqr{bdruinrefl} yields, after some
 symbolic algebra manipulations,    to
\be \Rui_\q^{b]}(x)=e^{-x \frac{\mu }{\sigma ^2}}\frac{ H(b-x)}{H(b)}, \; H(x)=\sqrt{2 \q  \sigma
   ^2+\mu ^2} \cosh \left(\frac{x \sqrt{2 \q  \sigma
   ^2+\mu ^2}}{\sigma ^2}\right)-\mu  \sinh \left(\frac{x
   \sqrt{2 \q  \sigma ^2+\mu ^2}}{\sigma ^2}\right)\ee
   see also \cite[Thm 1.1]{Ebe} for a  proof using martingale stopping.

\beXa Theorem \ref{p:time} becomes with $x>0$:
 \BEN[(1)] \im  the expected time to ruin when $\mu <0$ is
\beq \E_x
\le[ \tz \ri] =W(x) /\F(0)-\ovl W(x)=\frac{1}{-\g \; \mu} [
1-e^{-\g x}]-\frac{1}{\mu} [x -\fr{1-e^{-\g x}}{\g}]=-\frac{x}{\mu}.\eeq
{We can also check, as is well known}, that the last result holds asymptotically for any \lev process with $\k'(0) <0$,  i.e. that $\lim_{x\to \I}
\frac{\E_x [\tz ]}{x}
=-\fr 1{\k'(0)}.$

\im When $\mu > 0$, using $W^{*,2}(x)=\mu^{-2} \le(x(1+ e^{-\g x}) -2 \fr{1- e^{-\g x}}\g\ri), $ we find that the  expected time to ruin conditional on ruin occurring is:
\bea &&\E_x
\le[ \tz  \; \1_{\{\tz < \I\}} \ri]=\fr{\k''(0)}{2 \k'(0)} W(x) +\k'(0) W^{*2}(x) - \ovl W(x)\\&&= \frac{1}{ \mu \; \g} [
1-e^{-\g x}]-\frac{1}{\mu} [x -\fr{1-e^{-\g x}}{\g}]+\mu^{-1} \le(x(1+ e^{-\g x}) -2 \fr{1- e^{-\g x}}\g\ri)\no\\&&=\frac{x}{\mu} e^{-\g x},\eea
with maximum at $x^*=\g^{-1}=\fr{\s^2}{ 2 \mu}=\fr{  \k''(0)}{2  \k'(0)}$.

{This value   furnishes a reasonable initial
reserve, also since it coincides with the expected global infimum of a risk process started at $x^*$ is $0$. Indeed, assuming $\k'(0) >0$ and  differentiating the Wiener-Hopf factorization
$\E_0[e^{s \und X_{\I}}]=\k'(0) \fr {{s}}{{\k(s)}}$ yields
\bea \E_0[\und X_{\I}]=\k'(0) \lim_{s \to 0} \fr{\k(s) -s \k'(s)}{\k(s)^2}=
\k'(0) \lim_{s \to 0} \fr{ -s \k''(s)}{2 \k(s) \k'(s)}=\fr{ - \k''(0)}{2  \k'(0)}.\eea
}

\EEN
\eeXa


\beXa
 Optimizing the barrier under the classic de Finetti objective Theorem \ref{l:div} A) amounts to minimizing
$$W_\q'(x)=\frac  {1}{\s^2 \Dis}
\left[(\Dis-\mu  )
e^{(\Dis-\mu  )x/\s^2}+(\mu   +\Dis)e^{-(\mu  +\Dis)x/\s^2}\right].$$
Now the scale function verifies that
\beq \la {ha}
\fr{\s^{2}} 2 W_\q^{\prime \prime}(x) = \q W_\q(x) - \mu  W_\q'(x).\eeq
From this, it follows that if $\mu>0$, then
$b^*$ satisfies
\beq \la{go} W_\q(b^*)/W_\q'(b^*)=\mu/\q, \eeq and is explicitly given by \cite{gerber2004optimal}
\beq \la{dGSc}
e^{\fr{ 2 b^* \,\Dis }{\s^2}} =\left(\frac{\Dis + \mu}{\Dis - \mu}\right)^2 \Lra b^* = \fr{\s^2}{\Dis}\log \left(\frac{\Dis + \mu}{\Dis - \mu}\right)=\fr 2{\r_{1}-\r_{2}}\log\Big(\fr{-\r_{2}}
{\r_{1}}\Big)>0.
\eeq

Furthermore, as shown  by Jeanblanc and Shiryaev
\cite{jeanblanc1995optimization}, for $\mu >0$ it holds that $\frac{\s^2}{2} \left(V^{b^*]}\right)''(x) + \mu  \left(V^{b^*]}\right)'(x) -
\q \left(V^{b^*]}\right)(x) < 0$ for $x>b^*$, and this implies
that $\pi_{b^*}$ is the optimal strategy (among all admissible strategies).

If $\mu\leq 0$ on the other hand,
$W_\q'(x)^{-1}$ attains its maximum over $[0,\I)$ in $x=0$, and $b^*=0$ is optimal.
\end{Exa}

\beXa {\bf Optimal de Finetti dividends barrier until Parisian ruin}.
Recall the equation \eqr{bFP}
\bea \fr{ \Fqr }{\l} Z_{\q}(b,\Fqr) -  W_q(b) =\fr{ W_\q'(b)}{ \Fqr }.\eea

For Brownian motion, this yields

\bea { \Fqr } \left[\frac{e^{ \r_{1} x}}{
\r_{1}- \Fqr }-\frac  {e^{\r_{2} x}}{\r_{2}- \Fqr }\right] -  [
e^{\r_{1} x}-e^{\r_{2} x}] =\fr{  [\r_1
e^{\r_{1} x}-\r_2 e^{\r_{2} x}]}{ \Fqr } \Lra \eea
\bea
e^{\fr{ 2 \,\Dis  b^*}{\s^2}} = \left(\frac{\r_2}{\r_1}\right)^2\fr{ \Fqr -\r_1}{ \Fqr -\r_2} \Lra b^* = \fr 1{\r_{1}-\r_{2}}\log\Bigg((\fr{\r_{2}}
{\r_{1}})^2 \fr{ \Fqr -\r_1}{ \Fqr -\r_2}\Bigg)>0,
\eea
which converges to \eqr{dGSc} when $\l \to \I$.
\eeXa

\beXa The  SLG objective Theorem \ref{l:div}  B) is studied in \cite{lokka2008optimal,APP}.
The candidate optimal barrier \eqref{HS} will satisfy $k \D_\q^{(ZW)}(b)=W_\q'(b)$, which simplifies here to
$$\cosh( x \Dis/\s^2)- \fr \mu \Dis \sinh( x \Dis/\s^2)=k e^{-x \mu/\s^2}.$$

\end{Exa}

\subsection{Scale computations for processes with rational Laplace exponent}\label{subs:general}

Generalizing the previous example, we now assume the Laplace exponent is a rational function and that the equation $\k(\th)-\q = 0$ has distinct real roots $\r_\q^{(i)}$. From the partial fraction expansion of $1/(\k(\th)-\q)$, we easily obtain the $W$ scale function
$$
W_\q(x) = \sum_i A_i e^{\r_\q^{(i)} x}, \; \q>0
$$
where $A_i= 1/\k'(\r_\q^{(i)})$.  \Frt
$$\ovl W_\q(x) = \sum_i A_i \fr{e^{\r_\q^{(i)} x}-1}{\r_\q^{(i)} } = \sum_i A_i \fr{e^{\r_\q^{(i)} x}}{\r_\q^{(i)} } - \frac{1}{\q},
$$
by using $ \sum_i  \frac{A_i}{\th - \r_\q^{(i)}}=\fr{1}{\k(\th ) - \q }$ with $\th=0$. Then, from \eqref{oldZ} and \eqref{Z1}
$$ Z_\q(x) = \q  \sum_i A_i \fr{e^{\r_\q^{(i)} x}}{\r_\q^{(i)} }, \quad Z_{\q}^1(x) = \q \sum_i A_i \fr{e^{\r_\q^{(i)} x}}{(\r_\q^{(i)})^2} - \k'(0) \sum_i A_i \fr{e^{\r_\q^{(i)} x}}{\r_\q^{(i)} },$$
where $Z_{\q}^1(0)=0$ holds since $ \sum_i  \frac{A_i}{(\th - \r_\q^{(i)})^2}=\fr{\k'(\th )}{(\k(\th ) - \q)^2 }$ with $\th=0$
implies $ \sum_i  \frac{A_i}{(\r_\q^{(i)})^2}=\fr{\k'(0 )}{ \q^2 }$.
Similarly, from \eqref{Z} we obtain
\begin{align} \la{ZPh}
Z_\q(x,\th ) &= e^{\th x} + (\q- \k(\th ))\sum_i A_i \frac{e^{\r_\q^{(i)} x}-e^{\th x}}{\r_\q^{(i)}-\th } =
( \k(\th )-\q)\sum_i  \frac{A_i}{\th -\r_\q^{(i)}} e^{\r_\q^{(i)} x}\no \\
&= Z_\q(x)+ \th \sum_i {A_i} \fr{\fr{\k(\th)} \th- \frac{\q }{\r_\q^{(i)}}}{\th-\r_\q^{(i)}} {e^{\r_\q^{(i)} x}}.
\end{align}
For $\q=0$ the formulas are slightly different due to the fact that zero is one solution of $\k(\th) = 0$.

 \subsection{Cram\'{e}r-Lundberg model with exponential jumps \la{s:exp}}

We analyze now  the Cram\'{e}r-Lundberg model
 with exponential jump sizes with mean $1/\mu$, jump
rate $\lambda$, premium rate $c>0$,
and Laplace exponent
$\k(\th)=\th \le(c-\fr{\lambda}{\mu+\th}\ri)$, assuming $\k'(0) = c- \fr{\lambda} {\mu} \neq 0$. Let $\g=\mu - \l/c$ denote the adjustment coefficient, and let $\r=\fr \l{c \mu}$. Solving $\k(\th)-\q=0$ for $\th$ yields two distinct solutions $\r_{2} \leq 0 \leq \r_{1}=\Phi(\q)$ given by
\begin{align*}
\r_{1} =& \fr{1}{2c} \left(- \left(\mu c -\lambda - \q\right) + \sqrt{\left(\mu c -\lambda - \q \right)^2 + 4\mu \q c} \right),\\
\r_{2} =& \fr{1}{2c} \left(- \left(\mu c -\lambda - \q\right) - \sqrt{\left(\mu c -\lambda - \q \right)^2 + 4\mu \q c} \right).
\end{align*}

The $W$ scale function and is integral are:
$$  W_{\q}(x) = A_1 e^{\r_{1}x} + A_2 e^{\r_{2}x}, \; \ovl  W_{\q}(x)=\bc
\frac{1}{\k'(0)} [x -\r \fr{1-e^{-\g x}}{\g}], &\q= 0\\
A_1 \fr{e^{\r_{1} x}-1}{\r_{1}} + A_2 \fr{e^{\r_{2} x}-1}{\r_{2}},& \q >0\ec,$$
where $A_{1} = c^{-1} (\mu + \r_{1}) (\r_{1}-\r_{2})^{-1} = 1 / \k'(\r_{1})$ and $A_{2} = - c^{-1} (\mu + \r_{2})(\r_{1}-\r_{2})^{-1} = 1 / \k'(\r_{2})$.
Using the general results of the previous example, we find
\be
  Z_{\q}(x) =\bc 1&q=0\\ {\q} \left(\frac{A_1}{\r_{1}}
e^{\r_{1}x} + \frac{A_2}{\r_{2}}
e^{\r_{2}x}\right)= - \frac{c}{\mu} \left( \r_{2} A_1
e^{\r_{1}x} + \r_{1} A_2
e^{\r_{2}x} \right).&\q >0\ec.\ee

By tedious simplification of  \eqr{ZPh}, we find that
\begin{equation} \label{ZVe} Z_{\q}(x,\th) =
 Z_{\q}(x)+ \fr{\lambda}{c} \frac {\th }{\th + \mu} \frac{ e^{\r_{1}x}-
e^{\r_{2}x}} {\r_{1}-\r_{2}}, \; Z^1_\q(x)=\bc \r \fr{1-e^{-\g x}}{\g}, &\q= 0\\
   \fr{\lambda}{\mu c}  \frac{ e^{\r_{1}x}-
e^{\r_{2}x}} {\r_{1}-\r_{2}}&q >0 \ec
\end{equation}

\beXa Theorem \ref{p:time} becomes:
 \BEN[(1)] \im  When $\k'(0) <0$, we have $\Phi(0)=\z_0^{(1)}=-\g$ and hence
\begin{equation} \E_x
\le[ \tz\ri] = - \fr{1}{\g} W(x) - \ovl W(x) = - \fr{1}{\g c (1-\rho)} \left(1-\rho e^{-\g x} \right) - \frac{1}{\k'(0)} \left(x -\r \fr{1-e^{-\g x}}{\g}\right) = - \fr{x}{\k'(0)} - \fr{1}{\gamma}.\end{equation}

\im When $\k'(0) >0$, using $W^{*,2}(x)=\frac{\gamma x - 2 \rho }{\k'(0)^2 \gamma
   }+\frac{e^{-x \gamma } \rho  (\gamma
   \rho x +2)}{\k'(0)^2 \gamma }, $ we find that the  expected time to ruin conditional on ruin occurring is:
$$ \E_x \le[ \tz  \; \1_{\{\tz < \I\}} \ri]=\fr{\k''(0)}{2 \k'(0)} W(x) +\k'(0) W^{*2}(x) - \ovl W(x)= \frac{\rho }{c^2 \gamma} e^{-\g x}(\lambda x+c),$$ with maximum at $x^*=\frac{1} {\g }(2 -\r^{-1})$.
This value   furnishes a possible lower bound for the initial
reserve, which is positive if and only if $c < 2 \fr \lambda \mu \Eq p <  \fr \lambda \mu$.

\EEN
\eeXa
\beXa \la{ex:expb}
Let us recall now that  the function $ W_{\q}'(x)=H_D(x)^{-1}$ is
unimodal  with global minimum at
$$b^* = \frac{1}{\r_{1} - \r_{2}}
\begin{cases}\log
\frac{(\r_{2})^2(\mu +\r_{2})}{(\r_{1})^2(\mu +\r_{1})} \quad &\text{if $ W_{\q}''(0) <  0 \Eq (\q+\r)^2$}-
c\lambda\mu < 0\\ 0 & \text{if $ W_{\q}''(0)\geq 0 \Eq
(\q+\r)^2- c\lambda\mu \geq 0$}\end{cases}$$ since
$ W_{\q}''(0) \sim
{(\r_{1})^2(\mu +\r_{1})}-(\r_{2})^2(\mu +\r_{2})/(\r_{1})-\r_{2})
= (\q+\r)^2- c\lambda\mu$ (see also \eqr{e:secder}).  Furthermore, the optimal strategy is
always the barrier strategy at level $b^*$ \cite{APP}.

\eeXa

\ssec{Numerical optimization of dividends for the Azcue-Muller example\la{s:num}}

Consider the Cram\'er-Lundberg model perturbed by Gaussian component,
$X_{t}=x+ct-\sum_{i=1}^{N^{(\lambda)}_{t}}C_{i}+\sigma B_{t}$, where $C_{i}$ are iid pure Erlang claims, $E_{2,1}$ of order $n=2$ and $N^{(\lambda)}$ is an independent Poisson process with arrival rate $\lambda.$  The Laplace exponent  is
$\k (\theta) = c\theta - \lambda  + \lambda (\frac{\mu}{\mu+\theta})^2 +\frac{\sigma^2}{2} \theta$, and the equation $\kappa(\theta)-\delta=0$ has four roots. In what follows, the choice of parameters will be such that these roots are distinct. Since $\kappa$ is a rational function, the results of Subsection \ref{subs:general} can be used to obtain scale functions.

The interest in this example was awakened by Azcue and Muller~\cite{AM05}, who showed that the  barrier dividend strategy is not optimal  for certain parameter values. It was shown later that this is the case when the  barrier function has two local maxima, and  the last one is not the global maximum -- see \cite[Fig.1]{Loef08}.

 {It is natural to ask whether the barrier function \eqref{HDP} 
 can have the property of multi-modality which complicates the management of dividends. We did not find any such example in our experiments presented below}.

We present now some numerical experiments using a choice of parameters  close to \cite{Loef08}, namely
$\mu=1$, $\lambda=10$, $c =\frac{107}{5}$ and $\q =\frac{1}{10}$.
We consider
$\sigma=1.4$ and $\sigma=2$ as given in \cite{Loef08}. Note that, with these choice of parameters and in the absence of Brownian component, this example corresponds to the example given by Azcue and Muler~\cite{AM05} for which sufficient conditions for optimal barrier strategy do not hold.

Concerning the performance of barrier strategies under the model given above, see Figure \ref{fig:hdp1} and Figure \ref{fig:hdp2}, where we provide typical plots of the barrier function \eqref{HDP} of \eqref{DP}, for  different values of $\vt >0, \q >0, \th <0.$  Recall that, for $\th=0,$ \eqref{HDP} reduces to \eqref{th0} which is the barrier function of \eqref{e:DT}. \Fr  plots of \eqref{th0} are presented in Figure \ref{fig:hdt1} and Figure \ref{fig:hdt2}.
\begin{figure}[hbp]
\begin{center}
\includegraphics[width=6cm,height=4cm]{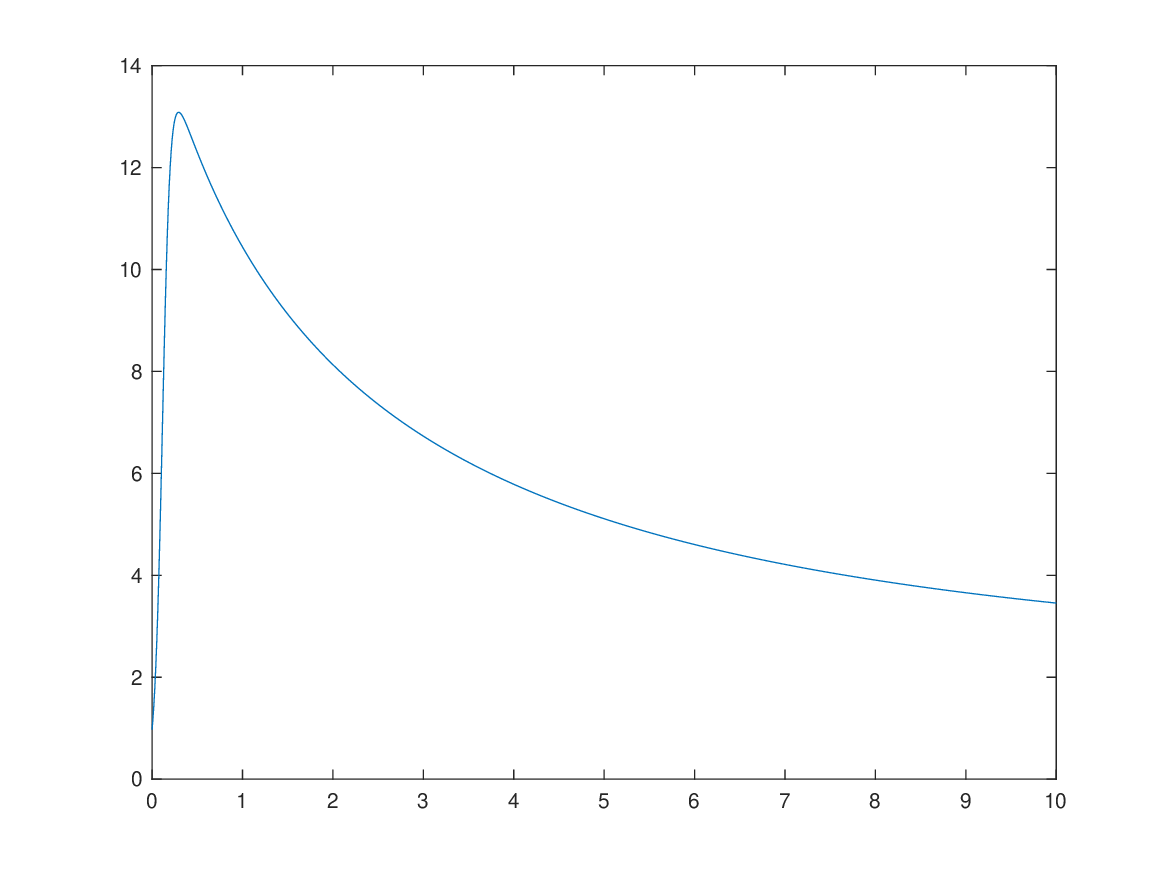}
\includegraphics[width=6cm,height=4cm]{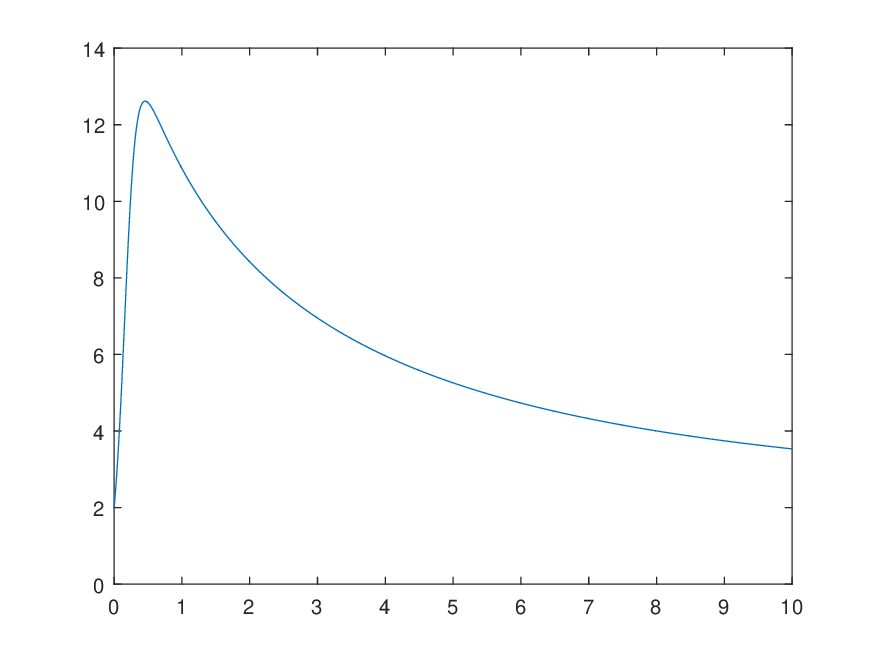}
\end{center}
\caption{Left: $\sigma=1.4$, $\theta=-0.01$, $\vartheta=1$ with $H_{DP}(0)=0.98$, $H_{DP}(\I)=2.5544$, Right: $\sigma=2$, $\theta=-0.01$, $\vartheta=1$ with $H_{DP}(0)=2$, $H_{DP}(\I)=2.5821$}
\label{fig:hdp1}
\end{figure}
\begin{figure}[htp]
\begin{center}
\includegraphics[width=6cm,height=4cm]{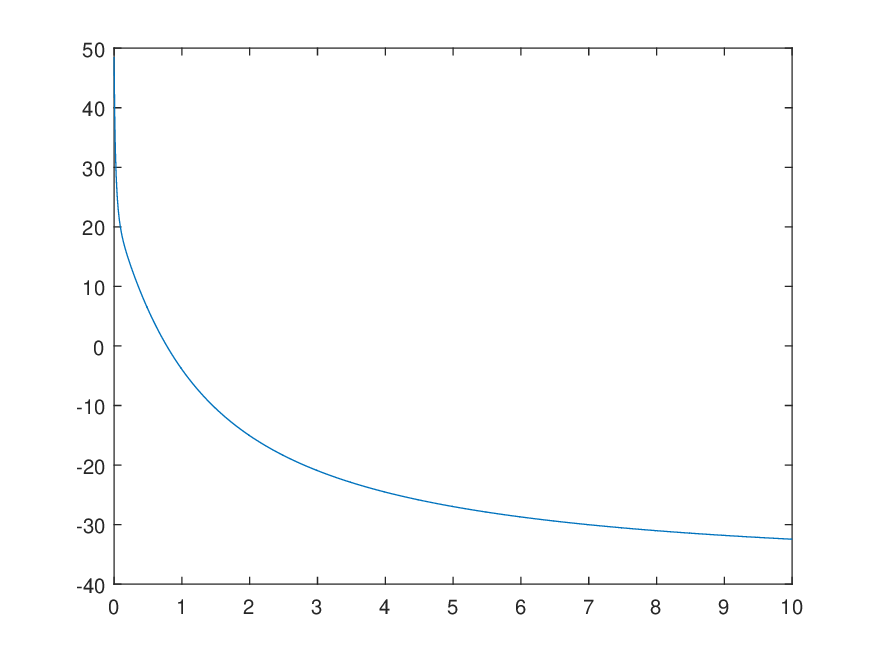}
\includegraphics[width=6cm,height=4cm]{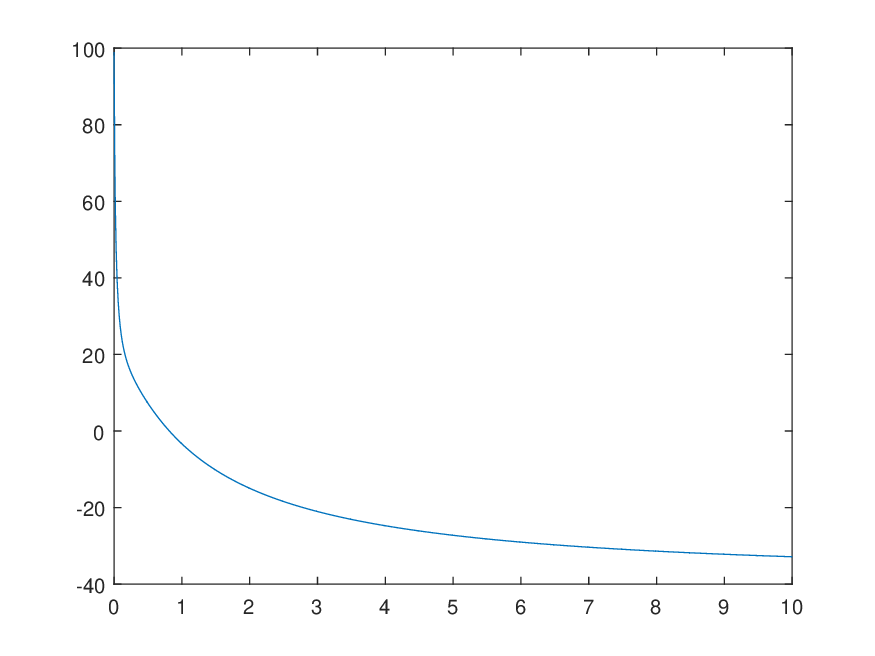}
\end{center}
\caption{Left: $\sigma=1.4$, $\theta=-0.5$, $\vartheta=50$ with $H_{DP}(0)=49$, $H_{DP}(\I)=2.5544$, Right: $\sigma=2$, $\theta=-0.5$, $\vartheta=50$ with $H_{DP}(0)=100$, $H_{DP}(\I)=2.5821$.}
\label{fig:hdp2}
\end{figure}
\begin{figure}[h]
\begin{center}
\includegraphics[width=6cm,height=4cm]{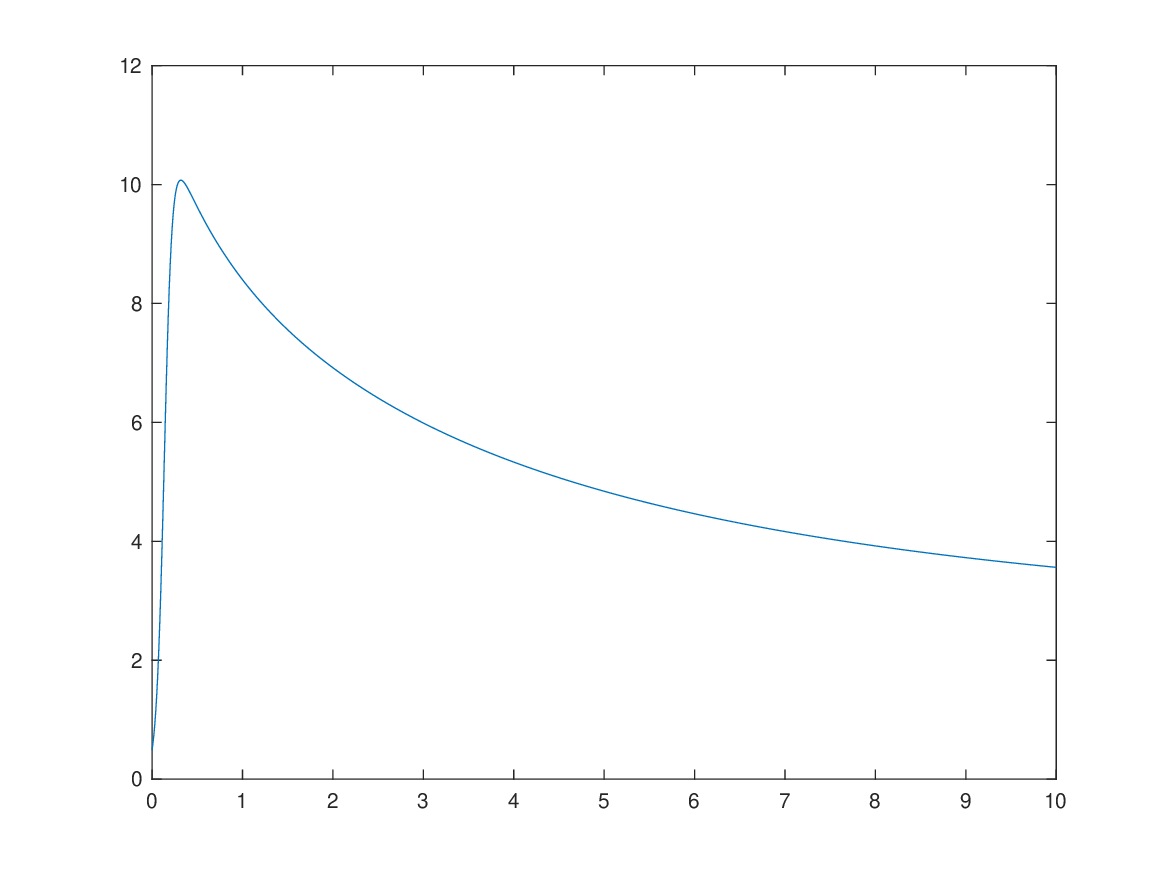}
\includegraphics[width=6cm,height=4cm]{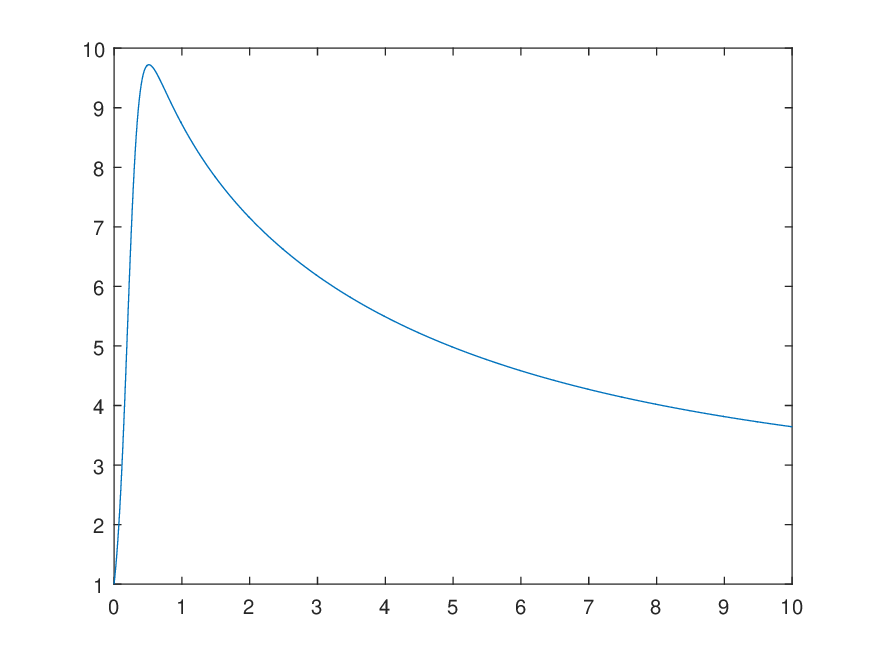}
\end{center}
\caption{Left: $\sigma=1.4$, $\vartheta=0.5$ with $H_{DT}(0)=0.49$, $H_{DT}(\I)=2.5544$, Right: $\sigma=2$, $\vartheta=0.5$ with $H_{DT}(0)=1$, $H_{DT}(\I)=2.5821$}
\label{fig:hdt1}
\end{figure}
\newpage
\begin{figure}[h]
\begin{center}
\includegraphics[width=6cm,height=4cm]{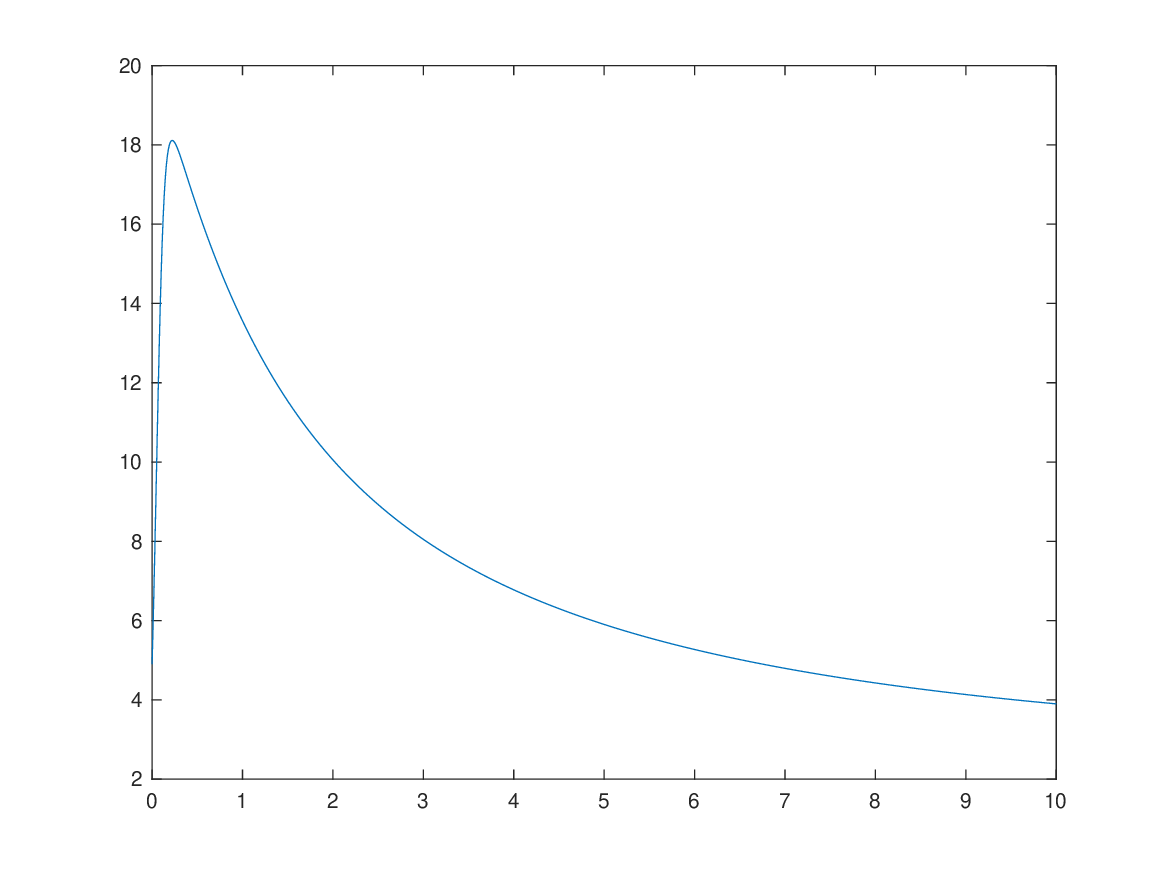}
\includegraphics[width=6cm,height=4cm]{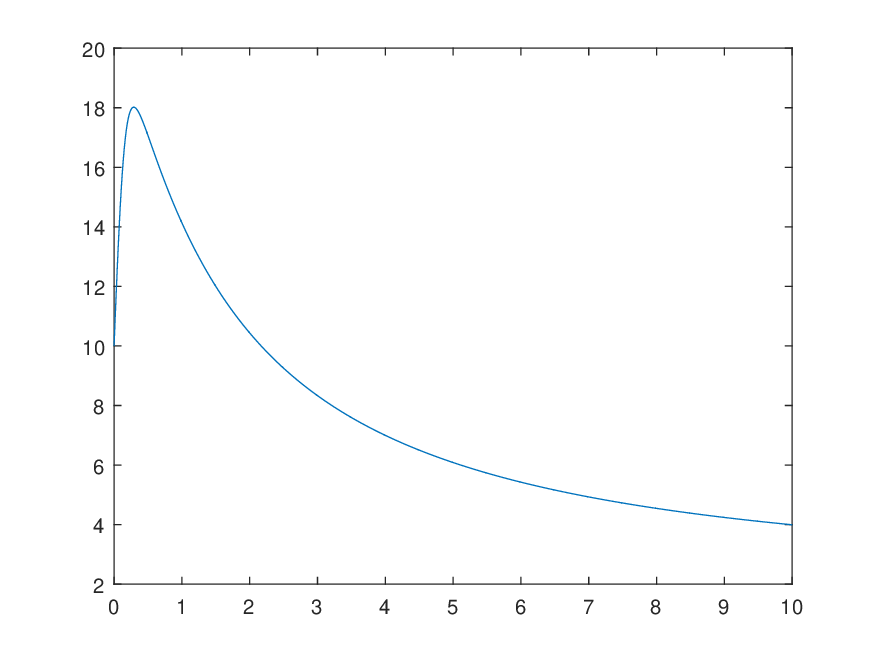}
\end{center}
\caption{Left: $\sigma=1.4$, $\vartheta=5$ with $H_{DT}(0)=4.9$, $H_{DT}(\I)=2.5544$, Right: $\sigma=2$, $\vartheta=5$ with $H_{DT}(0)=10$, $H_{DT}(\I)=2.5821$}
\label{fig:hdt2}
\end{figure}

\sec{Strong Markov processes with generalized drawdown  stopping \la{s:dd}}

In this section, $X_{t}$ will denote a one dimensional strong Markov process  {without positive jumps},
defined on a filtered probability space $(\Omega,  \{\mathcal{F}_t\}_{t \geq 0}, P)$.

Since many results  for spectrally negative \lev and diffusion \procs \; require  not much more than the \str , it was natural  to  attempt to extend such results
to spectrally negative strong Markov processes. As expected, everything worked out almost smoothly  for ``\lev -type cases'' like  random walks \cite{AV},   Markov additive processes \cite{IP}, L\'evy
processes with $\Omega$ state dependent killing \cite{IP}, and there are also some results for the more challenging case  of L\'evy processes with state dependent drift \cite{CPRY}.  In fact, the existence of some functions $W,Z$ satisfying \eqr{sRui}, \eqr{B} is clear in general, by smooth crossing and the \str.  However,
prior to the pioneering \cite{LLZ17b}, the classic and  drawdown  first passage  literatures   were restricted  mostly to parallel  treatments of the two particular cases of  diffusions and of
  spectrally negative \lev processes. \cite{LLZ17b}  showed that a direct
unified approach (inspired by \cite{Leh} in the case of diffusions) may achieve the same results for all time homogeneous Markov processes.

The crux of the approach  is to replace
$W,Z$ in the state dependent case by differential versions $\nu$ and $\de$, which were denoted in  \cite{LLZ17b} by $b,c$, in the context of the study of  drawdowns. Later,
 in \cite{ALL}, they were extended to generalized drawdown  times (which include first passage times). As will be clear from the discussion below, $\nu$ and $\de$ capture the behavior of excursions  of the process away  from its running maximum. Note however that $\nu$ is a measure, and determining when it admits a density  requires quite different  technical treatments for spectrally negative \lev and diffusion processes (see for example \cite[Lem.~8.2]{Kyp} which relates this to the challenging issue of the  differentiability of $W_q$); note also  that computing
 $W,Z,\nu,\de$  is still an open problem, even for simple classic processes like the \OU
process and the Feller branching diffusion with jumps. \cite{LLZ17b} (and \cite{ALL}) cut through this Gordian node by  restricting  to processes for which the limits defining $\nu,\de$  exist  -- see Assumptions \eqr{e:b}, \eqr{e:c},
and leaving to the user's responsibility   to  check this for their process; they also showed that  the known results   for diffusions and
  spectrally negative \lev processes were just particular cases of their general formulas -- see Section \ref{s:LLZ}.

  The results of \cite{LLZ17b,ALL}
 provide a unifying umbrella for \lev processes, diffusions, branching
\procs~(including with immigration), logistic branching
processes,  etc, under the caveat that beyond the \lev and diffusion cases, the user must establish the validity of Assumptions \eqr{e:b}, \eqr{e:c} and manage computing $\nu, \de$.

The end result is that for non-homogeneous spectrally negative Markov processes with classic first passage stopping we may provide extensions of the \ts exit equalities \eqr{sRui}, \eqr{Rui} and similar, involving now scale functions with one more variable
\begin{equation}\la{excursion} \sRui_{\q}^{b}(x,a)   =\fr{W_{\q}(x,a)}{W_{\q}(b,a)}, \quad  \Rui_{\q}^{b}(x,a)=Z_{\q}(x,a,\th) - W_\q(x,a) \frac{Z_{\q}(b,a,\th)}{W_\q(b,a)} .\end{equation}
For diffusions for example, $W_{\q}(x,a)$ is a certain Wronskian (see  \cite{Bor}) and for Langevin type processes with decreasing state-dependent drifts, $W_{\q}(x,a)$ solves  a certain renewal equation \cite{CPRY}.
So, formally the spectrally negative \Mar case is similar to the \lev one, up to adding one variable to the \fund~functions.  

Extensions to drawdown  stopping are possible as well \cite{LLZ17b,ALL}, but they are easier to state in terms of differential exit parameters $\nu,\de$ defined in \eqr{e:b}, \eqr{e:c} below.
 Before reviewing these extensions, we will introduce some objects of interest via  an illustrative example of  first passage \prob  for $(X,Y)$, with $Y$ a drawdown  process.
 In this case, simple geometric arguments (see Figure \ref{f:samplepath2})  reduce  the computation of   Laplace transforms of   exit times  of $(X,Y)$ from rectangles to those of simpler Laplace transforms  defined in \eqr{UbDL}, \eqr{DbUL}, which seem to be \fund~to this setup.

\ssec{Joint evolution  of a strong Markov process and its drawdown in a rectangle\label{s:strong}}

 In order to study  the process $(X,\Y)$, it is convenient to start with its evolution    in a rectangular region $R := [a,b] \times [0,d] \subset {\mathbb R} \times {\mathbb R}_+$,  where $a < b$ and $ d > 0$.

 A sample path of $(X,\Y)$, where $X$ is chosen to be the
standard Brownian motion,
and the region $R$ is depicted in Figure \ref{f:samplepath2}.
\figu{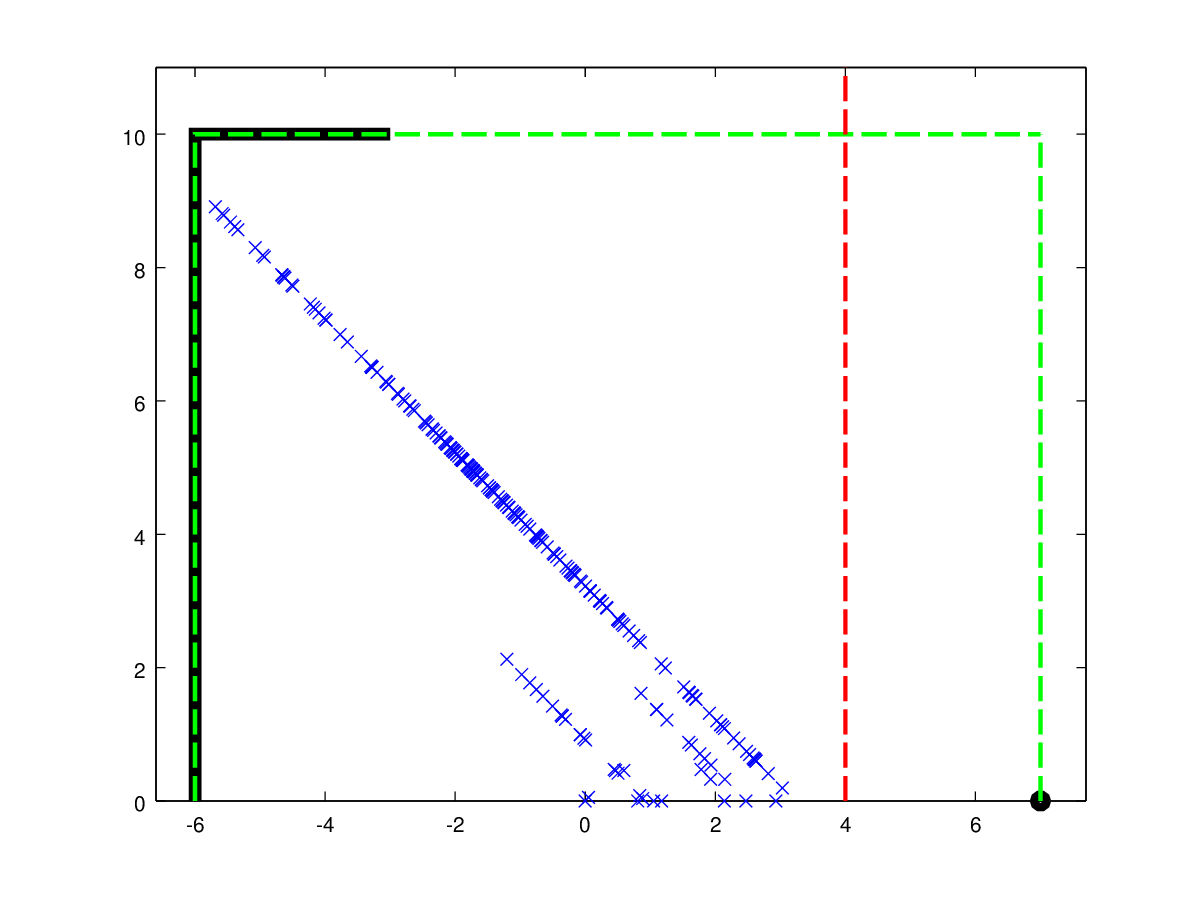}{A sample path of $(X,\Y)$ (sampled at
time step $\Delta t = 0.1$) when $X$ is a standard Brownian motion
with $X_0 =0.2$,
and the region $R$ with $d=10$, $a=-6$ and $b=7$. The dark
boundary shows the possible  exit points of $(X,\Y)$
from $R$. The base of the red line separates $R$ in two parts with different behavior}{0.8}

\beR \la{r:snM}

 As suggested by Figure \ref{f:samplepath2}, the study of the process $(X,\Y)$ may be reduced to one-dimensional problems:
 \BEN \im On the $y=0$ axis, we observe the maximum process  $\ovl X$. If \frt  downward excursions are excised, we obtain the so-called  upward ladder process $\T X(m)=m$ (the maximum studied as a function of itself), which is of course Markovian with generator $\pd{}{m}$. If \frt  time killing is present,  $\T X(m)$ becomes a killed drift subordinator, with \LE $\k(s)=s + \F_q$ (as a consequence of the \WH decomposition \cite{Kyp}).

 \im Away from the boundary
$y = 0$, the process
 oscillates during  negative excursions from the maximum on  line segments $l_{\ol{X}_t}$ where, for $c \in {\mathbb R}$,
$l_c := \{x \in {\mathbb R} \times {\mathbb R}_+: x_1 + x_2 = c \}$. Since $\ol{X}_t$ is fixed during such an \exc , we are dealing here essentially with the process $-X_t$.

 \im If the first excursion outside the rectangle kills the process, the ladder process becomes  a killed drift with  generator $\mG \f(m):=\f'(m)- \nu_{q}(d)\f(m)$ \cite{AACI,AVZ}, since the killing excursions are a Poisson process with rate $\nu_{q}(d)$.

  \im With generalized drawdown  defined in the next subsection (when  the upper boundary is replaced  by one determined by certain parametrizations $(\H d(m),d(m)), \H d(m)=m-d(m)$),
  the  generator of $\T X_m$ will have state dependent  killing:
\be\label{eq:gen}
   \mG \f(m):=\f'(m)- \nu_{q}(d(m))\f(m). \end{equation}\im Finally, in the spectrally negative \Mar case,
the  generator becomes:
\be\label{eq:genM}
   \mG \f(m):=\f'(m)- \nu_{q}(m,\H d(m))\f(m), \end{equation}   where  the killing rate
   \beq \nu_{q}(m,y)|_{y=\H d(m)}\eeq
  depending  of both the current position and the killing limit $y=\H d(m)$ is defined in \eqr{e:b} below.
\EEN
\eeR

     The fact that many functionals (ruin, dividends, tax, etc) of the original process may be  expressed
   as functionals of the killed ladder  process explains the prevalence of first order ODE's related to the generator \eqr{eq:gen} when working with spectrally negative processes.

 We see from the remarks above  that   $\nu_{\q}$ may serve as a more convenient alternative characteristic of a spectrally negative Markov process, replacing $W_{\q}$, and that  it may be used also in the case of  generalized drawdown  killing.

Define now
	\[
	T_R =T_{a,b,d}:=  \inf\{t: (X_t,\Y_t) \notin R \}= \t_d \wedge T_{a,-} \wedge T_{b,+}.
	\]
Several implications for $T_R$ are immediately clear from these dynamics:
for example, the process $(X,\Y)$ can leave $R$ only through
$\partial R \cap \{x \in {\mathbb R}\times {\mathbb R}_+: x_1 \leq b-d \}$
or through the point $(b,0)$ (see the shaded region in Figure
\ref{f:samplepath2}). Also,

\BEN \im
If $b-a \leq d$, it is impossible
for the process to leave $R$ through the upper drawdown  boundary of $\partial R$
and for these parameter values $T_R$ reduces to
$T_{a,-} \wedge T_{b,+}$. Here it suffices to know the survival/ruin functions \eqref{sRui}, \eqref{Rui} in order to obtain the Laplace transform of $T_R$.

 \im {If $a+ d \leq x$, it is impossible
for the process to leave $R$ through the left boundary} of $\partial R$,
and  $T_R$ reduces to
$T_{b,+} \wedge \t_d$.  Here it suffices to apply the spectrally negative drawdown  formulas
provided in   \cite{MP,LLZ17}.

\im In the remaining case $x \leq a+d \leq b$, both  drawdown  and classic exits are possible. For the latter case, see Figure \ref{f:samplepath2}. The key observation here is that  drawdown  [classic] exit occurs if and only if $X_t$ does [does not] cross the line $x_1=d+a$. The final answers will combine these two cases.

\EEN

 Two natural objects of interest in ``mixed drawdown /\fp" control over the rectangle are the ``\ts exit" times
\bea \t_{b,d}=\min(\td_d, \tb), \; \t_{a,d}=\min(\td_d, \ta).  \eea
 In terms of the two dimensional process $t\mapsto (X_t,\Y_t)$,  these are  the first exit times from
 the regions $(-\infty,b) \times [0,d]$ and $(a,\infty) \times [0,d].$

 We introduce now two Laplace transforms $UbD/DbU$(standing for up-crossing before drawdown/drawdown before  up-crossing) involving the ``\ts exit'' times, which are analogues of the killed \sur and ruin \pros:
\begin{equation} \la{d1}
\begin{aligned}
UbD^{b}_{\q,d}(x) &=\E_x \left[e^{-q \tb} ; \tb <  \td_d \ri]=\E_x \left[e^{-q \tb} ; \ol X_{\td_d} > b  \ri],\\ DbU^{b}_{\q,\th,d}(x)&=\E_x
\left[ e^{-q\td_{d}-\th (Y_{\td_d}-d)}
;  \td_{d}<\tb \right]=\E_x
\left[ e^{-q\td_{d}-\th (Y_{\td_d}-d)}
;   \ol X_{\td_d} <  b \right].
\end{aligned}
\end{equation}

By using $UbD/DbU$ we provide now Laplace transforms of $T_R$ and of the eventual overshoot at $T_R$.  One can break down the analysis of $T_R$ to nine cases, depending on which of the three exit boundaries
$T_{a,-}$, $T_{b,+}$ or $\td_d$ occurred, and on the three relations between $x$, $a$, $b$ and $d$ described above. The results are then the immediate applications of the strong Markov property.

\beP \label{Thm} Consider a spectrally negative
Markov process $X$ with differentiable scale function $W_{\q}$. Then, for $d \geq 0$ and $a  \leq   x \leq b$, we have: \begin{center}
\begin{equation}\label{e:mmainresult}
{\def\arraystretch{2}
  \begin{tabular}{ r | c  |c | c}
	&$a+d \leq x \leq b$ & $x  \leq a+d \leq  b$ & $b \leq a+ d$\\
\hline
$ \E_x \left[ e^{-q \tb } ; \tb \leq \min(\td_d, \ta) \right]=$ &
$UbD_{\q,d}^{b}(x)$
&
$\sRui_{\q}^{(a+d)}(x,a) UbD_{\q,d}^{b}(a+d)$
&
$\sRui_{\q}^{b}(x,a)$
\\\hline
$ \E_x \left[e^{-q \ta + \th(X_{\ta} -a) } ; \ta  \leq \min(\td_d, \tb) \right]=$
&
$0$
&
$\Rui_{q,\th}^{(a+d)}(x,a)$
&
$\Rui^{b}_{\q,\th}(x,a)$
\\\hline
$ \E_x\left[e^{-q \td_{d}  {- \th(Y_{\td_{d}} -d) }} ; \td_d  \leq \min(\tb, \ta )\right]=$
&
$DbU^{b}_{q,\th,d}(x)$
&
$\sRui_{\q}^{(a+d)}(x,a)DbU _{\q,\th,d}^{b}(a+d)$
 &
$0$
  \end{tabular}
}
\end{equation}
\end{center}
\eeP

\prf Note that in the third column the $d$ boundary is invisible and does not appear in the results, and in the first column the $a$ boundary is invisible and does not appear in the results. These two cases follow therefore by applying already known results.

The middle column holds by breaking the path at the first crossing  of $a+d$. The main points here are that \BEN \im the middle case
may happen only if $X_t$ visits $a$ before $a+d$;
\im  the first case (exit through $b$) and the third case (drawdown  exit) may happen only if $X_t$ visits first $a+d$,
 {with the drawdown  barrier being invisible}, and that
subsequently the lower first passage  barrier $a$ becomes invisible. \EEN The results follow then due to the smooth crossing upward and the strong Markov property.\\

We will leave open the question of how to compute the drawdown functions $UbD/DbU$ until   Subsection \ref{s:LLZ} where we will consider more general drawdown boundaries. However, we note here that for spectrally negative \lev processes they have simple formulas. In the \lev case  for example
\begin{equation} \la{UbDL} UbD_{\q,d}^{b}(x)=\E_x \left[e^{-q \tb } ; \tb \leq  \td_d \right]= e^{-(b-x) \fr{W_{\q}'(d)}{W_{\q}(d)}},
\end{equation}and the function $DbU$  may be obtained by integrating the \fund law \cite[Thm.~1]{MP},
\cite[Thm.~3.1]{LLZ17}\fn[5]{Note that \cite[Thm.~1]{MP} give a more complicated "sextuple law" with two cases, and that
\cite[Thm.~3.1]{LLZ17}  use an alternative to the function $Z_{\q}(x,\th)$, so that some computing is required to get \eqr{UbDL}, \eqr{funL} and \eqr{DelZW}.}
\begin{eqnarray} \la{dem} && \de_{\q,\th}(d,x,s):= {\E_x \left[
e^{ -\q \td_d - \th (Y_{\td_d}-d)}; \ovl X_{\td_d} \in \md s
\right]}=
\Big(  \nu_{\q }(d)\; e^{-  \nu_{\q }(d) (s-x)_+}\; \md s \Big) \T \de_{q,\th}(d)\no \\&& \Eq
 {\E_x \left[
e^{ -\q \td_d - \th (Y_{\td_d}-d)- \vt (\ovl X_{\td_d} -x)}
\right]}=
\fr{  \nu_{\q }(d)}{\vt+ \nu_{\q }(d) } \T \de_{q,\th}(d), \quad \for x.
\end{eqnarray}
where $\T \de_{q,\th}(d)$ is given by \eqref{de}. Integrating \eqr{funL} yields
\begin{equation} \la{DbUL}
DbU_{\q,\th,d}^{b}(x)=\Big(1-e^{-(b-x) \fr{W_{\q}'(d)}{W_{\q}(d)}}\Big)\T \de_{q,\th}(d).
\end{equation}

Note that the \fund law reflects the independence of the path before the last maximum and after,
conditional on the value of the last maximum. The exponential law of the last maximum is due to the \lev setup, and will be lost in  the \Mar case, where it will be replaced by the law of the first arrival in a ``\nonh Poisson process of killing excursions''.

\beC In the spectrally negative \lev case, Theorem \ref{Thm} holds with the first passage and drawdown  functions given by \eqref{sRui}, \eqref{Rui}, \eqref{UbDL}, \eqref{DbUL}.
\eeC

\ssec{Generalized drawdown stopping for processes without positive jumps} \la{s:gendd}
Generalized drawdown  times appear naturally  in the \AY solution of the Skorokhod embedding problem \cite{AY}, and in the  Dubbins-Shepp-Shiryaev, and Peskir-Hobson-Egami optimal stopping problems \cite{dubins1994optimal,peskir1998optimal,hobson2007optimal,egami2015excursion}.
 Importantly, they allow a unified treatment of classic first passage  and drawdown  times -- see \cite{AVZ,LVZ} (see also \cite{ALL} for a further generalization to taxed processes).   The idea is to replace  the upper side of the rectangle $R$ by a parametrized curve
\begin{equation*}
(x_1, x_2) =  (\H d(s), d(s)), \quad \H d(s)=s-d(s),
\end{equation*}
where $s=x_1+x_2$ represents the value of $\ol X_t$ during the excursion which intersects the upper boundary at $(x_1,x_2)$ (see Figure \ref{f:plD}). Alternatively, parametrizing by $x$ yields (note $Y_t \geq  d(\ol X_t) \Eq  Y_t \geq h(X_t)$)
$$y=h(x), \quad h(x)=\H{d}^{-1}(x) -x .$$
  \ninseps{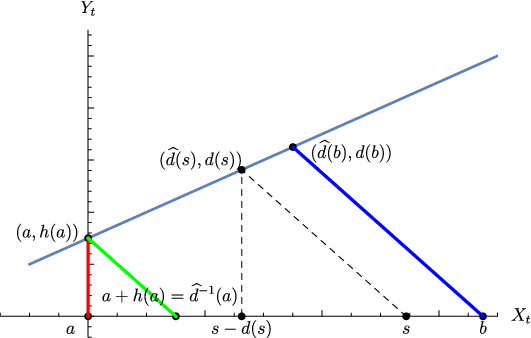}{Affine drawdown  exit of $(X,Y)$ with $a=0$,
 $d(s)=\fr 13 s+ 1$}{0.8}

  \beD \la{d:LVZ}  \cite{AY,LVZ}  For any 
   function $d(s)>0$  such that $\H d(s)=s-d(s) $ is nondecreasing, a {\bf generalized drawdown}  time is defined by
\beq \td_{\H d(\cdot)}:=\inf\{t\geq0:Y_t>  d(\ol X_t)\}=\inf \left\{t\geq 0: X_t < \H d(\ol {X}_t)\right\}.
 \label{ddg}\eeq
 \eeD

  Such times provide a natural unification of classic and drawdown  times. Introduce $$\T Y_{t}:=Y_t-  d(\ol X_t), \; {t\geq0} $$ 
   to be called  drawdown type
process.
Note that we have $\T Y_{0}=- \H d(X_{0})<0$,
and that the process $\T Y_t$ is in general  {non-Markovian}.   {However, it is
Markovian during each negative excursion of $X_t$, along one of the oblique lines
in the geometric decomposition sketched in Figure \ref{f:samplepath2}.}

\beXa  With affine functions
\begin{equation} d(s)= (1-\xi) s+ d \ \Eq \ \H d(s)= \xi s-d \ \Eq \ h(x)=\fr{(1-\xi) x+ d}{\xi}, \quad  \xi\in  {[0, 1]}, \label{aff}\end{equation}we obtain the   affine drawdown/regret times studied in \cite{AVZ}.

 Affine drawdown times reduce  to a classic drawdown  time \eqref{dd} when $\xi=1,  d(s)=  d$, and to a time of first passage below a level when $\xi=0, \H d(s)=-d,  d(s)=  s+d$. When  $\xi $ varies, we are dealing with the pencil of lines passing through $(x,y)=(-d,d)$.
In particular, for $\xi=1$ we obtain an infinite strip, and for $\xi=0, d=0,$ we obtain the positive quadrant  (this case corresponds to the classic ruin time).

One of the merits of affine drawdown times is that they allow unifying the classic first passage  theory with the  drawdown  theory \cite{AVZ}. 
 A second merit is that they are optimal
 for the variational problem considered below.

\eeXa

Introduce now  generalized drawdown  analogues of the drawdown  \sur and ruin \pros \eqref{d}, for which we will use the same notation:
\begin{align}\la{UbDdef}
UbD_{q,\H d (\cdot)}^{b}(x) &=\E_x \left[e^{-q \tb } ; \tb \leq  \td_{\H d (\cdot)} \ri],\\
DbU^{b}_{q,\th,\H d (\cdot)}(x)&=\E_x \left[ e^{-q\td_{\H d (\cdot)}-\th \T Y_{\td_{\H d (\cdot)}}} ; \td_{\H d (\cdot)}<\tb \right].\la{DbUdef}
\end{align}
An extension of Theorem \ref{Thm} to generalized drawdowns is straightforward:

\beP \label{Thmp} Consider a spectrally negative Markov process $X$ with differentiable scale function $W_q$. Then, for
$ a  \leq   x \leq b$ and $d(\cdot)$ satisfying the conditions of Definition \ref{d:LVZ},  we have:
{\small
\begin{center}
\begin{equation*}\label{e:mainresult}
{\def\arraystretch{2}
\begin{tabular}{ r | c  |c | c}
	&$a+ h(a) \leq x$ & $x  \leq a+ h(a) \leq  b$ & $b \leq a+ h(a)$\\
\hline
$  {\E_x \left[e^{-q \tb } ; \tb \leq \min(\td_{\H d(\cdot)}, \ta) \right]}=$ &
$UbD_{q,\H d(\cdot)}^{b}(x)$
&
$\sRui_q^{a+h(a)}(x,a)UbD_{q,\H d(\cdot)}^{b}(a+h(a))$
&
$\sRui_{q}^{b}(x,a)$
\\\hline
$  {\E_x \left[e^{-q \ta + \th(X_{\ta} -a) } ; \ta  \leq \min(\td_ {\H d(\cdot)}, \tb) \right]}=$
&
$0$
&
$\Rui_{q,\th}^{a+h(a)}(x,a)$
&
$\Rui^{b}_{q,\th}(x,a)$
\\\hline
$  {\E_x\left[e^{-q \td_{\H d(\cdot)} - \th(Y_{\td_{\H d(\cdot)}} -d) } ; \td_{\H d(\cdot)}  \leq \min(\tb, \ta)\right]}=$
&
$DbU^{b}_{q,\th,\H d(\cdot)}(x)$
&
$\sRui_q^{a+h(a)}(x,a)DbU_{q,\th,\H d(\cdot)}^{b}(a+h(a))$
 &
$0$
  \end{tabular}
}
\end{equation*}
\end{center}
}
\eeP

\ssec{First passage theory for upwards skip-free Markovian processes: $\nu$ and $\de$ replace $W,Z$
\la{s:LLZ}}

In this section, we review the functions $\nu_\q,\de_{\q,\th}$, essentially differential versions of the scale functions $W_\q,Z_{\q}$ of spectrally negative \lev theory, which serve to extend the spectrally negative \lev  theory to the spectrally negative \Mar case. They  were first constructed
in \cite{Leh,LLZ17b},    via
  an ``infinitesimal decomposition" approach into {\bf two sided infinitesimal exit problems} for $ X_t$ out of intervals $[x - d,x+\e]$. It was later observed in  \cite{ALL} that using  intervals $[x - d(x),x+\e]$  allows extending this to the
framework of generalized drawdown /\AY times  -- see Figure \ref{f:plD}.

The key step is assuming  the existence of differential versions of the ruin and \sur probabilities \eqref{sRui}, \eqref{Rui}:
\begin{Ass}
\label{As}For all $q,\th\geq0$ and $y\leq x  $ fixed,  assume  that $\sRui_{\q}^{b}(x,y)$  and
$\Rui_{q,\th}^{b}(x,y)$ are differentiable in $b$ at  $b=x$, and in particular that the following limits exist:
\begin{equation} \la{e:b}
\nu_{q}(x,y)   :=\lim_{\varepsilon \downarrow 0}\frac{1-
\sRui_{\q}^{x+\varepsilon}(x,y)}{\varepsilon} \quad \text{(total infinitesimal hazard rate)}
 \end{equation}and%
\begin{equation} \la{e:c}
{\de}_{q,\th}(x,y)   :=\lim_{\varepsilon \downarrow 0}
\frac{\Rui_{q,\th}^{x+\varepsilon}(x,y)}{\varepsilon} \quad \text{(infinitesimal spatial killing rate)}
\end{equation}\end{Ass}

\beR   It turns out that everything reduces
to  the differentiability of the two-sided  ruin and survival probabilities as functions of the upper limit. Informally, we may say that the pillar of first passage  theory for spectrally negative Markov processes is proving the existence of $\nu,\de$. 
\eeR

\beR
In the spectrally negative L\'evy case  \eqref{sRui}, \eqref{Rui} imply that $\nu_{q}(x,y)=\fr{W_{\q}'(x-y)}{W_{\q}(x-y)} = \nu_q(x-y)$, and ${\de}_{q,\th}(x,y)={\de}_{q,\th}(x-y)$ with ${\de}_{q,\th}(d)=\nu_\q(d)(Z_{\q}(d,\th)-W_\q(d)\fr{Z_{\q}'(d,\th)}{W_\q'(d)})=
\nu_\q(d) \T {\de}_{q,\th}(d)=$ (see \eqref{de}).
\eeR

A necessary condition for Assumption \ref{As} to hold is that,
\[
\t_{x}^{+}=0\text{ and }X_{\t_{x}^{+}}=x,\text{ }\mathbb{P}_{x}\text{-a.s. for
all }x\in\mathbb{R}\text{.}%
\]
In other words, $X$ must be upward regular\fn[4]{A process is called upward regular if $P_{y}%
(T_{x,+}<\infty)>0$, for all $y < x\in\mathbb{R}$.}
and upward creeping at every $x \in \R$.
Assumption  \ref{As} holds for processes $X$ that are  \textbf{upward skip-free}.


Assuming the existence of the limits in Assumption \ref{As}, \cite[(3.2), Thm.~3.1,Cor.~3.1]{LLZ17b} show how to compute the first passage  functions from their differential versions. The extension of this result with generalized drawdown  times  is \cite[Thm.~1]{ALL}:

\beP
\label{t:exit} Consider a Markov process
$X$ such that Assumption \ref{As} holds. Assume $d(\cdot) $ satisfies the conditions of Definition \ref{d:LVZ}, and $q,\th\geq0, b\in\mathbb{R}$.

A)  The ``upper first passage''  function \eqref{UbDdef} 
is given by
{ \begin{equation}
UbD(x)=UbD_{q,\H d (\cdot)}^{b}(x) =e^{-\int_x^b \nu_q(s,\H d (s))  d s},\la{UbDex}\end{equation}}
and satisfies the ODE
\begin{equation} \la{UbD}
UbD^{\prime}(y)-\nu_{q}(y,\H d(y))UbD(y)=0,\quad UbD(b)=1,
\end{equation}

B) The ``lower first passage''  function \eqref{DbUdef} 
is given by
{ \begin{equation} \label{DbUex}
DbU(x)=DbU^{b}_{q,\th,\H d (\cdot)}(x)=\int_x^b e^{-\int_x^y \nu_q(s, \H d (s)) \md s} \nu_q(y,\H d (y)) \de_{q, \th}({y},\H d (y)) d y,
\ee}
and satisfies the ODE
{ \begin{equation} \la{DbU}
DbU^{\prime}(y)-\nu_{q}(y,\H d(y))DbU(y)+\de_{\q,\th}(y,\H d(y))
=0,\quad DbU(b)=0.\end{equation}}
\eeP

\prf See \cite[(3.5)]{LLZ17b} for the case $ l(x)=x-d$, and \cite{ALL} for the general case.\qed

\beR

We view differential equations like  \eqr{UbD}, \eqr{DbU} as the fundamental object of spectrally negative first passage  theory,
 due to  their
 probabilistic interpretation as Kolmogorov equations of  the upward ladder process with excised negative excursions.
\eeR

\beR
In the spectrally negative L\'evy case, \eqref{UbDex} reduces by using \eqref{nu}   to
\begin{equation*}
UbD_{q,\H d (\cdot)}^{b}(x) =e^{-\int_x^b \nu_q(s,\H d (s))  d s} = e^{-\int_x^b \fr{W_{\q}'(d(s))}{W_{\q}(d(s))}  \md s},
\end{equation*}
and \eqref{DbUex} becomes
\begin{align*}
DbU^{b}_{q,\th,\H d (\cdot)}(x)&=\int_x^b e^{-\int_x^y \nu_q(d (s)) \md s} \nu_q(d (y)) \T \de_{q, \th}(d (y)) \md y\\
&=\int_x^b e^{-\int_x^y \fr{W_{\q}'(d(s))}{W_{\q}(d(s))} \md s} \fr{W_{\q}'(d(y))}{W_{\q}(d(y))} \left(Z_{\q}(d(y),\th)-W_\q(d(y))\fr{Z_{\q}'(d(y),\th)}{W_\q'(d(y))} \right) \md y.
\end{align*}
Furthermore, if we have classic drawdown    $d(s)=d\geq 0$, then we obtain \eqref{UbDL} and \eqref{DbUL}
\begin{align*}
UbD_{q,\H d (\cdot)}^{b}(x) &= e^{-(b-x) \fr{W_{\q}'(d)}{W_{\q}(d)}}= UbD_{\q,d}^{b}(x),\\
DbU^{b}_{q,\th,\H d (\cdot)}(x)&=\int_x^b e^{-\int_x^y \fr{W_{\q}'(d)}{W_{\q}(d)} \md s} \fr{W_{\q}'(d)}{W_{\q}(d)} \T \de_{q, \th}(d) \md y = \Big(1-e^{-(b-x) \fr{W_{\q}'(d)}{W_{\q}(d)}}\Big) \T \de_{q,\th}(d) = DbU^{b}_{q,\th,d}(x).
\end{align*}
\eeR

We may also express Proposition \ref{t:exit}    in terms of a generalized $W,Z$ basis.
\beR
(1) Introducing \begin{equation} W_{\q, d(\cdot)}(x,a):=e^{
\int_{a}^{x}\nu_{q}(s,\H d (s))\mathrm{d}s} \Eq \nu_{q}(x,\H d (x))=
 \fr{W_{\q, d (\cdot)}'(x)}{W_{\q, d (\cdot)}(x)},\end{equation}for some {\it arbitrary } $a  \leq x$,  we may rewrite \eqref{UbDex} as
\beq UbD^b_{q,\H d (\cdot)}(x)=\fr{W_{\q, d (\cdot)}(x,a)}
{W_{\q, d (\cdot)}(b,a)}.
\eeq

(2) We may  rewrite \eqref{DbUex} in an alternative form
\begin{equation} DbU_{q,\th,\H d (\cdot)}^b(x)=Z_{\q,d (\cdot)}(x,\th) - \fr{W_{\q, d (\cdot)}(x) }{W_{\q, d (\cdot)}(b)} Z_{\q,d (\cdot)}(b,\th).  \la{DbUZ}
\end{equation}
where we  put  ${\T \de}_{q,\th}(x,y):=\fr{\de_{q,\th}(x,y)}{\nu_{q}(x,y)}$ and
\begin{align} 
Z_{\q,d (\cdot)}(x,\th)&:={\T \de}_{q,\th}(x,\H d (x)) +W_{\q, d (\cdot)}(x)\int_{x}^{\I}\fr{{\T \de}'_{q,\th}(s,\H d (s))}{W_{\q, d (\cdot)}(s)}\mathrm{d}s \la{Zd}\\
&\Eq \ \fr{Z'_{q,d(\cdot)}(x,\th)}{W'_{q, d (\cdot)}(x)}=\int_{x}^{\I}\fr{{\T \de}'_{q,\th}(s,\H d (s))}{W_{\q, d (\cdot)}(s)}\mathrm{d}s\no\\
&\Eq \ {{\T \de}'_{q,\th}(x,\H d(x))}=- {W_{\q, d (\cdot)}(x)}\left(
\fr{Z'_{q,d(\cdot)}(x,\th)}{W'_{q, d (\cdot)}(x)}\right)'.\no
\end{align}

\eeR

\beR Note that while $\nu, \de$ are just functions of two variables, in the drawdown  framework  $W$ and $Z$ are functionals of the initial position and of the drawdown  function $d\cd$. \eeR

\prf B) It may be checked that substituting $Z_{\q,d (\cdot)}(x,\th)$ given by the first equality in \eqr{Zd} into \eqref{DbUZ} yields $DbU(x)= {\T \de}_{q,\th}(x,\H d (x))-{\T \de}_{q,\th}(b,\H d (b))\fr{W_{\q, d (\cdot)}(x)}{W_{\q, d (\cdot)}(b)}+ W_{\q, d (\cdot)}(x)\int_{x}^b \fr{{\T \de}'_{q,\th}(s,\H d (s))}{W_{\q, d (\cdot)}(s)}\mathrm{d}s$; but this is just an alternative way to express the solution of the ODE  \eqref{DbU}, obtained  by an integration by parts.\qed

\beR \la{r:refl} With   {\bf classic \fp} stopping $\H d(x)=a$, and we obtain
\begin{equation} \la{laws}\sRui_{\q}^{b}(x,a)   =\fr{W_{\q}(x,a)}{W_{\q}(b,a)}, \quad \Rui_{\q,\th}^{b}(x,a)=Z_{\q}(x,a,\th) - W_\q(x,a)\fr{Z_{\q}(b,a,\th)}{W_\q(b,a)},\end{equation} with  scale functions involving now just the  variable $a$ (the non-smooth first passage end), which reduce to the classic \lev formulas upon replacing $(x,a)$ by $x-a$.
\eeR

\beXa With {\bf fixed drawdown } stopping $ d(x)=d$,
 in the \lev spectrally negative case, it follows that $\nu_{\q}(d)=\fr{W_{\q}'(d)}{W_{\q}(d)} \Eq W_{q,d}(x)=e^{-x \nu_q(d)}$.
 We recover also the simple structure of the parameter $\de_{q,\th}$ \cite[Exa.~3.1]{LLZ17b}:
\begin{equation} \la{bW}  \de_{q,\th}(d)=
Z_{q}(d,\th) \nu_{\q}(d) -   Z_{q}'(d,\th)=\fr{W_{\q}'(d)}{W_{\q}(d)}{\T \de}_{q,\th}(d),\end{equation} with ${\T \de}_{q,\th}(d)=Z_{q}(d,\th) -   \fr {W_\q(d)}{ W'_\q(d)} Z_{q}'(d,\th),$ and \eqr{Zd} becomes    $$Z_{q,\th,d}(x)=\fr{e^{x \nu_q(d)}+{\T \de}_{q,\th}(d)}{1+{\T \de}_{q,\th}(d)}.$$

\eeXa

\vfe

\beR
Recall now that in the \lev context, the second scale function $Z$ \cite{AKP,Pisexit,IP}   may also be defined  via the solution of the {\bf non-smooth total discounted "regulation"/capital injections problem}.

Let $X^{[d\cd}_t=X_t+\L_t$  denote the process $X_t$ modified by  Skorokhod reflection at $d\cd$,  and  let $\E^{[d\cd}_x$  denote expectation for this process and let
  $ T_b^{[d\cd}$ denote the first passage  to $b$ of $X^{[d\cd}_t$.

It may be checked   by  Ivanovs-Palmowski proof of Theorem \eqr{l:refl} (see Remark \ref{r:proof})
that this keeps  being true when \gen drawdown  reflection at $d(\cdot)$ replaces reflection  at $0$, i.e. that the relation \eqr{DbUZ} is still equivalent to
\begin{equation} \label{refbailoutdd} 
\sRui_{\q,\th,[d\cd}^b(x):=\E^{[d\cd}_x \left[e^{-\q T_b^{[d\cd} - \th \L_{T_b^{[d\cd}}}\right]
=\bc \dfrac{Z_{\q,[d\cd}(x,\th)}{ Z_{\q,[d\cd}(b,\th)} & \th <\I\\\E_x \left[e^{-\q \tb} \1_{ \{ \tb < \t_{d\cd} \}}\right]=
\dfrac{W_{\q,d\cd}(x)}{ W_{\q,d\cd}(b)} &  \th =\I\ec\end{equation}

\eeR

\ssec{Optimal dividends problem with \gen drawdowns\la{s:deFM}}
Let $T_{\H d(\cdot)} = \td_{\H d(\cdot)} \wedge \ta  $
  denote the first passage time either below $a$, or below the drawdown  boundary for the process $X^{b]}_t$ reflected at $b$ with regulator $U_t$. One can consider the extension of de Finetti's optimal dividend problem \eqref{div}
   \begin{equation} \la{tdiv} V^{b]}(x)=V_{\q,\H d(\cdot)}^{b]}(x):=  \E_x \left[ \int_{0}^{T_{\H d(\cdot)} } e^{-q  t}dU_t \right]\end{equation}where $V$ depends now also on the function $\H d(\cdot)$.\fn[4]
 {This definition assumes that the initial point satisfies $X_0=\ovl X_0=x$, i.e. that the starting point is on the $x$ axis in Figure \ref{f:plD}.} 

By the strong Markov property, it holds that
 \begin{equation}
 V^{b]}(x)= {\E_x \left[e^{-q \tb } ; \tb \leq \min(\td_{\H d(\cdot)}, \ta) \right]} v(b),\quad
   v (b)=V_{\q,\H d(\cdot)}^{b]}(b)=  \E_b \left[ \int_{0}^{T_{\H d(\cdot)} } e^{-q  t}d\U^{b]}_t \right]. \la{VLd}\end{equation}

\beR \la{r:expo} The function $v(b)$ represents the expected discounted time until killing for the reflected process, when starting from $b$. This equals the time the process reflected at $b$ spends at point $(b,0)$ in Figure \ref{f:plD}, before a downward excursion beyond $\H d(b)$ kills the process. Furthermore, this time  is exponential with parameter $\nu_{\q}(b,\H d(b))$ (as a consequence of the fact  that the drawdown  process away from a running maximum is Markovian and the corresponding process of upward excursions  is Poisson, just as in the \lev case).
 Thus, the expectation  is the reciprocal of   $\nu_{\q}(b,\H d(b))$,
  and \begin{equation} v (b)=\nu_{\q}(b,\H d(b))^{-1}=\fr{W_{\q, d(\cdot)}(b)}
  {W_{\q, d(\cdot)}'(b)} \la{VLv}\end{equation}

\eeR

\beR
By \eqr{UbDex}, \eqr{VLv}   we arrive finally   to an explicit formula
for $V^{b]}(x)$:
\beq \la{V} V ^{b]}(x)=
\fr{e^{-\int_x^b \nu_q(y, \H d(y)) \md s}}
{\nu_q(b, \H d(b))}
\eeq
  expressing the expected dividends in terms of
  $\nu_q(y, \H d(y))$. Note that in the \lev case the equation \eqr{V}   simplifies to:
 \bea
 V^{b]}(x)=\fr{W_{\q}(d(x))}{W_{\q}(d(b))} \nu_{\q}(d(b))^{-1}\eea
 (using $x-\H d(x)=d(x)$), which checks with \cite[Lem. 3.1-3.2]{WangZhou}.

 The problem of choosing a drawdown  boundary to optimize dividends in \eqr{V} is tackled in \cite{AG} via Pontryaghin's maximum principle. The result depends of course of the process considered, but it always must use one of two types of segments: ``de Finetti segments"  of maximal slope, of direction
 $(\H d'(s), d'(s)=(0,1)$ and segments along which the equation
 \be {\partial }_{2}\nu_q(s,\H d(s))=const \ee is satisfied.

 \eeR

{\bf  For spectrally negative L\'evy process and affine drawdowns} $d(x)=(1-\xi)x+d$, $\H d(x)= \xi x-d$, $h(x)={ d(x)}/{\xi}$, the exit functions and $v(b)$ in \eqref{VLv} are simpler:
\beq \la{DbUaff} && W_{q, d(\cdot)}(x)=\le({W_q(d(x))}\ri)^{\fr 1{1-\xi}}, \quad  UbD_{q}^{b}(x,\H d(\cdot))=\le(\fr{W_q(d(x))}{W_q(d(b))}\ri)^{\fr 1{1-\xi}}, \no \\&&\nu_q(x,\H d(x))=\fr{1}{W_{q, d(\cdot)}(x) }  \fr{\md W_{q, d(\cdot)}(x)}{ \md x} = \fr{W_q'(d(x))}{W_q(d(x))}, \quad v(b)=\fr{W_q(d(b))}{W_q'(d(b))}
\eeq
see   \cite[Thm.~1.1]{AVZ}, with tax parameter $\g=0$,  and  \cite[Rem.~7]{AVZ}, with tax parameter $\g=1$.

 We may obtain in this case a more precise version of  Proposition \ref{Thmp}. Note first that when $a+h(a) > b$, the drawdown  constraint is invisible. The value function \eqref{tdiv} is therefore $\fr{W_{\q}(x-a)}{W_{\q}'(b-a)}$ (which can be maximized by minimizing $b \mapsto W_{\q}'(b)$ -- see Sec.~\ref{s:deF}).

   When $a+h(a) \leq b  $, combining  the discounted probability of reaching $b$ and the value $v(b)$ yields:
  \beP \label{Thmo} Consider a spectrally negative L\'evy process $X$ with three times differentiable scale function $W_{\q}$. Assume $d(x):= (1-\xi) x +d,$ where
$d \geq 0$, $\xi \leq 1$, $a  \leq   x \leq b$, $a+h(a) \leq b$.    Then:

A)  the expected discounted dividends are:
\beq \la{div3}V ^{b]}(x)=\bc
 \left(\frac{W_{q}(d(x))}
{W_{q}(d(b))}\right)^
{\frac{1}{1-\xi}} \frac{W_{q }(d(b))}{W_{q }'(d(b))}, &  {a+h(a)  \leq x },\\\frac{{W_{\q}(x-a )}}{{W_{\q}\left((h(a)\right)}}
 \left(\frac{W_{q}\left( (h(a)\right)}
{W_{q}(d(b))}\right)^
{\frac{1}{(1-\xi)}}  \frac{W_{q }(d(b))}{W_{q }'(d(b))}, & x   \leq a+h(a).
\ec \eeq

B)
 The barrier influence function  (which must be optimized in $b$) in the case $a+h(a)  \leq x$ is
\beq BI(b,d,\xi)=\fr{W_{\q}(d(b))^{ {1-\fr{1}{1-\xi}}}}
{ W'_{\q}(d(b))}
=\fr{W_{\q}(d(b))^{ {-\fr{\xi}{1-\xi}}}}
{ W'_{\q}(d(b))}. \la{BI}
  \eeq
  The critical points $b^*$ for fixed $d,\xi$ satisfy \footnote{When $\xi=d=0$, we recover in the \cP case the equation $W_{\q}''(b)=0$.}
{\beq \fr{W_{\q}'' W_{\q}}{(W_{\q}')^2}(d(b^*))+ \fr{\xi}{1-\xi}=0. \la{b*critical}\eeq}

  For local maxima at $b^* >0$ to exist, it is necessary that $\fr{W_{\q}'' W_{\q}}{(W_{\q}')^2}(0)+ \fr{\xi}{1-\xi}<0$ and that
   $\Bigg(W_{\q} W_{\q}' W_{\q}''' + W_{\q}'' \Big(W_{\q}' \Big)^2  - 2 W_{\q} \Big(W_{\q}'' \Big)^2\Bigg)(d(b_*))>0$.

  C)
 The barrier influence function  in the case $ x   \leq a+h(a)$ is
\beq W_{\q}\left(h(a)\right)^{ \fr{\xi}{1-\xi}}
 BI(b,d,\xi).
  \eeq

\eeP

\prf A) The  first case,  in which barrier $a$ is invisible, holds by \cite[Thm.~1.1]{AVZ} (by plugging there $\g=0$).\fn[4]{Note that the limiting case $\xi=1$ is consistent by L'Hospital's theorem with our previous $UbD^b_q(x,d)$ defined in \eqref{UbDL}.}

The second case holds by the \str. Note  that until $\ovl X_t$ visits  $a+h(a)$,
 {the upper drawdown  barrier is invisible}, and the classic formula for smooth passage applies. Subsequently,
 we are in the first case, with starting point $x=a+h(a)$, applying the first case and using $d(a+ h(a))=h(a)$ (see Figure \ref{f:plD}).

B) For the critical points, note that the sign of $-BI'$ coincides with that of $\fr{W_{\q}'' W_{\q}}{(W_{\q}')^2}(d(b))+ \fr{\xi}{1-\xi}$, and that $W'$ is positive. \qed

\beR
To compare value functions when $\xi,d  $ vary, let us choose the fixed point  $x=a=0$.
It may be easily checked that for any $\xi=0,d \geq 0$ $ V^{b]}(0)= \fr{W_{q}(d)}{W_{q}'(b_0)}$, where $b_0$ is the argmax of $BI(b)$ when $\xi=0$ (using the translation invariance of \lev processes).

Also,  the ``de Finetti solution" $\xi=0$ always beats $\xi>0$ at equal $d$, due to the singularity  of $BI(b)$ \eqr{BI} at $0$ when $\xi>0$, which  makes immediate stopping optimal. 
\iffa
Since ${W_{q}(d)}$ is increasing, it follows
that without extra constraints,  with affine drawdown  boundary, the optimal solution is trivially $d=\I, \xi=0 \Eq b^*=1, V^{b^*]}(x)=\I$. Other solutions become thus of interest only under a constraint $d(a) \leq d_0$.

Furthermore, $\xi>0$    becomes interesting  once an upper  bound on the derivative $d'(s)$ or on the total ``regret/risk area'' is placed -- see Figure \ref{f:plD}.
\eeR

Let us provide an example.
\beXa {\bf Brownian motion}\label{ex:BM}
Consider Brownian motion
with drift  $X(t) = \s B_t + \mu  t$ and   affine drawdown stopping. The scale function $W_q$ is given in \eqref{BMWsf}.

Assume that $x\geq a+h(a)=a+ \frac{d(a)}{\xi} = \frac{a+d}{\xi}$ so that the barrier influence function is given by \eqref{BI}. By Theorem \ref{Thmo}, the critical point $b^*$ satisfies \eqref{b*critical} which by using \eqref{ha} reduces to
\begin{equation*}
   \fr{\xi}{1-\xi}\fr{\s^2}2  \left( \nu_q' (d (b^*)) \right)^2 - \mu \; \nu_q (d (b^*))  + q=0.
\end{equation*}
Solving the quadratic equation implies that $b^*$ satisfies
\beq \la{gon} \fr{\mu}{2 \q} + \sqrt{\left(\fr{\mu}{2 \q}\right)^2-\fr{\s^2 \xi}{2\q(1-\xi)}}\nu_\q(d (b^*))=1,\eeq
which reduces when $\xi=0$ to \eqref{go}.
\end{Exa}

\sec{Chronology \la{s:ch}} \BEN[A)] \im Ruin theory for the {Cram\'{e}r-Lundberg} or {compound Poisson}  risk model was born in Lundberg's treaty \cite{lundberg1903approximerad}.
\im The extension to the \lev case was achieved in the landmark paper ``Problem of destruction and resolvent of a terminating process with independent increments'', where  the formula  \bea \sRui_{\q}^{b]}(x,a)=\E_x \left[ e^{-\q \tb} 1_{\left\{ \tb <  \tz \right\}}\right]  =\frac { W_{\q}(x)}  { W_{\q}(b)}\eea
 for the "smooth"  two-sided exit problem (TSE) \cite[Thm.~3]{Suprun} is provided\fn[5]{Informally, $ W_{\q}$ may be viewed as an  analog of the transfer function for discrete systems.}. The  \LT of $ W_{\q}$ was computed in \cite[(33)]{Suprun}. Also,  \cite[Thm.~2]{Suprun} provided the formula of the resolvent density for the process killed outside an interval $[a,b]$\fn[4]{Under the  \CL risk model, \cite{dickson1992distribution}
derived independently the particular case $\q=0$ of the resolvent
formula -- see also  Gerber and Shiu \cite[(6.5-6.6)]{gerber1998time}, who extend  Dickson's resolvent formula  to $\q >0$.}.
$$u_\q(x,y) =  \frac {W_{\q}(x-a) } {W_{\q}(b-a)} W_{\q} (b-y) -W_{\q} (x-y).$$

\im
\cite[(4)-(7)]{Ber97} introduced the notation $ W_{\q}$ and the name scale function for spectrally negative \lev processes. The central object of the paper  is now  $ W_{\q}$ (instead of Suprun's resolvent).       Probabilistic proofs of other problems are provided, by reducing them to smooth TSE. The non-smooth
  two-sided first passage problem is solved in
  \cite[Cor.~1]{Ber97},   and \cite[Thm.~2]{Ber97}
  determined the  decay parameter $\r$ of the process killed upon exiting an interval, and showed that the quasi-stationary distribution is $W_{-\r}$. The subsequent landmark textbook   \cite{Ber} offers a comprehensive  treatment of \lev processes, including the beautiful excursion theory.

  \im A first treatment of the optimal discounted dividends problem in the
classical compound Poisson model can be found in Section 6.4 of Buhlmann (1970) \cite{buhlmann2007mathematical}. The resulting formula $\frac {W_{\q}(b) } {W_{\q}'(b)}$ for dividends at $b$, when starting from $b$, is a consequence of the fact that  the discounted dividends have an exponential law of rate $\frac  {W_{\q}'(b)}{W_{\q}(b) }$.

  \im \cite{lin2003classical} studies the \GS function (a generalization of the ruin probability) for a \cP process with a constant barrier and discovers the ``dividends-penalty'' identity connecting  it to the scale function, denoted by $h$, and to the \GS function without barrier.
  \im
  \cite{AKP} introduced the second scale function  $ Z_{\q}$,
  initially for  relating to $ W_{\q}$ the solution of the ruin problem
   $ 
   \Rui_\q(x):=\E_x \le[e^{-\q \tz } \; \1_{\{\tz < \I\}} \ri]=  Z_{\q}(x) -  W_{\q}(x) \,  \fr {\q}{\Fq}.$
A case could be made for using $\Rui_\q(x)$ rather than $ Z_{\q}(x)$  as the second "alphabet letter" in first passage formulas. In fact, the former, being bounded,  is more convenient to compute numerically. However, it turned out that $ Z_{\q}(x)$ leads often to simpler  results and proofs, due to the fact that $e^{-\q t}  Z_{\q} (X_t)$ is a martingale \cite[Rem 5]{AKP}, \cite{nguyen2005some}. \im \cite{PDR,Pisexit} solved  in terms of $W,Z$ several first passage problems  for reflected processes.
  \im \cite{zhou2007exit} remarks that previous excursion theory proofs can often be replaced by simple applications of the strong Markov property, and of "$\e$ approximation" arguments in the non \cP case.

  \im \cite{Kyp} provided a comprehensive textbook on \lev processes  and applications.
   \im \cite{KL}  solved the TSE  for refracted processes (which are skip-free, but not L\'evy), in terms  of extensions of $W$ and $Z$.
    \im \cite{APP15,IP} introduced
    the two variables  extension $Z_{\q}(x,\th)$, which is useful for example for computing the \GS function
    $  \Rui_{\q,\th}^{b}(x):=\E_x\left[e^{- \q \tz+ \th X_{\tz}} \1_{\{\tz<\tb\}} \right]=Z_{\q}(x,\th) -   \fr { W_{\q}(x)}{  W_{\q}(b)}Z_{\q}(b,\th)
$ see Theorem \ref{l:s} A). The first paper showed also that this function was the unique ``smooth'' $\q$-harmonic extension of $e^{x \th}, x \leq 0$.

   \im \cite{Iva,IP} showed that  the known formulas on spectrally negative \lev processes  apply for spectrally negative Markov additive processes.
\im \cite{AIZ,BPPR,li2015two,AZ} ibidem for exponential Parisian processes.

\im \cite{LP,LZ17,vidmar2018first} ibidem for Omega  models
(processes with state dependent killing).
\im \cite{AV} ibidem for skip-free discrete state-space random walks.
 \im \cite{vidmar2018temporal,vidmar2018exit} ibidem for     positive \ssr \ Markov processes with one-sided jumps.
 \im \cite{APY} study  exponential Parisian processes with non-smooth reflection restricted to a buffer (this boundary regime interpolates between reflecting and stopping).
\im \cite{LLZ17b,ALL,AG,avram2019w} initiate the study of  time-
homogeneous strong Markov processes with one-sided jumps.
\EEN

\sec{List of notations \la{s:n}}
\renewcommand\arraystretch{1.65}
\begin{tabular}{p{0.54\textwidth}|p{0.46\textwidth}}
\hline
$\ta,\ \tb,\ \td_{d(\cdot)},\ \und \td_{d(\cdot)},\ \ta^{b]},\ \tb^{[a}\  $& times of first passage \eqr{fpt},  drawdown ,  draw-up \eqr{dd}, first passage with reflection \eqr{Trb}, \eqr{Tra}\\
\hline
$ \sRui_{q}^{b}(x,a)  =\mathbb{E}_{x}\left[  e^{-q\tb } \1_{\left\{
\tb< \ta \right\}  }\right] = \dfrac{W_{\q}(x-a)}{W_{\q}(b-a)}$& survival probability \eqr{sRui}, \eqr{twos}\\\hline
$\Rui_{q,\th}^{b}(x,a)    =\mathbb{E}_{x}\left[  e^{-q\ta+\th(X_{\ta}-a
)} \1_{\left\{  \ta<\tb\right\}  }\right]$& ruin probability \eqr{Rui}, \eqr{RZ}\\
\hline
$\und X_t=\inf_{0 \leq s\leq t} X_s, \quad \ovl X_t = \sup_{0 \leq s\leq t} X_s,  $& infimum and supremum processes \eqr{Xd}
  \\\hline
  $\L^{[a }_t=-(\und X_t-a)_-, \quad   \U_t=\U^{b]}_t=\le(\ovl X_t -b\ri)_+$
  & minimal  ``Skorokhod regulators" \eqr{Xd}
  \\\hline
   $X_t^{[a}=X_t  + \; \L_t, \quad  X^{b]}_t=X_t   -  \U_t$&regulated processes \eqr{Xd}
   \\\hline
  $  Y_t=\ovl{ X}_t-X_t, \quad \H Y_t=X_t -\und X_t$&  drawdown   and draw-up processes \eqr{ddp}, \eqr{dup}\\\hline
$\k(\th),\ \F(\q), \ W_\q(x),\  Z_{\q}(x,\th),\ W_{\q,\r}(x),\ Z_{\q,\r}(x,\th)$& Levy exp. \eqr{LE}, its inverse \eqr{Fq},  sc. functions \eqr{WLT}, \eqr{Zt}, \eqr{Z}, Parisian sc.functions \eqr{Z2}\\\hline
 $\nu_\q(s)=\fr{ W_{\q}'(s_+)}{ W_{\q}(s)}$ &rate  of down excursions larger than $s$ \eqr{nu}\\\hline
 $u_\q(x),\ u_\q^{|a}(x,y),\  u_\q^{|a,b]}(x,y), \ u_\q^{[a,b]}(x,y),\ u_\q^{[a,b|}(x,y)$&resolvents of free and constrained processes
  \\\hline
 $Z_{\q}^{(1)}(x)=\fr{\partial Z_{\q}(x,\th)}{\partial \th}_{\th=0}=\ovl   Z_{\q}(x) -\k'(0_+) \ovl  W_{\q}(x)$&\GS function for $w(x)=x$ \eqr{Z1}\\\hline
  $\sRui_{\q,\th}^b(x,a])=\Ea_x \left[e^{-\q \tba  - \th \L_{\tba }}\right]
= \dfrac{Z_{\q}(x-a,\th)}{Z_{\q}(b-a,\th)}$& discounted cumulative bailouts \eqr{B}, \eqr{refbailout}
\\\hline
$V^{b]}(x) = \Eb_x\le[\int_0^{\tzb} e^{-\q  t}d  \U_t\ri]=\frac{ W_{\q}(  x)}{W_{\q}^ {\prime}(b)}, $ & expected discounted dividends until $\tzb$ \eqr{div}\\\hline $V^{[0,b]}(x) = \Ezb_x\le[\int_0^{\I} e^{-\q  t}d  \U_t\ri]=\frac{Z_{\q}(x)}{  Z'_{\q}(b)}$ & expected discounted dividends with double reflection \eqr{divSLG}\\\hline $V_{w}^{b]}(x)
={G}_w(x) +   W_\q(x)\fr {1-{G}_w'(b)}{ W'_\q(b)}  $&modified de Finetti objective \eqr{deFmod}\\\hline
${\de}_{\q,\th}(x,d,s)= {\E_x \left[
e^{ -\q \td_d - \th (Y_{\td_d}-d)}; \ovl X_{\td_d} \in \md s
\right]}
$& joint law of maximum and drawdown
at  drawdown time \eqr{funL}\\\hline
$\T \de_{q,\th}(d)=\E_x \left[
e^{ -\q \td_d - \th  (Y_{\td_d}-d)}\right]=Z_{\q}(d,\th)-W_\q(d)\fr{Z_{\q}'(d,\th)}{W_\q'(d)}
, \for x$&drawdown function \eqr{de}\\\hline $DP_{\q,\th,\vt}^{b]}(x):=\Eb_x\left[e^{ -\q T_0^{b]}  + \th X_{T_0^{b]}}-\vt  \U_{T_0^{b]}}} \right], DP_{\q,\theta,\vt}^{\vdots 0,b]}(x)$& dividends-penalty functions \eqr{DP},\eqr{DP2}\\\hline$DB_{\q,\th,\vt}^{[0,b]}(x)=\Ezb_x\left[e^{  -\vt \U_{e_\q}- \th \L_{e_\q} }\right]$& dividends-bailouts function \eqr{DB} \\
\hline $
 B^{[0,b|}(x)=\E^{[0,b|}_x\le[\int_0^{\tbz} e^{-\q t}d \L_t\ri]=\frac{Z_{\q}( x)}{Z_{\q}(b)} G_{\q}^B(b)-G_{\q}^B(x), G_{\q}^B(x)= \ovl{Z}_\q(x)+\frac{\k'(0_+)}{\q}, G_{\q,\r}^B(x)=
\frac{\r}{\q+\r}G_{\q}^B (x)$&expected  (Parisian) discounted bailouts   until $\tbz$ \eqr{VL}, \eqr{G},\eqr{ParBail},\eqr{ParBailG}\\\hline $
 B^{[0,b]}(x)=\Ezb_x\le[\int_0^{\I} e^{-\q t}d \L_t\ri]=\frac{Z_{\q}( x)}{Z_{\q}'( b)} (G_{\q}^B)'(b)-G_{\q}^B(x)$&expected  discounted bailouts   with double reflection \eqr{VLSLG}\\\hline $V_{S,k}^{[0,b]}(x)=V^{[0,b]}(x)- k B^{[0,b]}(x)$&Shreve-Lehoczky-Gaver objective \eqr{value-SLG}\\\hline
$UbD_{q,\H d (\cdot)}^{b}(x) =\E_x \left[e^{-q \tb } ; \tb \leq  \td_{\H d (\cdot)} \ri] $&up before drawdown \eqr{UbDdef}\\\hline
$DbU^{b}_{q,\th,\H d (\cdot)}(x)=\E_x \left[ e^{-q\td_{\H d (\cdot)}-\th \T Y_{\td_{\H d (\cdot)}}} ; \td_{\H d (\cdot)}<\tb \right]$& drawdown  before up \eqr{DbUdef}
\\\hline
\end{tabular}

  \ssec{A summary of asymptotic relations for spectrally negative \lev \procs\la{s:asy}}
\BEN \im
When $\k'(0_+) >0$, $\Fq$ is the asymptotically dominant singularity  of $ W_{\q}(x) \sim \fr{e^{x \Fq}}{\k'(\Fq)}= \F'(\q) e^{x \Fq}$ as $x\to \infty$. Furthermore,
by \eqr{Wu} $  W_{\q}(x)= \F'(\q) e^{\Fq x} -u_\q(-x).$

\im Recalling $Z_{\q}(x,\th)=  \le(\k(\th)-\q\ri)\int_0^\I e^{- \th y}  W_{\q}(x+y) dy$ \eqref{Zt}, \itf \beq \la{ZtoW}
 \lim_{x\to\infty}\frac{Z_{\q}(x,\th)}{W_{\q}(x)}=
\le(\k(\th)-\q\ri) \lim_{x\to\infty}\int_0^\I e^{- \th y} \fr{ W_{\q}(x+y)}{ W_{\q}(x)} dy= \frac{\k(\th)-\q} {\th-\Fq}.\eeq
When $\th =0$, this yields
 \beq \la{ZoW}
 \lim_{x\to\infty}\frac{Z_{\q}(x)}{W_{\q}(x)}=
 \frac{\q} {\Fq}, \qu
 \lim_{x\to\infty}\frac{Z_{\q}(x)}{Z_{\q}(x,\th)}=
 \frac {\Fq-\th}{\q- \k(\th)}\frac {\q}{\Fq}.\eeq

\im Recalling   $Z_{\q}(x,\th)=\fr{\k(\th)-\q}{\th-\Fq}  W_{\q}(x)+ \Rui_{\q,\th}(x)$ \eqref{Zt}, \itf
 \begin{equation}
 \la{ZW}
\lim_{\th \to \I} Z_{\q}(x,\th) \fr{\th-\Fq}{\k(\th)-\q} = W_{\q}(x)\ee
and
\be \la{ZWF}
 \lim_{\l \to \I} Z_{\q}(x, \Fqr ) \fr{ \Fqr }{\l}=\lim_{\l \to \I} Z_{\q}(x, \Fqr ) \fr{ \Fqr -\Fq}{\l} =W_{\q}(x)
\end{equation}

\EEN

%
%

%

\bs \bs

{\bf Acknowledgement.} Many thanks to  Hansjoerg Albrecher,  Ester Frostig, Jevgenijs Ivanovs, Bin Li, Ronnie Loeffen, Zbigniew Palmovski, Jos\'e-Luis Perez, Martijn Pistorius, Matija Vidmar and Xiaowen Zhou for  useful
discussions, and for their invaluable contributions to this field. D.~Grahovac acknowledges the support of University of Osijek grant ZUP2018-31.

\small
\bibliographystyle{amsalpha}
\newcommand{\etalchar}[1]{$^{#1}$}

\end{document}